\documentclass{amsart}

\usepackage[all]{xy}

\theoremstyle{plain}
\newtheorem{thm}{Theorem}[section]

\newtheorem{question}[thm]{Question}

\newtheorem{subthm}{Theorem}[subsection]
\newtheorem{subprop}[subthm]{Proposition}
\newtheorem{sublem}[subthm]{Lemma}
\newtheorem{subcor}[subthm]{Corollary}

\theoremstyle{definition}
\newtheorem{dfn}[thm]{Definition}
\newtheorem{subdfn}[subthm]{Definition}

\theoremstyle{remark}
\newtheorem{rem}[thm]{Remark}

\newtheorem{subrem}[subthm]{Remark}
\newtheorem{subfacts}[subthm]{Facts}

\newcommand{\HH}{\mathrm{H}}

 \newcommand{\unter}[2]{\genfrac{}{}{0pt}{}{#1}{#2}}

\begin{document}

\title{Universal deformation rings and dihedral defect groups}

\author{Frauke M. Bleher}
\address{Department of Mathematics\\University of Iowa\\
Iowa City, IA 52242-1419}
\email{fbleher@math.uiowa.edu}
\thanks{The author was supported in part by  
NSF Grant DMS01-39737 and NSA Grant
H98230-06-1-0021.}

\begin{abstract}
Let $k$ be an algebraically closed field of characteristic $2$, and let $W$ be the ring of infinite Witt vectors over $k$. Suppose $G$ is a finite group, and  $B$ is a block of $kG$ with dihedral defect group $D$ 
which is Morita equivalent to the principal $2$-modular block of a finite simple group.
We determine the universal deformation ring $R(G,V)$ for every $kG$-module $V$ which belongs to $B$ and has stable endomorphism ring $k$. It follows that $R(G,V)$ is always isomorphic to a subquotient ring of $WD$. Moreover, we obtain an infinite series of examples of universal deformation rings which are not complete intersections.
\end{abstract}

\maketitle


\section{Introduction}
\label{s:intro}
\setcounter{equation}{0}
\setcounter{figure}{0}

In this paper we determine the universal deformation rings $R(G,V)$ associated to 
certain mod $2$ representations $V$ of finite groups $G$ which belong to blocks  of $G$ having dihedral defect group $D$. There are three reasons for making this calculation. The first is that it provides more evidence for an affirmative answer to a question posed in \cite{bc}, which would bound $R(G,V)$ in terms of the group ring of $D$ over the Witt vectors (see Question \ref{qu:main}). The second reason is that we produce an infinite family of $R(G,V)$ which are not complete intersections, but which do satisfy the dimension versus depth condition conjectured in \cite{bc5}  (see Question \ref{qu:lci}). This in turn leads to interesting questions in number theory (see the discussion after Question \ref{qu:lci}). The final reason is to describe a general method for applying results in modular and ordinary representation theory due to Brauer, Erdmann and others to the computation of universal deformation rings.

To make this more precise, let $k$ be an algebraically closed field of characteristic $p>0$, let $W=W(k)$ be the ring of infinite Witt vectors over $k$, and let $F$ be the fraction field of $W$. Let $\Gamma$ be a profinite group, and suppose $V$ is a finite dimensional vector space over $k$ with a continuous $\Gamma$-action. If all continuous $k\Gamma$-module endomorphisms of $V$ are given by scalar multiplications, an argument of Faltings (see \cite{lendesmit}) shows that $V$ has a universal deformation ring $R(\Gamma,V)$. The topological ring $R(\Gamma,V)$ is universal with respect to deformations of $V$ over commutative local $W$-algebras with residue field $k$ which are the projective limits of their discrete Artinian quotients. For more information on deformation rings see \cite{lendesmit}, \cite{maz1} and \S\ref{s:prelimdef}. In number theory, deformation rings are at the center of work by many authors concerning Galois representations, modular forms, elliptic curves and diophantine geometry (see e.g. \cite{cornell}, \cite{wiles,taywiles}, \cite{breuil} and their references).

In \cite{lendesmit}, de Smit and Lenstra show that $R(\Gamma,V)$ is the inverse limit of the universal deformation rings $R(G,V)$ when $G$ runs over all finite discrete quotients of $\Gamma$ through which the $\Gamma$-action on $V$ factors. Thus to answer questions about the ring structure of $R(\Gamma,V)$, it is natural to first consider the case when $\Gamma=G$ is finite. In this case, $V$ has a universal deformation ring $R(G,V)$ under the weaker condition that the stable endomorphism ring $\underline{\mathrm{End}}_{kG}(V)$ is of dimension $1$ over $k$ (see \cite[Prop. 2.1]{bc}); we assume $\underline{\mathrm{End}}_{kG}(V)=k$ in what follows. 

In \cite{bc}, the author and T. Chinburg determined $R(G,V)$ for $V$ belonging to blocks with cyclic defect groups, i.e. blocks of finite representation type. 
In \cite{bl}, the author considered $V$ belonging to blocks with Klein four defect groups and
described their universal deformation rings.
This is a natural progression from \cite{bc}, since the blocks with Klein four defect groups have tame representation type and are the only such blocks having abelian (non-cyclic) defect groups. 
The results obtained in \cite{bc} led to the following question.

\begin{question}
\label{qu:main} 
Let $B$ be a block of $kG$ with defect group $D$, and suppose $V$ is a finitely generated $kG$-module with stable endomorphism ring $k$ such that the unique $($up to isomorphism$)$ non-projective indecomposable summand of $V$ belongs to $B$. Is the universal deformation ring $R(G,V)$ of $V$ isomorphic to a subquotient ring of the group ring $WD$?
\end{question}

It is shown in \cite{bc} and \cite{bl} that the answer to this question is positive in case $D$ is cyclic or a Klein four group. Moreover, it follows that in all these cases $R(G,V)$ is a complete intersection ring (see \cite[Thm. 7.2]{bc5}). 
In case $k$ has characteristic $2$, $G$ is the symmetric group $S_4$ and $E$ is the unique (up to isomorphism) non-trivial simple $kS_4$-module, the author and T. Chinburg showed in \cite{bc4.9,bc5} that $R(S_4,E)\cong W[t]/(t^2,2t)$. Hence $R(S_4,E)$ is not a complete intersection ring, which answers a question of M. Flach \cite{flach}. Note that $E$ belongs to the principal block of $kS_4$ which has as defect groups dihedral groups of order $8$. A new proof of this result has been given by Byszewski in \cite{cantrememberhisname} using only elementary obstruction calculus. In \cite{bc4.9,bc5}, the author and T. Chinburg also considered the question of whether deformation rings which are not complete intersections arise from arithmetic. It was shown in particular that there are infinitely many real quadratic fields $L$ such that the Galois group $G_{L,\emptyset}$ of the maximal totally unramified extension of $L$ surjects onto $S_4$ and $R(G_{L,\emptyset},E)\cong R(S_4,E)\cong  W[t]/(t^2,2t)$ is not a complete intersection, where $E$ is viewed as a module of $G_{L,\emptyset}$ via inflation. 

In this paper, we expand these results further by considering entire families of blocks of tame representation type with defect $d\ge 3$. More precisely, our goal is to determine the structure of the universal deformation rings $R(G,V)$ in case $D$ is dihedral of order at least $8$ and $B$ is Morita equivalent to the principal $2$-modular block of a finite simple group. In particular, $B$ contains precisely three isomorphism classes of simple modules. Note that in \cite{brauer2}, Brauer has proved that a block with dihedral defect groups contains at most three simple modules up to isomorphism; hence we look at the largest case. It follows by the classifications by Gorenstein-Walter \cite{gowa} and by Erdmann \cite{erd} that except for possibly one remaining family, all blocks $B$ with dihedral defect groups containing precisely three isomorphism classes of simple modules are Morita equivalent to the principal $2$-modular block of some finite simple group (see Remark \ref{rem:oherdmann}). 

A summary of our main results is as follows. The precise statements can be found in Propositions \ref{prop:psl1}, \ref{prop:psl2} and \ref{prop:a7} and in Theorem \ref{thm:main} and Corollary \ref{cor:lci}.
\begin{thm}
\label{thm:bigmain}
Suppose $k$ has characteristic $2$.
Let $B$ be a block of $kG$ with dihedral defect group $D$ of order $2^d$ where $d\ge 3$, which is Morita equivalent to the principal $2$-modular block of a finite simple group. Let $V$ be a finitely generated $B$-module with stable endomorphism ring $k$ and universal deformation ring $R(G,V)$. Then either
\begin{enumerate}
\item[i.] $R(G,V)/2R(G,V)\cong k$, in which case $R(G,V)$ is isomorphic to a quotient ring of $W$, or
\item[ii.] $R(G,V)/2R(G,V)\cong k[t]/(t^{2^{d-2}})$, in which case  $R(G,V)\cong W[[t]]/(p_d(t)(t-2),2\,p_d(t))$ for a certain monic polynomial $p_d(t)\in W[t]$ of degree $2^{d-2}-1$ whose non-leading coefficients are all divisible by $2$.
\end{enumerate}
In all cases, $R(G,V)$ is isomorphic to a subquotient ring of $WD$. Given the block $B$, each of the cases $(i)$ and $(ii)$ occurs for infinitely many $V$. In case $(ii)$, $R(G,V)$ is not a complete intersection.
\end{thm}

This Theorem  also gives an affirmative answer to the following question from \cite{bc5} for $B$, $D$ and $V$ as in the statement of the Theorem  (see \cite[Thm. 7.2]{bc5}).

\begin{question}
\label{qu:lci}
Suppose $B$, $D$ and $V$ are as in Question $\ref{qu:main}$ and $D$ has nilpotency $r$. Is it the case that $\mathrm{dim}(R(G,V)) - \mathrm{depth}(R(G,V)) \le r-1$?
\end{question}

Theorem \ref{thm:bigmain} provides an infinite series of  finite  groups $G$ and mod $2$ representations $V$ for which $R(G,V)$ is not a complete intersection (see Corollary \ref{cor:lci}). This raises the question of whether one can use such $G$ and $V$ to construct further examples of deformation rings arising from arithmetic which are not complete intersections, in the following sense.  As in \cite{bc5}, one can ask whether there are number fields $L$ together with a finite set of places $S$ of $L$ such that $G$ is a quotient of  the Galois group $G_{L,S}$ of the maximal unramified outside $S$ extension $L_S$ of $L$ which has the following property. There should be a surjection $\psi:G_{L,S}\to G$ which induces an isomorphism $R(G_{L,S},V)\to R(G, V)$ of deformation rings when $V$ is viewed as a representation for $G_{L,S}$ via $\psi$. It was shown in \cite{bc5} that a sufficient condition for $R(G_{L,S},V)\to R(G, V)$ to be an isomorphism is that $\mathrm{Ker}(\psi)$ has no non-trivial pro-$2$ quotient. This is equivalent to the requirement that if $L'$ is the fixed field of $\mathrm{Ker}(\psi)$ acting on $L_S$, then each ray class group of $L'$ associated to a conductor involving only places over $S$ should have odd order. 

As mentioned earlier, the above arithmetic problem was considered in \cite{bc5} when $G=S_4$ and $V$ is irreducible of dimension $2$. 
It is more challenging to treat the cases produced by Theorem \ref{thm:bigmain} in an analogous way, but we think this raises interesting questions in Galois theory. For example, if $G=A_7$ and $V$ is an irreducible mod $2$ representation of degree $14$, can one find $L$ and $S$ and $\psi:G_{L,S}\to G$ as above for which $\mathrm{Ker}(\psi)$ has no non-trivial pro-$2$ quotient, and hence $R(G_{L,S},V)\cong R(G, V)$ is not a complete intersection? Is this possible when $L=\mathbb{Q}$? Another interesting case provided by Theorem \ref{thm:bigmain} is when $G$ is isomorphic to $\mathrm{PSL}_2(\mathbb{F}_q)$ where $q$ is an odd prime power and $8$ divides $\# G$. For example, if $q=\ell^2$ where $\ell\not\equiv \pm 1\mod 24$ (resp. $\ell\not\equiv 1,4,16\mod 21$), it was proved in \cite{takehito} (resp. in \cite{detwew}) that $\mathrm{PSL}_2(\mathbb{F}_q)$ occurs regularly over $\mathbb{Q}(t)$, implying that there are $\mathrm{PSL}_2(\mathbb{F}_q)$ extensions of any number field $L$ for such $q$.
On the other hand, it was shown in \cite{wiese} that for any prime $\ell$ there are infinitely many positive integers $r$ such that for $q=\ell^r$, $\mathrm{PSL}_2(\mathbb{F}_q)$ occurs as a Galois group over $\mathbb{Q}$.
We are looking for $\mathrm{PSL}_2(\mathbb{F}_q)$ extensions of number fields which satisfy additional constraints on their ray class groups. Generalizing the techniques of \cite{bc5} to treat such questions is beyond the scope of this paper, but we believe this will lead to interesting new number theoretic results.

We now describe the steps used to prove Theorem \ref{thm:bigmain}.

We first determine which indecomposable $B$-modules have stable endomorphism ring $k$, using that $B$ is Morita equivalent to a special biserial algebra. This enables us to use the description of indecomposable modules of special biserial algebras as so-called string and band modules. We go through each component $\mathfrak{C}$ of the stable Auslander-Reiten quiver of $B$ starting with modules in $\mathfrak{C}$ of minimal length. As it turns out, for each block $B$ there are infinitely many isomorphism classes of indecomposable $B$-modules with stable endomorphism ring $k$. More precisely, there are always two components $\mathfrak{C}$ of type $\mathbb{Z}A_\infty^\infty$ which consist entirely of modules with stable endomorphism ring $k$. The other modules with stable endomorphism ring $k$ either lie at the ends of $3$-tubes, or they form a single $\Omega$-orbit in one or two components $\mathfrak{C}'$ of type $\mathbb{Z}A_\infty^\infty$. 

After we have found all indecomposable $B$-modules $V$ with stable endomorphism ring $k$, we then determine their universal deformation rings modulo $2$, i.e. $R(G,V)/2R(G,V)$. To do so, we concentrate first on one  block $B$ of a given defect $d\ge 3$. We then make use of the fact that all the blocks $B$ of defect $d$ are stably equivalent by a stable equivalence of Morita type over $k$. This leads to the universal deformation rings $R(G,V)$ in case (i) of Theorem \ref{thm:bigmain}. 

For $V$ as in case (ii), our strategy is to use Brauer's results on the ordinary characters belonging to $B$ to find the largest quotient of $R(G,V)$ which is flat over $W$, namely $W[[t]]/(p_d(t))$ for $p_d(t)$ as in part (ii) of Theorem \ref{thm:bigmain}. We then use ring theory to show $R(G,V)$ must have the form $W[[t]]/(p_d(t)\,(t-2\gamma),\alpha\,2^m\, p_d(t))$ for some $\gamma\in W$, $\alpha\in\{0,1\}$ and $m\ge 1$. To determine $\gamma$, $\alpha$ and $m$, we take advantage of the fact that if $\overline{U}$ is the universal mod $2$ deformation of $V$ then $\underline{\mathrm{End}}_{kG}(\overline{U})=k$, so that $R(G,\overline{U})$ is well defined.  To compute $R(G,\overline{U})$, we use that $\overline{U}$ lies at the end of a $3$-tube of the stable Auslander-Reiten quiver of $B$. Using results from Brauer and Erdmann we see that the vertices of $\overline{U}$ are Klein four groups. Moreover, suppose $K$ is one such Klein four group, and let $N_G(K)$ be the normalizer of $K$ in $G$. Then the Green correspondent $f\overline{U}$ is a $kN_G(K)$-module which is induced from a module belonging  to the end of a $3$-tube of the stable Auslander-Reiten quiver of a block $b_1$ which is Morita equivalent to $kA_4$. The results of \cite{bl} imply that $R(G,\overline{U})=k$, which leads to $\alpha=1$ and $m=1$.

The paper is organized as follows.
In \S\ref{s:prelimdef}, we recall the definitions of deformations and deformation rings  
and prove that stable equivalences of Morita type preserve deformation rings (see Lemmas \ref{lem:stabmordef} and \ref{lem:weakstabmordef}). 
We also prove some results which help determine universal deformation rings that are certain quotient rings of $W[[t]]$ (see Lemmas \ref{lem:trythis}, \ref{lem:Wlift}, \ref{lem:herewegoagain} and \ref{lem:dihedraltrick}).
In \S\ref{s:dihedralsylow}, we use the classifications by Gorenstein-Walter \cite{gowa} and by Erdmann \cite{erd} to describe all $2$-modular blocks $B$ of a finite group $G$ which have dihedral defect groups and are Morita equivalent to the principal $2$-modular block of a finite simple group. We also describe results by Brauer \cite{brauer2} about the ordinary irreducible characters of $G$ belonging to $B$. 
In \S\ref {s:defmod2} and \S\ref{s:stableend}, we determine which indecomposable $B$-modules $V$ have stable endomorphism ring $k$ and find their universal deformation rings modulo $2$ (see Propositions \ref{prop:psl1}, \ref{prop:psl2} and \ref{prop:a7}). 
In \S\ref{s:udr}, we  
analyze the $B$-modules belonging to the boundaries of $3$-tubes and
use the results from \cite{brauer2} about the ordinary irreducible characters belonging to $B$ to determine the universal deformation rings of $V$ (see Theorem \ref{thm:main}). In particular, this implies Theorem \ref{thm:bigmain}. 
Since we make use of the fact that all the blocks we consider in this paper are Morita equivalent to special biserial algebras, we recall in \S\ref{s:prelimstring} the basic definitions
 of special biserial algebras and string algebras and describe their indecomposable modules and their Auslander-Reiten quivers.

The author would like to thank T. Chinburg for many valuable discussions which led to the elegant proof of the isomorphism class of the full deformation ring $R(G,V)$ in case (ii) of Theorem \ref{thm:bigmain} using ordinary characters. The author would also like to thank the University of Iowa for awarding her a Faculty Scholarship which gave her the opportunity to devote a full semester to research at the University of Pennsylvania in the Fall of 2005 during which time the main parts of this paper were written. Finally, the author would like to thank the Department of Mathematics at the University of Pennsylvania for their hospitality during this time.


\section{Preliminaries: Universal deformation rings}
\label{s:prelimdef}
\setcounter{equation}{0}
\setcounter{figure}{0}

In this section, we give a brief introduction to versal and universal deformation rings and deformations. For more background material, we refer the reader to \cite{maz1} and \cite{lendesmit}.

Let $k$ be an algebraically closed field of characteristic $p>0$, and let $W$ be the ring of infinite Witt vectors over $k$. Let $\hat{\mathcal{C}}$ be the category of all complete local commutative Noetherian rings with residue field $k$. The morphisms in $\hat{\mathcal{C}}$ are continuous $W$-algebra homomorphisms which induce the identity map on $k$. Let $\mathcal{C}$ be the full subcategory of $\hat{\mathcal{C}}$ of Artinian objects.


\subsection{Universal and versal deformation rings}
\label{ss:udr}
\setcounter{equation}{0}
\setcounter{figure}{0}

Suppose $G$ is a finite group and $V$ is a finitely generated $kG$-module. 
A lift of $V$ over an object $R$ in $\hat{\mathcal{C}}$ is a pair $(M,\phi)$ where $M$ is a finitely generated $RG$-module which is free over $R$, and $\phi:k\otimes_R M\to V$ is an isomorphism of $kG$-modules. Two lifts $(M,\phi)$ and $(M',\phi')$ of $V$ over $R$ are isomorphic if there is an isomorphism $\alpha:M\to M'$ with $\phi=\phi'\circ (k\otimes\alpha)$. The isomorphism class $[M,\phi]$ of a lift $(M,\phi)$ of $V$ over $R$ is called a deformation of $V$ over $R$, and the set of such deformations is denoted by $\mathrm{Def}_G(V,R)$. The deformation functor
$$\hat{F}_V:\hat{\mathcal{C}} \to \mathrm{Sets}$$ 
sends an object $R$ in $\hat{\mathcal{C}}$ to $\mathrm{Def}_G(V,R)$ and a morphism $f:R\to R'$ in $\hat{\mathcal{C}}$ to the map $\mathrm{Def}_G(V,R) \to
\mathrm{Def}_G(V,R')$ defined by $[M,\phi]\mapsto [R'\otimes_{R,f} M,\phi']$, where  $\phi'=\phi$ after identifying $k\otimes_{R'}(R'\otimes_{R,f} M)$ with $k\otimes_R M$.

In case there exists an object $R(G,V)$ in $\hat{\mathcal{C}}$ and a deformation $[U(G,V),\phi_U]$ of $V$ over $R(G,V)$ such that for each $R$ in $\hat{\mathcal{C}}$ and for each lift $(M,\phi)$ of $V$ over $R$ there is a unique morphism $\alpha:R(G,V)\to R$ in $\hat{\mathcal{C}}$ such that $\hat{F}_V(\alpha)([U(G,V),\phi_U])=[M,\phi]$, then $R(G,V)$ is called the universal deformation ring of $V$ and $[U(G,V),\phi_U]$ is called the universal deformation of $V$. In other words, $R(G,V)$ represents
the functor $\hat{F}_V$ in the sense that $\hat{F}_V$ is naturally isomorphic to $\mathrm{Hom}_{\hat{\mathcal{C}}}(R(G,V),-)$. In case the morphism $\alpha:R(G,V)\to R$ relative to the lift $(M,\phi)$ of $V$ over $R$ is not unique, $R(G,V)$ is called the versal deformation ring of $V$ and $[U(G,V),\phi_U]$ is called the versal deformation of $V$.

By \cite{maz1}, every finitely generated $kG$-module $V$ has a versal deformation ring.
By a result of Faltings (see \cite[Prop. 7.1]{lendesmit}), $V$ has a universal deformation ring in case $\mathrm{End}_{kG}(V)=k$.

\medskip

The following two results were proved in \cite{bc}, where $\Omega$ denotes the Heller operator for $kG$ (see for example \cite[\S 20]{alp}).

\begin{subprop}
\label{prop:stablend}
{\rm (\cite[Prop. 2.1]{bc}).}
Suppose $V$ is a finitely generated $kG$-module with stable endomorphism ring $\underline{\mathrm{End}}_{kG}(V)=k$.  Then $V$ has  a universal deformation ring $R(G,V)$.
\end{subprop}

\begin{sublem} 
\label{lem:defhelp}
{\rm (\cite[Cors. 2.5 and 2.8]{bc}).}
Let $V$ be a finitely generated $kG$-module with stable endomorphism ring $\underline{\mathrm{End}}_{kG}(V)=k$.
\begin{enumerate}
\item[i.] Then $\underline{\mathrm{End}}_{kG}(\Omega(V))=k$, and $R(G,V)$ and $R(G,\Omega(V))$ are isomorphic.
\item[ii.] There is a non-projective indecomposable $kG$-module $V_0$ $($unique up to isomorphism$)$ such that $\underline{\mathrm{End}}_{kG}(V_0)=k$, $V$ is isomorphic to $V_0\oplus P$ for some projective $kG$-module $P$, and $R(G,V)$ and $R(G,V_0)$ are isomorphic.
\end{enumerate}
\end{sublem}

The following result can be proved similarly to \cite[Prop. 5.2]{bc}, using \cite[Lemmas 5.3 and 5.4]{bc}.
\begin{subprop}
\label{prop:induceddef}
Let $L$ be a subgroup of $G$ and let $U$ be a finitely generated indecomposable $kL$-module with $\underline{\mathrm{End}}_{kL}(U)= k$. Suppose there exists an indecomposable $kG$-module $V$ with $\underline{\mathrm{End}}_{kG}(V)= k$ and a projective $kG$-module $P$ such that 
\begin{equation}
\label{eq:thatstheone}
\mathrm{Ind}_L^G U = V\oplus P.
\end{equation}
Assume further that
\begin{equation}
\label{eq:ohyeah!}
\mathrm{dim}_k \mathrm{Ext}_{kL}^1(U,U) = \mathrm{dim}_k \mathrm{Ext}_{kG}^1(V,V).
\end{equation}
Then $R(G,V)$ is isomorphic to $R(L,U)$.
\end{subprop}


\subsection{Stable equivalences of Morita type}
\label{ss:stabmor}
\setcounter{equation}{0}
\setcounter{figure}{0}

Suppose $G$ (resp. $H$) is a finite group, and let $A$ (resp. $B$) be a block of $WG$ (resp. $WH$). For $R\in\mathrm{Ob}(\mathcal{\hat{C}})$, define $RA$ (resp. $RB$) to be the block algebra in $RG$ (resp. $RH$) corresponding to $A$ (resp. $B$), i.e. $RA=R\otimes_W A$ (resp. $RB=R\otimes_W B$). Let $\Gamma$ be $A$ or $B$. Then $\Gamma$ is a $W$-algebra that is projective as a $W$-module.  Moreover, $\Gamma$ is a symmetric $W$-algebra in the sense that $\Gamma$ is isomorphic to its $W$-linear dual $\check{\Gamma}=\mathrm{Hom}_W(\Gamma,W)$, as $\Gamma$-$\Gamma$-bimodules.
In the following,  $\Gamma\mbox{-mod}$ denotes the category of finitely generated left  $\Gamma$-modules, and $\Gamma\mbox{-\underline{mod}}$ denotes the $W$-stable category, i.e. the quotient category of $\Gamma\mbox{-mod}$ by the subcategory of relatively $W$-projective modules. Recall that a $\Gamma$-module is called relatively $W$-projective if it is isomorphic to a direct summand of $\Gamma\otimes_W M$ for some $W$-module $M$.

In \cite{bl1}, it was shown that a split-endomorphism two-sided tilting complex (as introduced by Rickard \cite{rickard2}) for the derived categories of bounded complexes of finitely generated modules over $A$, resp. $B$, preserves the versal deformation rings of bounded complexes of finitely generated modules over $kA$, resp. $kB$.
It follows from a result by Rickard (see \cite{rickard1} and \cite[Prop. 6.3.8]{kozim}) that a derived equivalence 
between the derived categories of bounded complexes $D^b(A\mbox{-{mod}})$ and $D^b(B\mbox{-{mod}})$ induces
a stable equivalence between $A\mbox{-\underline{mod}}$ and $B\mbox{-\underline{mod}}$ which is itself induced by a bimodule. More precisely we get the following definition of stable equivalence of Morita type going back to Brou\'{e} \cite{broue1}.

\begin{subdfn}
\label{def:stabeq} 
Let $M$ be a $B$-$A$-bimodule and $N$ an $A$-$B$-bimodule. We say $M$ and $N$ induce a stable equivalence of Morita type between $A$ and $B$, if $M$ and $N$ are projective both as left and as right modules, and if 
\begin{eqnarray}
\label{eq:stab}
N\otimes_BM&\cong& A\oplus P \quad\mbox{ as $A$-$A$-bimodules, and} \\ \nonumber
M\otimes_AN&\cong& B\oplus Q \quad\mbox{ as $B$-$B$-bimodules}, 
\end{eqnarray}
where $P$ is a projective $A$-$A$-bimodule, and $Q$ is a projective $B$-$B$-bimodule.
In particular, $M\otimes_A-$ and $N\otimes_B-$ induce mutually inverse equivalences between the $W$-stable module categories $A\mbox{-\underline{mod}}$ and $B\mbox{-\underline{mod}}$.
\end{subdfn}

We now prove that stable equivalences of Morita type preserve versal deformation rings.

\begin{sublem}
\label{lem:stabmordef}
Let $A$ and $B$ be blocks of group rings over $W$ as above. Suppose that $M$ is a $B$-$A$-bimodule and $N$ is an $A$-$B$-bimodule which induce a stable equivalence of Morita type between $A$ and $B$. Let $V$ be a finitely generated $kA$-module, and define $V'=(k\otimes_WM)\otimes_{kA}V$. Then $R(G,V)$ is isomorphic to $R(H,V')$ in $\hat{\mathcal{C}}$.
\end{sublem}

\begin{proof}
Suppose $R\in\mathrm{Ob}(\mathcal{C})$ is Artinian. Then $M_R=R\otimes_W M$ is projective as left $RB$-module and as right $RA$-module, and $N_R=R\otimes_W N$ is projective as left $RA$-module and as right $RB$-module.
Since $M_R\otimes_{RA} (N_R) \cong R\otimes_W(M\otimes_AN)$, we have, using (\ref{eq:stab}),
$$M_R\otimes_{RA} N_R\cong RB\oplus P_R \quad\mbox{ as $RA$-$RA$-bimodules},$$
where $P_R=R\otimes_W P$ is a projective $RA$-$RA$-bimodule.
Similarly, we get 
$$N_R\otimes_{RB} M_R\cong RA \oplus Q_R \quad\mbox{ as $RB$-$RB$-bimodules},$$
where $Q_R=R\otimes_W Q$ is a projective $RB$-$RB$-bimodule.

Let $X=Q_k\otimes_{kA}V$ and $Y=P_k\otimes_{kB}V'$. Then since $Q_k$ is a projective $kA$-$kA$-bimodule, it follows that $X$ is a projective $kA$-module. Similarly, $Y$ is a projective $kB$-module. Let $X_R$ be the projective $RA$-cover of $X$ and $Y_R$ the projective $RB$-cover of $Y$, which exist since we assume $R$ to be Artinian. Then $X_R$ defines a lift $(X_R,\pi_X)$ of $X$ over $R$ and $Y_R$ defines a lift $(Y_R,\pi_Y)$ over $R$. Because $R$ is commutative local, every lift of $X$ (resp. $Y$) over $R$ is isomorphic to $(X_R,\pi_X)$ (resp. $(Y_R,\pi_Y)$).

Let now $(U,\phi)$ be a lift of $V$ over $R$. Then $U$ is a finitely generated $RA$-module. Define $U'=M_R\otimes_{RA}U$. 
Since $M$ is a finitely generated projective right $RA$-module and since $U$ is a finitely generated free $R$-module, it follows that $U'$ is a finitely generated projective, and hence free, $R$-module. Moreover
\begin{equation}
\label{eq:lala}
U'\otimes_Rk=(M_R\otimes_{RA}U)\otimes_Rk = M_k\otimes_{kA}(U\otimes_Rk)\xrightarrow{M_k\otimes\phi} M_k\otimes_{kA}V =V'.
\end{equation}
This means that $(U',\phi')=(M_R\otimes_{RA}U,M_k\otimes\phi)$ is a lift of $V'$ over $R$. We therefore obtain for all $R\in\mathrm{Ob}(\mathcal{C})$ a well-defined map 
$$\tau_R:\mathrm{Def}_G(V,R)\to\mathrm{Def}_H(V',R).$$
We need to show that $\tau_R$ is bijective. Let $U''=N_R\otimes_{RB}U'$ and $\phi''=N_k\otimes \phi'$. Similarly to (\ref{eq:lala}), it follows that $(U'',\phi'')=(N_R\otimes_{RB}U',N_k\otimes \phi')$ is a lift of $V''=N_k\otimes_{kB}V'$ over $R$. We have 
\begin{eqnarray}
\label{eq:olala}
(U'',\phi'')&=&(N_R\otimes_{RB}(M_R\otimes_{RA}U),N_k\otimes(M_k\otimes\phi))\\
&\cong&((RA\oplus Q_R)\otimes_{RA}U, (kA\oplus Q_k)\otimes\phi) \nonumber\\
&\cong&(U\oplus (Q_R\otimes_{RA}U), \phi\oplus(Q_k\otimes\phi)) \nonumber\\
&\cong& (U\oplus X_R,\phi\oplus\pi_X)\nonumber
\end{eqnarray}
as lifts of $V''$ over $R$, where the last isomorphism follows, since $(Q_R\otimes_{RA}U, Q_k\otimes\phi)$ is a lift of $X=Q_k\otimes_{kA}V$ over $R$. Moreover,
$$V''=N_k\otimes_{kB}(M_k\otimes_{kA}V)\cong (kA\oplus Q_k)\otimes_{kA}V\cong V \oplus (Q_k\otimes_{kA}V)=V\oplus X.$$
Hence it follows by \cite[Prop. 2.6]{bc} that $\tau_R$ is injective. 

Now let $(L,\psi)$ be a lift of $V'=M_k\otimes_{kA}V$ over $R$.  Then $(L',\psi')= (N_R\otimes_{RB}L,N_k\otimes\psi)$ is a lift of $V''=N_k\otimes_{kB}V'\cong V\oplus X$ over $R$. By \cite[Prop. 2.6]{bc}, there exists a lift $(U,\phi)$ of $V$ over $R$ such that $(L',\psi')\cong (U\oplus X_R,\phi\oplus\pi_X)$. Arguing similarly as in (\ref{eq:olala}), we then have that $(L',\psi')$ is isomorphic to $(U'',\phi'')=(N_R\otimes_{RB}U',N_k\otimes\phi')$ where $(U',\phi')=(M_R \otimes_{RA}U,M_k\otimes\phi)$. Therefore, $(M_R\otimes_{RA}L' ,M_k\otimes\psi')\cong (M_R\otimes_{RA}U'',M_k\otimes\phi'')$.
Arguing again similarly as in (\ref{eq:olala}), we have
\begin{eqnarray*}
(M_R\otimes_{RA}L' ,M_k\otimes\psi') 
&\cong&(L\oplus Y_R,\psi\oplus\pi_Y),
\mbox{ and }\\
(M_R\otimes_{RA}U'',M_k\otimes\phi'')
&\cong&(U' \oplus Y_R,\phi'\oplus\pi_Y).
\end{eqnarray*}
Thus by \cite[Prop. 2.6]{bc}, it follows that $(L,\psi)\cong (U',\phi')$, i.e. $\tau_R$ is surjective.
Since the deformation functors $\hat{F}_V$ and $\hat{F}_{V'}$ are continuous, this implies that they are naturally isomorphic. Hence the versal deformation rings $R(G,V)$ and $R(H,V')$ are isomorphic in $\hat{\mathcal{C}}$.
\end{proof}

Using the same type of arguments as in the proof of Lemma \ref{lem:stabmordef}, but restricting our attention to Artinian objects $R$ in $\mathcal{C}$ that are $k$-algebras, we get the following weaker result about versal deformation rings modulo $p$.

\begin{sublem}
\label{lem:weakstabmordef}
Let $A$ and $B$ be blocks of group rings over $W$ as above. Suppose that $M_k$ is a $kB$-$kA$-bimodule and $N_k$ is a $kA$-$kB$-bimodule which induce a stable equivalence of Morita type between $kA$ and $kB$. Let $V$ be a finitely generated $kA$-module, and define $V'=M_k\otimes_{kA}V$. Then $R(G,V)/pR(G,V)\cong R(H,V')/pR(H,V')$.
\end{sublem}

\begin{subrem}
\label{rem:stabmor}
Using the notations of Lemma \ref{lem:stabmordef} (resp. of Lemma \ref{lem:weakstabmordef}), suppose that the stable endomorphism ring $\underline{\mathrm{End}}_{kG}(V)=k$. Then it follows that also $\underline{\mathrm{End}}_{kH}(V')=k$. By Lemma \ref{lem:defhelp}(ii), there exists a non-projective indecomposable $kH$-module $V'_0$ (unique up to isomorphism) which is a direct summand of $V'$ with $\underline{\mathrm{End}}_{kH}(V'_0)=k$ and $R(H,V')\cong R(H,V'_0)$. In the situation of Lemma \ref{lem:stabmordef}, it then follows that $R(G,V)\cong R(H,V'_0)$. In the situation of Lemma \ref{lem:weakstabmordef}, we have at least $R(G,V)/pR(G,V)\cong R(H,V'_0)/pR(H,V'_0)$
\end{subrem}

\begin{subrem}
\label{rem:stablegood}
Let $\Lambda$ be a finite dimensional $k$-algebra, and denote by $\mathrm{mod}_{\mathcal{P}}(\Lambda)$ the full subcategory of $\Lambda\mbox{-mod}$ whose objects are the modules which have no non-zero projective summands.
Suppose $\Lambda'$ is another finite dimensional $k$-algebra, and $F:\Lambda\mbox{-\underline{mod}}\to \Lambda'\mbox{-\underline{mod}}$ is a stable equivalence. 
Let
$$0\to A \xrightarrow{\genfrac{(}{)}{0pt}{}{\mbox{\tiny $f$}}{\mbox{\tiny $s$}}} B \coprod P \xrightarrow{(g,t)} C\to 0$$
be an almost split sequence in $\Lambda\mbox{-mod}$ where $A,B,C$ are in $\mathrm{mod}_{\mathcal{P}}(\Lambda)$, $B$ is non-zero and $P$ is projective. Then, by \cite[Prop. X.1.6]{ars}, for any morphism $g':F(B)\to F(C)$ with $F(\underline{g})=\underline{g}'$ there is an almost split sequence
$$0\to F(A) \xrightarrow{\genfrac{(}{)}{0pt}{}{\mbox{\tiny $f'$}}{\mbox{\tiny $u$}}} F(B) \coprod P' \xrightarrow{(g',v)} F(C)\to 0$$
in $\Lambda'\mbox{-mod}$ where $P'$ is projective and $F(\underline{f})=\underline{f}'$.

Moreover, by \cite[Cor. X.1.9 and Prop. X.1.12]{ars}, if $\Lambda$ and $\Lambda'$ are selfinjective with no blocks of Loewy length $2$, then the stable Auslander-Reiten quivers of $\Lambda$ and $\Lambda'$ are isomorphic stable translation quivers, and $F$ commutes with $\Omega$.
\end{subrem}

\subsection{Universal deformation rings that are quotient rings of $W[[t]]$}
\label{ss:prelimdefmore}
\setcounter{equation}{0}
\setcounter{figure}{0}

In this subsection, we provide a few results that help determine universal deformation rings that are certain quotient rings of $W[[t]]$. As before, let $G$ be a finite group. The first result deals with universal deformation rings modulo $p$.

\begin{sublem}
\label{lem:trythis}
Let $Y$ be a finitely generated uniserial $kG$-module 
satisfying $\mathrm{End}_{kG}(Y)=k$ 
and $\mathrm{Ext}^1_{kG}(Y,Y)=k$. Suppose $Y$ has descending radical series $(T_1,T_2,\ldots, T_\ell)$ where $\ell\ge 1$ and $T_1,\ldots,T_\ell$ are simple $kG$-modules, not necessarily distinct.
Assume there exists an integer $s\ge 1$ such that the projective cover $P_{T_1}$ has the form $P_{T_1}= \begin{array}{cc}\multicolumn{2}{c}{T_1}\\U_1&U_2\\ \multicolumn{2}{c}{T_1} \end{array}$ where $U_1$ and $U_2$ are uniserial $kG$-modules, $U_1$ may be zero, and $U_2$ has descending radical series
$$(T_2,\ldots,T_\ell,T_1,T_2,\ldots,T_\ell,\ldots,T_1,T_2,\ldots,T_\ell)$$
of length $\ell p^s-1$.
Define $\overline{U}$ to be the uniserial $kG$-module $\overline{U}=\begin{array}{c}T_1\\U_2\end{array}$, and suppose $\mathrm{Ext}^1_{kG}(\overline{U},Y)=0$. Then the universal deformation ring of $Y$ modulo $p$ is $\overline{R}=R(G,Y)/pR(G,Y)\cong k[t]/(t^{p^s})$, and the universal mod $p$ deformation of $Y$ over $\overline{R}$ is represented by the $kG$-module $\overline{U}$.
\end{sublem}

\begin{proof}
By assumption, $\mathrm{Ext}^1_{kG}(Y,Y)=k$, which implies that $\overline{R}\cong k[t]/(t^r)$ for some $r\ge 1$. The module $\overline{U}$ is a uniserial $kG$-module of length $\ell p^s$ with descending radical series
$$(T_1,\ldots,T_\ell,T_1\ldots,T_\ell,\ldots,T_1,\ldots,T_\ell) = (T_1,\ldots,T_\ell)^{p^s}.$$
If we let $t$ act as the shift down by $\ell$, it follows that $\overline{U}$ is a free $k[t]/(t^{p^s})$-module which is a lift of $Y$ over $k[t]/(t^{p^s})$. Hence there is a $k$-algebra homomorphism
$$\phi:\overline{R}\to k[t]/(t^{p^s})$$
corresponding to $\overline{U}$. Since $\overline{U}$ is indecomposable as a $kG$-module, it follows that $\phi$ is surjective. We now show that $\phi$ is a $k$-algebra isomorphism. Suppose this is false. Then there exists a surjective $k$-algebra homomorphism $\phi_1:\overline{R}\to k[t]/(t^{p^s+1})$ such that $\pi\phi_1=\phi$ where $\pi:k[t]/(t^{p^s+1})\to k[t]/(t^{p^s})$ is the natural projection. Let $\overline{U}_1$ be a lift of $Y$ over $k[t]/(t^{p^s+1})$ relative to $\phi_1$. Then $\overline{U}_1$ is a lift of $\overline{U}$ over $k[t]/(t^{p^s+1})$ with $t^{p^s}\overline{U}_1\cong Y$. Thus we have a short exact sequence of
$k[t]/(t^{p^s+1})\,G$-modules
\begin{equation}
\label{eq:thesos}
0\to t^{p^s}\overline{U}_1\to \overline{U}_1\to \overline{U}\to 0.
\end{equation}
We now show that this sequence cannot split as a sequence of $kG$-modules. Suppose it splits. Then $\overline{U}_1\cong Y\oplus \overline{U}$ as $kG$-modules. Let $z=\left(\begin{array}{c}y\\u\end{array}\right)\in Y\oplus \overline{U} \cong \overline{U}_1$. Then $t$ acts on $z$ as multiplication by the matrix $$A_t=\left(\begin{array}{cc}0&\alpha\\0&\mu_t\end{array}\right)$$ 
where $\alpha:\overline{U}\to Y$ is a surjective $kG$-module homomorphism, and $\mu_t$ is multiplication by $t$ on $\overline{U}$. Since $t^{p^s}\overline{U}_1\cong Y$, there exists a non-zero $z=\left(\begin{array}{c}y\\u\end{array}\right)\in Y\oplus \overline{U} \cong \overline{U}_1$ with $(A_t)^{p^s} z\neq 0$. But, since $\mathrm{End}_{kG}(Y)=k$, $\alpha$ corresponds to the isomorphism $\overline{U}/t\overline{U}\cong Y$ which means that the kernel of $\alpha$ is $t\overline{U}$. Thus
$$(A_t)^{p^s} \left(\begin{array}{c}y\\u\end{array}\right) = \left(\begin{array}{c}\alpha(\mu_t^{p^s-1}(u))\\\mu_t^{p^s}(u)\end{array}\right)
= \left(\begin{array}{c}0\\0\end{array}\right),$$
which gives a contradiction. Hence the short exact sequence (\ref{eq:thesos}) does not split as a sequence of $kG$-modules. Since
$\mathrm{Ext}^1_{kG}(\overline{U},Y)=0$ by assumption and $t^{p^s}\overline{U}_1\cong Y$, this is impossible. Therefore, $\overline{U}_1$ does not exist, which means that $\phi$ is a $k$-algebra isomorphism. Thus $\overline{R}\cong k[t]/(t^{p^s})$, and the universal mod $p$ deformation of $Y$ over $\overline{R}$ is represented by the $kG$-module $\overline{U}$.
\end{proof}

The next result analyzes when a finitely generated $kG$-module can be lifted to a $WG$-module which is free as a  $W$-module. In the following, $F$ denotes the fraction field of $W$.

\begin{sublem}
\label{lem:Wlift}
Let $V$ be a finitely generated $kG$-module such that there is a non-split short exact sequence of $kG$-modules
\begin{equation}
\label{eq:ohno}
0\to Y_2\to V\to Y_1\to 0
\end{equation}
with $\mathrm{Ext}^1_{kG}(Y_1,Y_2)=k$.
Assume that for $i\in\{1,2\}$, there exists a $WG$-module $X_i$ which is a lift of $Y_i$ over $W$; in particular $X_i$ is free as $W$-module. Suppose that 
\begin{equation}
\label{eq:ohno2}
\mathrm{dim}_F\;\mathrm{Hom}_{FG}(F\otimes_WX_1,F\otimes_WX_2) =
\mathrm{dim}_k\;\mathrm{Hom}_{kG}(Y_1,Y_2)-1.
\end{equation}
Then there exists a $WG$-module $X$ which is a lift of $V$ over $W$; in particular, $X$ is free as $W$-module. 
\end{sublem}

\begin{proof}
Since $X_1$ and $X_2$ are free as $W$-modules, we see, using spectral sequences, that 
\begin{eqnarray*}
\HH^i(G,\mathrm{Hom}_W(X_1,X_2))&=&\mathrm{Ext}^i_{WG}(X_1,X_2), \mbox{ and}\\
\HH^i(G,\mathrm{Hom}_k(Y_1,Y_2))&=&\mathrm{Ext}^i_{kG}(Y_1,Y_2)
\end{eqnarray*}
for all $i\ge 0$. From the short exact sequence
$$0\to \mathrm{Hom}_W(X_1,X_2) \xrightarrow{\cdot p}  \mathrm{Hom}_W(X_1,X_2) \xrightarrow{\mathrm{mod}\; p} \mathrm{Hom}_k(Y_1,Y_2)\to 0$$
we obtain a long exact cohomology sequence
$$\xymatrix @R=.6pc{
0\ar[r] & \mathrm{Hom}_{WG}(X_1,X_2)\ar[r]^{\mu^0_*}& \mathrm{Hom}_{WG}(X_1,X_2) \ar[r]^{\pi^0_*} & \mathrm{Hom}_{kG}(Y_1,Y_2) 
\ar `r/6pt[d] `[ld] `^d[ll] _{\delta^0} [ddll]  
\\&&&&\\  
&\mathrm{Ext}^1_{WG}(X_1,X_2)\ar[r]^{\mu^1_*}& \mathrm{Ext}^1_{WG}(X_1,X_2) \ar[r]^{\pi^1_*} & \mathrm{Ext}^1_{kG}(Y_1,Y_2)\ar[r]&\cdots
}$$
where for $i\in\{0,1\}$, $\mu^i_*$ stands for multiplication by $p$ and $\pi^i_*$ stands for reduction modulo $p$. To prove that there is a lift of $V$ over $W$, it is enough to prove that $\pi^1_*$ is surjective.

Let $n=\mathrm{dim}_F\;\mathrm{Hom}_{FG}(F\otimes_WX_1,F\otimes_WX_2)$.
Since $X_1$ and $X_2$ are free $W$-modules of finite rank, it follows that $$F\otimes_W\mathrm{Hom}_{WG}(X_1,X_2) \cong \mathrm{Hom}_{FG}(F\otimes_WX_1,F\otimes_WX_2)$$
as $F$-modules. Thus $\mathrm{Hom}_{WG}(X_1,X_2)$ is a free $W$-module of rank $n$. Hence $\mathrm{Im}(\pi^0_*)= k^n$ and thus $\mathrm{Ker}(\delta^0)=k^n$, as $W$-modules. On the other hand by (\ref{eq:ohno2}), $\mathrm{Hom}_{kG}(Y_1,Y_2)= k^{n+1}$. Thus $\mathrm{Im}(\delta^0)=k$, which implies $\mathrm{Ker}(\mu^1_*) =k$. 
Hence $\mathrm{Ext}^1_{WG}(X_1,X_2)$ is a non-zero finitely generated $W$-module, and so  the image of $\pi^1_*$, which is isomorphic to the cokernel of $\mu^1_*$, is non-trivial. Since $\mathrm{Ext}^1_{kG}(Y_1,Y_2)=k$, it follows that $\pi^1_*$ is surjective.
\end{proof}

For the remainder of this section, we consider the case $p=2$ and prove some results which help determine the universal deformation ring $R$ provided $R/2R\cong k[t]/(t^{2^n})$ for some positive integer $n$. The first result is a generalization of \cite[Lemma 2.2]{bc4.9}.

\begin{sublem}
\label{lem:herewegoagain}
Suppose $k$ has characteristic $p=2$.
Let $d\ge 3$ be an integer, and let $f(t)$ be a monic polynomial in $W[t]$ of degree $2^{d-2}-1$ where all non-leading coefficients lie in $2W$.
Suppose $R$ is a complete local Noetherian $W$-algebra with residue field $k$ for which there is a continuous surjection $\tau:R \to W[[t]]/(f(t))$ and an isomorphism $\mu:R/2R \to k[s]/(s^{2^{d-2}})$ of $W$-algebras. Then $R$ is isomorphic to $W[[t]]/(f(t)(t-2\gamma),\alpha2^mf(t))$ as a $W$-algebra, where $\gamma\in W$, $\alpha\in\{0,1\}$ and $0 < m\in\mathbb{Z}$. 
\end{sublem}

\begin{proof}
It follows from the assumptions that there is a continuous $W$-algebra surjection $\psi:W[[t]] \to R$.  Then $\tau\circ \psi:W[[t]] \to W[[t]]/(f(t))$ is a surjective $W$-algebra homomorphism.  Hence modulo 2, $t$ is sent to a generator of the maximal ideal of $W[[t]]/(2,f(t))\cong k[[t]]/(t^{2^{d-2}-1})$, which means that $t$ is sent to $u(t)\cdot t$ for some unit $u(t)$ in $W[[t]]/(2,f(t))$. Since we can lift $u(t)$ to a unit in $W[[t]]$, we see that $\tau(\psi(t))=v(t)\cdot t +2b(t)$ for certain $v(t)\in W[[t]]^*$ and $b(t)\in W[[t]]$. Composing $\psi$ with the continuous $W$-algebra automorphism $W[[t]] \to W[[t]]$  which sends $t$ to $v(t)^{-1}(t - 2b(t))$, we can assume that $\tau (\psi(t)) = t$.  Then $\mathrm{ker}(\tau \circ \psi)$ is the ideal  $(f(t))= W[[t]] \cdot f(t)$.  Hence the kernel $J$ of $\psi$ is contained in the ideal $(f(t))$.  Moreover, $J$ is properly contained in $(f(t))$, since otherwise $J = (f(t))$ and then $R \cong W[[t]]/(f(t))$ as $W$-algebras. But this is impossible since $R/2R \cong k[s]/(s^{2^{d-2}})$. 

The maximal ideal of $W[[t]]$ is generated by $2$ and $t$.  So the maximal ideal
of $R/2R$ 
is generated by the image of $t$ under the surjection $W[[t]] \to R/2R$ induced by $\psi:W[[t]] \to R$.
However,
the isomorphism $\mu:R/2R \to k[s]/(s^{2^{d-2}})$ shows that as a $W$-module, $R/2R$
is generated by $1$ together with the powers $\xi,\xi^2,\ldots,\xi^{2^{d-2}-1}$ for any generator $\xi$ of the maximal ideal of $R/2R$. 
So $R/2R$ is generated as a $W$-module by the images of $1,t,t^2,\ldots, t^{2^{d-2}-1}$.
Hence the image of $W \oplus Wt \oplus\cdots \oplus Wt^{2^{d-2}-1}\subset W[[t]]$ under $\psi:W[[t]] \to R$ must be all of $R$ since $R$ is complete.  Thus $\psi(t^{2^{d-2}}) = \psi(a_0 + a_1 t +\cdots +a_{2^{d-2}-1}t^{2^{d-2}-1})$
for some $a_0, a_1,\ldots,a_{2^{d-2}-1} \in W$.  This means that $t^{2^{d-2}} - (a_0 + a_1 t +\cdots + a_{2^{d-2}-1}t^{2^{d-2}-1}) = j \in J$.  But $J \subseteq (f(t))$, so 
\begin{equation}
\label{eq:ohwell!}
t^{2^{d-2}} - (a_0 + a_1 t +\cdots + a_{2^{d-2}-1}t^{2^{d-2}-1}) = f(t)\cdot q(t)
\end{equation}
for a unique $q(t)\in W[[t]]$. Let $q(t)=c_0+t+\sum_{i=1}^\infty c_it^i$. Suppose there exists $i\ge 1$ with $c_i\neq 0$. Letting $r=\mathrm{min}\{\mathrm{ord}_2(c_i)\;|\;i\ge 1\}$, we can rewrite $q(t)=c_0+t+2^rtw(t)$ for some $w(t)\in W[[t]]$ which 
is not congruent to the zero power series modulo $2$.
Comparing the coefficients modulo $2^{r+1}$ of the terms of degree at least $2^{d-2}$ on both sides of (\ref{eq:ohwell!}), we see that $2^r w(t)$ must be congruent to the zero power series modulo $2^{r+1}$, which is impossible. Hence 
$q(t)=t+c_0$ for some $c_0\in W$.   Moreover, $c_0$ is in $2W$, since otherwise $c_0 \in W^*$ and hence $t+c_0 \in W[[t]]^*$. But then $f(t) = (t+c_0)^{-1} j\in J$, which is impossible since we showed that $J$ is properly contained in $(f(t))$.  So $c_0 = -2\gamma$
for some $\gamma \in W$.  

This means $f(t)(t-2\gamma)W[[t]] \subseteq J \subset f(t) W[[t]]$. Hence $J=f(t)J'$ where 
$$(t-2\gamma) =(t-2\gamma)W[[t]] \subseteq J'\subset W[[t]]$$
and $J'\neq W[[t]]$.
Therefore, $J'/(t-2\gamma)$ is a proper ideal of $W[[t]]/(t-2\gamma)\cong W$, and hence either zero or generated by a positive power of $2$. It follows that $J'=(t-2\gamma,\alpha2^m)$ where $\alpha\in\{0,1\}$ and $m\in\mathbb{Z}^+$. Thus $J=(f(t) (t-2\gamma), \alpha 2^mf(t))$.
\end{proof}

We now specialize to a particular $f(t)$.
\begin{subdfn}
\label{def:seemtoneed}
Suppose $k$ has characteristic $p=2$. Let $F$ be the fraction field of $W$, and
fix an algebraic closure $\overline{F}$ of $F$.
Let $d\ge 3$ be an integer, 
and let $\zeta_{2^\ell}$ be a fixed primitive $2^\ell$-th root of unity in $\overline{F}$ for $2\le \ell\le d-1$. 
\begin{enumerate}
\item[i.] Define
$$p_d(t)=\prod_{\ell=2}^{d-1} \mathrm{min.pol.}_F(\zeta_{2^\ell}+\zeta_{2^\ell}^{-1}),$$
and let $R'=W[[t]]/(p_d(t))$.
\item[ii.] Let $Z=\langle \sigma\rangle$ be a cyclic group of order $2^{d-1}$, and let $\tau:Z\to Z$ be the group automorphism sending $\sigma$ to $\sigma^{-1}$. Then $\tau$ can be extended to a $W$-algebra automorphism of the group ring $WZ$ which will again be denoted by $\tau$. Let $T(\sigma^2)=1+\sigma^2+\sigma^4+\cdots +\sigma^{2^{d-1}-2}$, and define
$$S'= (WZ)^{\langle \tau\rangle}/\left(T(\sigma^2),\sigma T(\sigma^2)\right).$$
\end{enumerate}
\end{subdfn}

\begin{subrem}
\label{rem:ohyeah}
The minimal polynomial $\mathrm{min.pol.}_F(\zeta_{2^\ell}+\zeta_{2^\ell}^{-1})$ for $\ell\ge 2$ is as follows:
\begin{eqnarray*}
\mathrm{min.pol.}_F(\zeta_{2^2}+\zeta_{2^2}^{-1})(t)&=&t,\\ \mathrm{min.pol.}_F(\zeta_{2^{\ell}}+\zeta_{2^{\ell}}^{-1})(t)
&=&\left(\mathrm{min.pol.}_F(\zeta_{2^{\ell-1}}+\zeta_{2^{\ell-1}}^{-1})(t)\right)^2-2 \qquad \mbox{for $\ell\ge 3$}.
\end{eqnarray*}

The $W$-algebra $R'$ from Definition \ref{def:seemtoneed} is a complete local Noetherian ring with residue field $k$. Moreover, 
\begin{eqnarray*}
F\otimes_W R'&\cong& \prod_{\ell=2}^{d-1} F(\zeta_{2^\ell}+\zeta_{2^\ell}^{-1})\quad\mbox{ as $F$-algebras,}\\
k\otimes_W R'&\cong& k[t]/(t^{2^{d-2}-1})\quad\mbox{ as $k$-algebras.}
\end{eqnarray*}
Additionally, for any sequence $\left(r_\ell\right)_{\ell=2}^{d-1}$ of odd integers, $R'$ is isomorphic to the $W$-subalgebra of
$$\prod_{\ell=2}^{d-1} W[\zeta_{2^\ell}+\zeta_{2^\ell}^{-1}]$$
generated by the element 
$\left(\zeta_{2^\ell}^{r_\ell}+\zeta_{2^\ell}^{-r_\ell}\right)_{\ell=2}^{d-1}$.
\end{subrem}

\begin{sublem}
\label{lem:dihedraltrick}
Using the notations of Definition $\ref{def:seemtoneed}$, there is a continuous $W$-algebra isomorphism $\rho:R'\to S'$ with $\rho(t)=\sigma+\sigma^{-1}$. 

Suppose now that $D_{2^d}$ is a dihedral group of order $2^d$.
If in Lemma $\ref{lem:herewegoagain}$ the polynomial $f(t)$ is taken to be equal to $p_d(t)$, $\alpha=1$ and $m=1$, then the ring $R$ in this Lemma is isomorphic to a subquotient algebra of the group ring $WD_{2^d}$.
\end{sublem}

\begin{proof}
The ring of invariants $(WZ)^{\langle \tau\rangle}$ is a free $W$-module with basis
$$\{1, \sigma+\sigma^{-1}, \sigma^2+\sigma^{-2},\ldots, \sigma^{2^{d-2}-1}+\sigma^{-2^{d-2}+1},\sigma^{2^{d-2}}\}.$$
Hence the $W$-rank of this module is $2^{d-2}+1$. Expanding expressions of the form
$(\sigma+\sigma^{-1})^\ell$ for various $\ell$, one sees that as $W$-algebra $(WZ)^{\langle \tau\rangle}$ is generated by $(\sigma+\sigma^{-1})$ and $\sigma^{2^{d-2}}$. 
In $S'$, the residue class of $\sigma^{2^{d-2}}$ can be expressed as a polynomial in $(\sigma+\sigma^{-1})$, which means that $S'$ is generated as a $W$-algebra by the residue class of $(\sigma+\sigma^{-1})$. Since the residue class of $(\sigma^{2^{d-2}-1}+\sigma^{-2^{d-2}+1})$ can also be expressed as a polynomial in $(\sigma+\sigma^{-1})$ in $S'$ and since $S'$ has no torsion, we conclude that $S'$ is a free $W$-module of rank $2^{d-2}+1-2 = 2^{d-2}-1$.

Define $\hat{\rho}:W[[t]]\to S'$ to be the continuous $W$-algebra homomorphism sending $t$ to the residue class of $(\sigma+\sigma^{-1})$. Then $\hat{\rho}$ is surjective. 
Using the description of the minimal polynomial of $(\zeta_{2^\ell}+\zeta_{2^\ell}^{-1})$, $2\le \ell\le d-1$, over $F$ in Remark \ref{rem:ohyeah}, we see that 
$$[\mathrm{min.pol.}_F(\zeta_{2^\ell}+\zeta_{2^\ell}^{-1})]\,(\sigma+\sigma^{-1})
=\sigma^{2^{\ell-2}}+\sigma^{-2^{\ell-2}}.$$
Hence we have
\begin{equation}
\label{eq:ole1}
p_d(\sigma+\sigma^{-1})=\prod_{\ell=2}^{d-1} (\sigma^{2^{\ell-2}}+\sigma^{-2^{\ell-2}}),
\end{equation}
and we see by induction that the latter is equal to 
\begin{equation}
\label{eq:ole2}
(\sigma+\sigma^{-1})+(\sigma^3+\sigma^{-3})+\cdots + (\sigma^{2^{d-2}-1}+\sigma^{-2^{d-2}+1})=\sigma T(\sigma^2)
\end{equation}
which is zero in $S'$. Thus $p_d(t)$ lies in the kernel of $\hat{\rho}$. This means that
we obtain a surjective continuous $W$-algebra homomorphism
$$\rho:R'=W[[t]]/(p_d(t)) \to S'.$$
Since both $R'$ and $S'$ are free over $W$ of rank $2^{d-2}-1$, it follows that $\rho$ is a continuous $W$-algebra isomorphism. In particular, $R'$ is isomorphic to a subquotient algebra of the group ring $WD_{2^d}$ when $D_{2^d}$ is a dihedral group of order $2^d$.

Suppose now that in Lemma $\ref{lem:herewegoagain}$ the polynomial $f(t)$ is taken to be equal to $p_d(t)$, $\alpha=1$ and $m=1$. Then the ring $R$ in this Lemma is isomorphic to $W[[t]]/(p_d(t)(t-2),2\,p_d(t))$ as a $W$-algebra. Thus to finish the proof of Lemma \ref{lem:dihedraltrick}, it suffices to show that the ring $W[[t]]/(p_d(t)(t-2))$ is isomorphic to a subquotient algebra of $WD_{2^d}$. Define
$$\Theta=(WZ)^{\langle\tau\rangle}/\left(T(\sigma^2)- \sigma T(\sigma^2)\right).$$
Then $\Theta$ is isomorphic to a subquotient algebra of $WD_{2^d}$ and it is generated as a $W$-algebra by the residue class of $(\sigma+\sigma^{-1})$. Moreover, $\Theta$ is a free $W$-module of rank $2^{d-2}$, since the ideal $\left(T(\sigma^2)- \sigma T(\sigma^2)\right)$ is generated over $W$ by $T(\sigma^2)- \sigma T(\sigma^2)$. Define a continuous $W$-algebra homomorphism $\theta:W[[t]]\to \Theta$ by sending $t$ to the residue class of $(\sigma+\sigma^{-1})$.  Then, using (\ref{eq:ole1}) and (\ref{eq:ole2}), we see that
\begin{eqnarray*}
\theta(p_d(t)(t-2)) &=& p_d(\sigma+\sigma^{-1})((\sigma+\sigma^{-1})-2) \\
&=& \sigma T(\sigma^2)((\sigma+\sigma^{-1})-2)\\
&=& 2\left[ T(\sigma^2)- \sigma T(\sigma^2)\right]
\end{eqnarray*}
which is zero in $\Theta$. Hence $p_d(t)(t-2)$ lies in the kernel of $\theta$. This means that we obtain a surjective continuous $W$-algebra homomorphism
$$\overline{\theta}:W[[t]]/(p_d(t)(t-2)) \to \Theta.$$
Since both $W[[t]]/(p_d(t)(t-2))$ and $\Theta$ are free over $W$ of rank $2^{d-2}$, it follows that $\overline{\theta}$ is a continuous $W$-algebra isomorphism. Thus $W[[t]]/(p_d(t)(t-2))$ is isomorphic to a subquotient algebra of $WD_{2^d}$, which completes the proof of Lemma \ref{lem:dihedraltrick}.
\end{proof}


\section{Blocks with dihedral defect groups}
\label{s:dihedralsylow}
\setcounter{equation}{0}
\setcounter{figure}{0}

Let $k$ be an algebraically closed field of characteristic $p=2$, and let $W$ be the ring of infinite Witt vectors over $k$. 
Let $G$ be a finite group, and let $B$ be a block of $kG$ with dihedral defect groups which is Morita equivalent to the principal block of a finite simple group. From the classification by Gorenstein and Walter of the groups with dihedral Sylow $2$-subgroups in \cite{gowa}, it follows that there are three families of such blocks $B$, up to Morita equivalence:
\begin{enumerate}
\item[i.] the principal $2$-modular blocks of $\mathrm{PSL}_2(\mathbb{F}_q)$ where $q\equiv 1 \mod 4$, 
\item[ii.] the principal $2$-modular blocks of $\mathrm{PSL}_2(\mathbb{F}_q)$ where $q\equiv 3 \mod 4$, and 
\item[iii.] the principal $2$-modular block of the alternating group $A_7$.
\end{enumerate}
Note that in all cases (i) - (iii), $B$ contains precisely 3 isomorphism classes of simple modules. 

\begin{rem}
\label{rem:oherdmann}
In \cite{erd}, Erdmann classified all blocks with dihedral defect groups. It follows that if we consider all such blocks $B_0$ containing precisely 3 isomorphism classes of simple modules, there is one more family attached to case (iii) containing a Morita equivalence class of possible blocks for each defect $d\ge 3$. By \cite[\S X.4]{erd}, the blocks in this family having defect $d\ge 4$ cannot be excluded as possible blocks with dihedral defect groups. 

However, if we assume that $B_0$ is Morita equivalent to the principal block of $kH$ for some finite (not necessarily simple) group $H$, we can exclude this family as follows. Since $k$ has characteristic $2$, we can assume that $H$ has no normal subgroup of odd order. By \cite{gowa}, it then follows that $H$ is isomorphic to either a subgroup of $\mathrm{P\Gamma L}_2(\mathbb{F}_q)$ containing $\mathrm{PSL}_2(\mathbb{F}_q)$ for some odd prime power $q$, or to the alternating group $A_7$. Using a theorem by Clifford \cite[Hauptsatz V.17.3]{hup}, we see that the only possibility for $B_0$ to be Morita equivalent to a block in the bigger family attached to case (iii) occurs when $B_0$ has defect $d=3$, i.e. $B_0$ is Morita equivalent to the principal $2$-modular block of the alternating group $A_7$.
\end{rem}

The blocks in (i), (ii) and (iii) are all Morita equivalent to basic algebras of special biserial algebras. (For the relevant background on special biserial algebras we refer to \S\ref{s:prelimstring}.)
In \S \ref{ss:psl1}, \S \ref{ss:psl2} and \S \ref{ss:a7}, we give the quivers and relations for the basic algebras of these blocks,
together with their projective indecomposable modules and their decomposition matrices. 
In \S\ref{ss:ordinarydihedral}, we then state some results from \cite{brauer2} about the ordinary irreducible characters of $G$ which belong to $B$.

\subsection{The principal $2$-modular block of $\mathrm{PSL}_2(\mathbb{F}_q)$ when $q\equiv 1 \mod 4$}
\label{ss:psl1}
\setcounter{equation}{0}
\setcounter{figure}{0}

Let $G$ be a finite group, and let $B$ be a block of $kG$ which is Morita equivalent to the principal block of $k\mathrm{PSL}_2(\mathbb{F}_q)$ where $q\equiv 1 \mod 4$. Suppose that $2^d$ is the order of the defect groups of $B$, i.e. the order of the Sylow $2$-subgroups of $\mathrm{PSL}_2(\mathbb{F}_q)$. Then, by \cite{erd}, $B$ is Morita equivalent to the special biserial algebra $\Lambda=kQ/I$ where $Q$ is given in Figure \ref{fig:psl1A} and
$$I=\langle \gamma\beta,\delta\eta,(\eta\delta\beta\gamma)^{2^{d-2}} - (\beta\gamma\eta\delta)^{2^{d-2}}\rangle.$$
\begin{figure}[ht] \hrule \caption{\label{fig:psl1A} The quiver $Q$ for blocks as in \S$\ref{ss:psl1}$.}
$$\xymatrix @R=-.2pc {
&0&\\
Q= \quad 1\; \bullet \ar@<.8ex>[r]^(.7){\beta} \ar@<.9ex>[r];[]^(.3){\gamma}
& \bullet \ar@<.8ex>[r]^(.44){\delta} \ar@<.9ex>[r];[]^(.56){\eta} & \bullet\; 2}$$
\vspace{1ex}
\hrule
\end{figure}
We denote the irreducible $\Lambda$-modules by $S_0,S_1,S_2$, or, using short-hand, by $0,1,2$.
The radical series of the projective indecomposable $\Lambda$-modules (and hence of the projective indecomposable $B$-modules) are described in  Figure \ref{fig:psl1B} 
\begin{figure}[ht] \hrule \caption{\label{fig:psl1B} The radical series of the projective indecomposable modules for blocks as in \S$\ref{ss:psl1}$.}
$$P_0=\begin{array}{ccc}&0&\\
1&&2\\0&&0\\2&&1\\0&&0\\[-1ex]\vdots&&\vdots\\0&&0\\1&&2\\0&&0\\2&&1\\&0&\end{array},\qquad P_1=\begin{array}{c}1\\0\\2\\0\\1\\[-1ex]\vdots\\0\\2\\0\\1\end{array}, \qquad P_2=\begin{array}{c}2\\0\\1\\0\\2\\[-1ex]\vdots\\0\\1\\0\\2\end{array}.$$
\hrule
\end{figure}
where the radical series length of each of these modules is $2^d+1$. 
The decomposition matrix of $B$ is given in Figure \ref{fig:decompsl1}.
\begin{figure}[ht]  \hrule \caption{\label{fig:decompsl1} The decomposition matrix for blocks as in \S$\ref{ss:psl1}$.}
$$\begin{array}{ccc}
&\begin{array}{c@{}c@{}c}\varphi_0\,&\,\varphi_1\,&\,\varphi_2\end{array}\\[1ex]
\begin{array}{c}\chi_1\\ \chi_2\\ \chi_3 \\ \chi_4\\ \chi_{5,i}\end{array} &
\left[\begin{array}{ccc}1&0&0\\1&1&0\\1&0&1\\1&1&1\\2&1&1\end{array}\right]
 &\begin{array}{c}\\ \\ \\ \\ 1\le i\le 2^{d-2}-1.\end{array}\end{array}$$
\vspace{1ex}
\hrule
\end{figure}

\subsection{The principal $2$-modular block of $\mathrm{PSL}_2(\mathbb{F}_q)$ when $q\equiv 3 \mod 4$}
\label{ss:psl2}
\setcounter{equation}{0}
\setcounter{figure}{0}

Let $G$ be a finite group, and let $B$ be a block of $kG$ which is Morita equivalent to the principal block of $k\mathrm{PSL}_2(\mathbb{F}_q)$ where $q\equiv 3 \mod 4$. Suppose that $2^d$ is the order of the defect groups of $B$. Then, by \cite{erd}, $B$ is Morita equivalent to the special biserial algebra $\Lambda=kQ/I$ where $Q$ is given in Figure \ref{fig:psl2A} and
$$I=\langle \delta\beta,\lambda\delta,\beta\gamma,\kappa\gamma,\eta\kappa, 
\gamma\eta,\gamma\beta-\lambda\kappa,\kappa\lambda-(\delta\eta)^{2^{d-2}},
(\eta\delta)^{2^{d-2}} - \beta\gamma\rangle.$$
\begin{figure}[ht] \hrule \caption{\label{fig:psl2A} The quiver $Q$ for blocks as in \S$\ref{ss:psl2}$.}
$$Q=\vcenter{\xymatrix  {
 0\,\bullet \ar@<.7ex>[rr]^{\beta} \ar@<.8ex>[rr];[]^{\gamma}\ar@<.7ex>[rdd]^{\kappa} \ar@<.8ex>[rdd];[]^{\lambda}
&&\bullet\ar@<.7ex>[ldd]^{\delta} \ar@<.8ex>[ldd];[]^{\eta}\,1\\&&\\ &
\unter{\mbox{\normalsize $\bullet$}}{\mbox{\normalsize $2$}}& }}$$
\vspace{1ex}
\hrule
\end{figure}
We denote the irreducible $\Lambda$-modules by $S_0,S_1,S_2$, or, using short-hand, by $0,1,2$.
The radical series of the projective indecomposable $\Lambda$-modules (and hence of the projective indecomposable $B$-modules) are described in  Figure \ref{fig:psl2B}
\begin{figure}[ht] \hrule \caption{\label{fig:psl2B} The radical series of the projective indecomposable modules for blocks  as in \S$\ref{ss:psl2}$.}
$$P_0=\begin{array}{cc}\multicolumn{2}{c}{0}\\
1&2\\\multicolumn{2}{c}{0}\end{array},\qquad P_1=\begin{array}{cc}\multicolumn{2}{c}{1}\\0&2\\&1\\&2\\[-1ex]&\vdots\\&1\\&2\\\multicolumn{2}{c}{1}\end{array}, \qquad 
P_2=\begin{array}{cc}\multicolumn{2}{c}{2}\\0&1\\&2\\&1\\[-1ex]&\vdots\\&2\\&1\\\multicolumn{2}{c}{2}\end{array}.$$
\hrule
\end{figure}
where for $i\in\{1,2\}$, $\mathrm{rad}(P_i)/\mathrm{soc}(P_i)$ is isomorphic to the direct sum of $S_0$ and a uniserial module of length $2^{d-1}-1$.
The decomposition matrix of $B$ is given in Figure \ref{fig:decompsl2}.
\begin{figure}[ht]   \hrule\vspace{1ex}\caption{\label{fig:decompsl2}The decomposition matrix for blocks  as in \S$\ref{ss:psl2}$.}
$$\begin{array}{ccc}
&\begin{array}{c@{}c@{}c}\varphi_0\,&\,\varphi_1\,&\,\varphi_2\end{array}\\[1ex]
\begin{array}{c}\chi_1\\ \chi_2\\ \chi_3 \\ \chi_4\\ \chi_{5,i}\end{array} &
\left[\begin{array}{ccc}1&0&0\\0&1&0\\0&0&1\\1&1&1\\0&1&1\end{array}\right]
&\begin{array}{c}\\ \\ \\ \\ 1\le i\le 2^{d-2}-1.\end{array}
\end{array}$$
\vspace{1ex}
\hrule
\end{figure}

\subsection{The principal $2$-modular block of $A_7$}
\label{ss:a7}
\setcounter{equation}{0}
\setcounter{figure}{0}

\begin{figure}[ht] \hrule \caption{\label{fig:a7A} The quiver $Q$ for blocks as in \S$\ref{ss:a7}$.}
$$\xymatrix @R=-.2pc {
&1&0&\\
 Q= \quad&\ar@(ul,dl)_{\alpha} \bullet \ar@<.8ex>[r]^{\beta} \ar@<.9ex>[r];[]^{\gamma}
& \bullet \ar@<.8ex>[r]^(.46){\delta} \ar@<.9ex>[r];[]^(.54){\eta} & \bullet\;2}$$
\vspace{1ex}
\hrule
\end{figure}
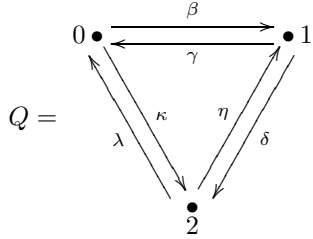

Let $G$ be a finite group, and let $B$ be a block of $kG$ 
which is Morita equivalent to the principal block of $kA_7$. Suppose that $2^d$ is the order of the defect groups of $B$. Then $d=3$, and by \cite{erd},
$B$ is Morita equivalent to the special biserial algebra $\Lambda=kQ/I$ where $Q$ is given in Figure \ref{fig:a7A} and
$$I=\langle\beta\alpha,\alpha\gamma,\gamma\beta,\delta\eta,\eta\delta\beta\gamma-\beta\gamma\eta\delta,
\alpha^2
-\gamma\eta\delta\beta\rangle.$$
We denote the irreducible $\Lambda$-modules by $S_0,S_1,S_2$, or, using short-hand, by $0,1,2$.
The radical series of the projective indecomposable $\Lambda$-modules (and hence of the projective indecomposable $B$-modules) are described in  Figure \ref{fig:a7B} 
\begin{figure}[ht] \hrule \caption{\label{fig:a7B} The radical series of the projective indecomposable modules for blocks  as in \S$\ref{ss:a7}$.}
$$P_0=\begin{array}{cc}\multicolumn{2}{c}{0}\\
1&2\\0&0\\2&1\\\multicolumn{2}{c}{0}\end{array},\qquad P_1=\begin{array}{cc}\multicolumn{2}{c}{1}\\1&0\\ 
&2\\&0\\
\multicolumn{2}{c}{1}\end{array}, \qquad 
P_2=\begin{array}{c}2\\0\\1\\0\\2\end{array}.$$
\hrule
\end{figure}
where $\mathrm{rad}(P_1)/\mathrm{soc}(P_1)$ is isomorphic to the direct sum of  $S_1$ and a uniserial module with radical series $(S_0,S_2,S_0)$. 
The decomposition matrix of $B$ is given in Figure \ref{fig:decoma7}.
\begin{figure}[ht]  \hrule\vspace{1ex}\caption{\label{fig:decoma7}The decomposition matrix for blocks  as in \S$\ref{ss:a7}$.}
$$\begin{array}{ccc}
&\begin{array}{c@{}c@{}c}\varphi_0\,&\,\varphi_1\,&\,\varphi_2\end{array}\\[1ex]
\begin{array}{c}\chi_1\\ \chi_2\\ \chi_3 \\ \chi_4\\ \chi_{5,i}\end{array} &
\left[\begin{array}{ccc}1&0&0\\1&1&0\\1&0&1\\1&1&1\\0&1&0\end{array}\right]
&\begin{array}{c}\\ \\ \\ \\ i=1.\end{array}
\end{array}$$
\vspace{1ex}
\hrule
\end{figure}

\subsection{Ordinary characters for blocks with dihedral defect groups}
\label{ss:ordinarydihedral}
\setcounter{equation}{0}
\setcounter{figure}{0}

Let $G$ be a finite group and let $B$ be a block of $kG$ with dihedral defect group $D$ of order $2^d$ where $d\ge 3$. Moreover, assume that $B$ contains exactly three isomorphism classes of simple $kG$-modules. This means that in the notation of \cite[\S4]{brauer2} we are in Case $(aa)$ (see \cite[Thm. 2]{brauer2}).

Let $F$ be the fraction field of $W$, and let $\zeta_{2^\ell}$ be a fixed primitive $2^\ell$-th root of unity in an algebraic closure of $F$ for $2\le\ell\le d-1$. Let
$$\chi_1,\chi_2,\chi_3,\chi_4,\qquad \chi_{5,i}, 1\le i\le 2^{d-2}-1,$$
be the ordinary irreducible characters of $G$ belonging to $B$. Let $\sigma$ be an element of order $2^{d-1}$ in $D$. By \cite{brauer2}, there is a block $b_\sigma$ of $kC_G(\sigma)$ with $b_\sigma^G=B$ which contains a unique $2$-modular character $\varphi^{(\sigma)}$ such that the following is true. There is an ordering of $(1,2,\ldots,2^{d-2}-1)$ such that for $1\le i\le 2^{d-2}-1$ and $r$ odd,
\begin{equation}
\label{eq:great1}
\chi_{5,i}(\sigma^r)=(\zeta_{2^{d-1}}^{ri}+\zeta_{2^{d-1}}^{-ri})\cdot \varphi^{(\sigma)}(1).
\end{equation}

Note that $W$ contains all roots of unity of order not divisible by $2$. Hence by \cite{brauer2} and by \cite{fong}, the characters 
$\chi_1,\chi_2,\chi_3,\chi_4$
correspond to simple $FG$-modules. On the other hand, the characters 
$\chi_{5,i}$, $i=1,\ldots,2^{d-2}-1$,
fall into $d-2$ Galois orbits $\mathcal{O}_2,\ldots, \mathcal{O}_{d-1}$ under the action of $\mathrm{Gal}(F(\zeta_{2^{d-1}}+\zeta_{2^{d-1}}^{-1})/F)$. Namely for $2\le\ell\le d-1$, $\mathcal{O}_{\ell}=\{ \chi_{5,2^{d-1-\ell}(2u-1)} \;|\; 1\le u\le 2^{\ell-2}\}$. The field generated by the character values of each $\xi_\ell\in\mathcal{O}_\ell$ over $F$ is $F(\zeta_{2^\ell}+\zeta_{2^\ell}^{-1})$. Hence by \cite{fong}, each $\xi_\ell$ corresponds to an absolutely irreducible $F(\zeta_{2^\ell}+\zeta_{2^\ell}^{-1})G$-module $X_\ell$. 
By \cite[Satz V.14.9]{hup},  this implies that for $2\le\ell\le d-1$, the Schur index of each $\xi_\ell\in\mathcal{O}_\ell$ over $F$ is $1$. Hence we obtain $d-2$ non-isomorphic simple $FG$-modules $V_2,\ldots,V_{d-1}$
with characters $\rho_2,\ldots, \rho_{d-1}$ satisfying
\begin{equation}
\label{eq:goodchar1}
\rho_\ell =\sum_{\xi_\ell\in\mathcal{O}_\ell}\xi_\ell = 
\sum_{u=1}^{2^{\ell-2}} \chi_{5,2^{d-1-\ell}(2u-1)} 
\qquad\mbox{for $2\le \ell \le d-1$.}
\end{equation}
By \cite[Hilfssatz V.14.7]{hup}, $\mathrm{End}_{FG}(V_\ell)$ is a commutative $F$-algebra isomorphic to the field generated over $F$ by the character values of any $\xi_\ell\in\mathcal{O}_\ell$. This means
\begin{equation}
\label{eq:goodendos}
\mathrm{End}_{FG}(V_\ell)\cong F(\zeta_{2^\ell}+\zeta_{2^\ell}^{-1})\qquad\mbox{for $2\le \ell \le d-1$.}
\end{equation}

By \cite{brauer2}, the characters $\chi_{5,i}$ have the same degree $x$ for $1\le i\le 2^{d-2}-1$. The characters $\chi_1,\chi_2,\chi_3,\chi_4$ have height $0$ and $\chi_{5,i}$, $1\le i\le 2^{d-2}-1$, have height $1$. Hence $x=2^{a-d+1}x^*$ where $\#G=2^a\cdot g^*$ and $x^*$ and $g^*$ are odd. Since the centralizer $C_G(\sigma)$ contains $\langle \sigma \rangle$, we have $\# C_G(\sigma)=2^{d-1}\cdot 2^b\cdot m^*$ where $b\ge 0$ and $m^*$ is odd. Suppose $\varphi^{(\sigma)}(1)=2^c\cdot n^*$ where $c\ge 0$ and $n^*$ is odd. Note that if $\psi$ is an ordinary irreducible character of $C_G(\sigma)$ belonging to the block $b_\sigma$, then by \cite[p. 61]{serre}, $\psi(1)$ divides $(\# C_G(\sigma))/(\#\langle\sigma\rangle)=2^b\cdot m^*$. Because $\psi(1)=s_\psi\cdot \varphi^{(\sigma)}(1)$ for some positive integer $s_\psi$, we have $c\le b$. 

Let $C$ be the conjugacy class in $G$ of $\sigma$, and let $t(C)\in WG$ be the class sum of $C$. We want to determine the action of $t(C)$ on $V_\ell$ for $2\le \ell\le d-1$. For this, we identify $\mathrm{End}_{FG}(V_\ell)\cong F(\zeta_{2^\ell}+\zeta_{2^\ell}^{-1})$ with $\mathrm{End}_{F(\zeta_{2^\ell}+\zeta_{2^\ell}^{-1})G}(X_\ell)$ for one particular absolutely irreducible $F(\zeta_{2^\ell}+\zeta_{2^\ell}^{-1})G$-constituent $X_\ell$ of $V_\ell$ with character $\xi_\ell$. By (\ref{eq:goodchar1}), we can choose $\xi_\ell=\chi_{5,2^{d-1-\ell}}$. Then, under this identification, for $2\le \ell \le d-1$, the action of $t(C)$ on $V_\ell$ is given as multiplication by 
\begin{eqnarray}
\label{eq:longone}
\frac{\# C}{\xi_\ell(1)} \cdot \xi_\ell(\sigma)&=&\frac{\# C}{\xi_\ell(1)} \cdot\varphi^{(\sigma)}(1)\cdot (\zeta_{2^{d-1}}^{2^{d-1-\ell}}+\zeta_{2^{d-1}}^{-2^{d-1-\ell}})\\
&=& \frac{[G:C_G(\sigma)]}{x} \cdot\varphi^{(\sigma)}(1)\cdot (\zeta_{2^{d-1}}^{2^{d-1-\ell}}+\zeta_{2^{d-1}}^{-2^{d-1-\ell}})\nonumber\\
&=& 2^{c-b} \frac{g^*\cdot n^*}{m^*\cdot x^*}\cdot (\zeta_{2^{d-1}}^{2^{d-1-\ell}}+\zeta_{2^{d-1}}^{-2^{d-1-\ell}})\nonumber
\end{eqnarray}
where, as shown above, $c\le b$. Note that $\frac{g^*\cdot n^*}{m^*\cdot x^*}$ is a unit in $W$, since $g^*\cdot n^*$ and $m^*\cdot x^*$ are odd. Since $t(C)\in WG$, we must have $c\ge b$, i.e. $c=b$.
Therefore, (\ref{eq:longone}) implies that there exists a unit $\omega$ in $W$ such that for $2\le \ell \le d-1$, the action of $t(C)$ on $V_\ell$ is given as multiplication by
\begin{equation}
\label{eq:thatsit}
\omega\cdot (\zeta_{2^{d-1}}^{2^{d-1-\ell}}+\zeta_{2^{d-1}}^{-2^{d-1-\ell}})
\end{equation}
when we identify $\mathrm{End}_{FG}(V_\ell)$ with $\mathrm{End}_{F(\zeta_{2^\ell}+\zeta_{2^\ell}^{-1})G}(X_\ell)$ for an absolutely irreducible $F(\zeta_{2^\ell}+\zeta_{2^\ell}^{-1})G$-constituent $X_\ell$ of $V_\ell$ with character $\chi_{5,2^{d-1-\ell}}$.


\section{Universal deformation rings modulo $2$}
\label{s:defmod2}
\setcounter{equation}{0}
\setcounter{figure}{0}

As in \S\ref{s:dihedralsylow}, let $k$ be an algebraically closed field of characteristic $p=2$, let $G$ be a finite group, and let $B$ be a block of $kG$ with dihedral defect groups containing precisely three isomorphism classes of simple $kG$-modules. Suppose $2^d$ is the order of the defect groups of $B$. In this section, we determine the universal deformation ring modulo $2$ for all finitely generated $B$-modules with stable endomorphism ring $k$. Since the case $d=2$ has been done in \cite{bl}, we assume throughout this section that $d\ge 3$.

In \S \ref{ss:defmod2psl1}, we first look at blocks $B$ that are Morita equivalent to blocks as in \S\ref{ss:psl1}. 
In \S \ref{ss:defmod2rest}, we then show how stable equivalences of Morita type can be used to get analogous results for blocks $B$ that are Morita equivalent to blocks as in \S\ref{ss:psl2} and \S\ref{ss:a7}.

\subsection{Universal deformation rings modulo $2$ for blocks as in $\S\ref{ss:psl1}$}
\label{ss:defmod2psl1}
\setcounter{equation}{0}
\setcounter{figure}{0}

The objective of this subsection is to prove the following result:

\begin{subprop}
\label{prop:psl1}
Let $B$ be a block of $kG$ which is Morita equivalent to $\Lambda=kQ/I$ with $Q$ and $I$ as in $\S\ref{ss:psl1}$. Suppose $2^d$ is the order of the defect groups of $B$ with $d\ge 3$. Denote the three simple $B$-modules by $T_0$, $T_1$ and $T_2$, where $T_i$ corresponds to $S_i$, for $i\in\{0,1,2\}$, under the Morita equivalence.
Let $\mathfrak{C}$ be a component of the stable Auslander-Reiten quiver of $B$ containing a module with endomorphism ring $k$. Then $\mathfrak{C}$ contains a simple module, or a uniserial module of length $4$.
\begin{enumerate}
\item[i.] Suppose $\mathfrak{C}$ contains $T_0$. Then $\mathfrak{C}$ and $\Omega(\mathfrak{C})$ are both of type $\mathbb{Z}A_\infty^\infty$. All modules $M$ in $\mathfrak{C}\cup\Omega(\mathfrak{C})$ have stable endomorphism ring equal to $k$ and $R(G,M)/2R(G,M)\cong k$.
\item[ii.] Let $i\in\{1,2\}$, and suppose $\mathfrak{C}$ contains $T_i$. Then $\mathfrak{C}$ is a $3$-tube with $T_i$ belonging to its boundary, and $\mathfrak{C}=\Omega(\mathfrak{C})$. If $M$ lies in $\mathfrak{C}$ having stable endomorphism ring $k$, then $M\in\{T_i,\Omega^2(T_i),\Omega^4(T_i)\}$ up to isomorphism, and $R(G,M)/2R(G,M)\cong k$.
\item[iii.] Suppose $\mathfrak{C}$ contains a uniserial module $Y$ of length $4$. Then $\mathfrak{C}$ and $\Omega(\mathfrak{C})$ are both of type $\mathbb{Z}A_\infty^\infty$, and $\mathfrak{C}=\Omega(\mathfrak{C})$ exactly when $d=3$. If $M$ lies in $\mathfrak{C}\cup\Omega(\mathfrak{C})$ having stable endomorphism ring $k$, then $M$ is isomorphic to $\Omega^j(Y)$ for some integer $j$, and $R(G,M)/2R(G,M)\cong k[t]/(t^{2^{d-2}})$. 
\end{enumerate}
The only components of the stable Auslander-Reiten quiver of $B$ containing modules with stable endomorphism ring $k$ are the ones in $(i)-(iii)$.
\end{subprop}

\begin{subrem}
\label{rem:psl1uni}
If $B$ is as in Proposition $\ref{prop:psl1}$, then there are precisely four uniserial $B$-modules of length $4$:
$$Y_1=\begin{array}{c}T_1\\T_0\\T_2\\T_0\end{array}, \;
\Omega^2(Y_1)= \begin{array}{c}T_0\\T_2\\T_0\\T_1\end{array}, \;
Y_2=\begin{array}{c}T_2\\T_0\\T_1\\T_0\end{array},\;
\Omega^2(Y_2)= \begin{array}{c}T_0\\T_1\\T_0\\T_2\end{array}.$$
\end{subrem}

To prove Proposition \ref{prop:psl1}, we need several Lemmas. 

\begin{sublem}
\label{lem:psl1}
Let $\Lambda=kQ/I$ with $Q$ and $I$ as in $\S\ref{ss:psl1}$. 
Let $\mathfrak{C}_0$ be the component of the stable Auslander-Reiten quiver of $\Lambda$ containing $S_0$.
Let $M$ be an indecomposable $\Lambda$-module with $\mathrm{End}_\Lambda(M)=k$. Then $M$ either lies in $\mathfrak{C}_0$, or $M$ is isomorphic to $S_1$ or $S_2$, or $M$ is uniserial of length $4$.
\end{sublem}

\begin{proof}
Suppose first that $S_0$ is a direct summand of $\mathrm{top}(M)$. Then $S_0$ cannot be a direct summand of $\mathrm{soc}(M)$. 
The modules
$S_0$, $\begin{array}{c}0\\1\end{array}$, $\begin{array}{c}0\\2\end{array}$,
$\begin{array}{c}0\\1\\0\\2\end{array}$, $\begin{array}{c}0\\2\\0\\1\end{array}$ and $\begin{array}{cc}\multicolumn{2}{c}{0}\\1&2\end{array}$ have endomorphism ring equal to $k$. If $M$ is uniserial of length at least 5 or $M$ is not uniserial of length at least 4, then there is an endomorphism of $M$ factoring non-trivially through either $S_0$, $\begin{array}{c}0\\1\end{array}$ or $\begin{array}{c}0\\2\end{array}$, which is a contradiction.

Similarly, we see that if $S_0$ is a direct summand of $\mathrm{soc}(M)$, then $S_0$ cannot be a direct summand of $\mathrm{top}(M)$ and $M$ must be isomorphic to $S_0$, $\begin{array}{c}1\\0\end{array}$, $\begin{array}{c}2\\0\end{array}$,
$\begin{array}{c}2\\0\\1\\0\end{array}$, $\begin{array}{c}1\\0\\2\\0\end{array}$ or $\begin{array}{cc}1&2\\\multicolumn{2}{c}{0}\end{array}$.

If $S_0$ is neither a direct summand of $\mathrm{top}(M)$ nor of $\mathrm{soc}(M)$, then $M$ can only be isomorphic to $S_1$, $S_2$, $\begin{array}{c}1\\0\\2\end{array}$ or $\begin{array}{c}2\\0\\1\end{array}$, since $P_1$ and $P_2$ are both uniserial. 

Of all the above possibilities for $M$, only $S_1$, $S_2$, and the uniserial modules of length $4$ do not lie in $\mathfrak{C}_0$.
\end{proof}

The proofs of the following three Lemmas are rather technical and will be deferred to
\S\ref{s:stableend}. Let $\Lambda=kQ/I$ with $Q$ and $I$ as in $\S\ref{ss:psl1}$. 

\begin{sublem}
\label{lem:psl1endos0}
Let $\mathfrak{C}_0$ be the component of the stable Auslander-Reiten quiver of $\Lambda$ containing $S_0$, and let $M$ be a $\Lambda$-module belonging to $\mathfrak{C}_0\cup\Omega(\mathfrak{C}_0)$. Then $\mathfrak{C}_0$ and $\Omega(\mathfrak{C}_0)$ are both of type $\mathbb{Z}A_\infty^\infty$, and $\underline{\mathrm{End}}_\Lambda(M)=k$ and $\mathrm{Ext}^1_\Lambda(M,M)=0$. 
\end{sublem}

\begin{sublem}
\label{lem:psl1ci}
Let $i\in\{1,2\}$ and let $\mathfrak{C}_i$ be the component of the stable Auslander-Reiten quiver of $\Lambda$ containing $S_i$. Then $\mathfrak{C}_i$ is a $3$-tube with $S_i$ belonging to its boundary, and $\mathfrak{C}_i=\Omega(\mathfrak{C}_i)$. If $M$ belongs to $\mathfrak{C}_i$ and has stable endomorphism ring $k$, then $M$ is in the $\Omega$-orbit of $S_i$ and $\mathrm{Ext}^1_\Lambda(M,M)=0$.
\end{sublem}

\begin{sublem}
\label{lem:psl1uni}
Let $\mathfrak{C}$ be a component of the stable Auslander-Reiten quiver of $\Lambda$ containing a uniserial module $X$ of length $4$. Then $\mathfrak{C}$ and $\Omega(\mathfrak{C})$ are both of type $\mathbb{Z}A_\infty^\infty$, and $\mathfrak{C}=\Omega(\mathfrak{C})$ exactly when $d=3$. If $M$ belongs to $\mathfrak{C}\cup\Omega(\mathfrak{C})$ and has stable endomorphism ring $k$, then $M$ is in the $\Omega$-orbit of $X$ and $\mathrm{Ext}^1_\Lambda(M,M)=k$.
\end{sublem}

\medskip

\noindent
\textit{Proof of Proposition $\ref{prop:psl1}$.}
The last statement of Proposition \ref{prop:psl1} will also be proved in \S\ref{s:stableend}. We now prove the remaining statements. 
Let $\mathfrak{C}$ be a component of the stable Auslander-Reiten quiver of $B$ containing a module with endomorphism ring $k$. By Lemma \ref{lem:psl1}, $\mathfrak{C}$ contains a simple module or a uniserial module of length $4$.
Part (i) (resp. part (ii)) of Proposition \ref{prop:psl1} follows from Lemma \ref{lem:psl1endos0} (resp. Lemma \ref{lem:psl1ci}). 
Because of Lemma \ref{lem:psl1uni}, Remark \ref{rem:psl1uni} and Lemma \ref{lem:defhelp}, to prove part (iii) we only need to show 
$R(G,M)/2R(G,M)\cong k[t]/(t^{2^{d-2}})$ if $M$ is either $Y_1$ or $Y_2$.
We show this for $M=Y_1$. (The case $M=Y_2$ is proved similarly.) Let $\overline{R}= R(G,Y_1)/2R(G,Y_1)$. By Lemmas \ref{lem:psl1} and \ref{lem:psl1uni}, we have $\mathrm{End}_{kG}(Y_1)=k$ and $\mathrm{Ext}^1_{kG}(Y_1,Y_1)=k$. The projective indecomposable $kG$-module $P_{T_1}$ has the form $\begin{array}{c}T_1\\U\\T_1\end{array}$ where $U$ is uniserial of length $4\cdot 2^{d-2}-1$ with descending radical series
$$(T_0,T_2,T_0,T_1,T_0,T_2,T_0,\ldots,T_1,T_0,T_2,T_0).$$
If $\overline{U}=\begin{array}{c}T_1\\U\end{array}$, then 
$$\mathrm{Ext}^1_{kG}(\overline{U},Y_1)=
\underline{\mathrm{Hom}}_{kG}(\Omega(\overline{U}),Y_1)=
\underline{\mathrm{Hom}}_{kG}(T_1,Y_1)=0.$$
By Lemma \ref{lem:trythis}, this implies $\overline{R}\cong k[t]/(t^{2^{d-2}})$, and the universal mod $2$ deformation of $Y$ over $\overline{R}$ is represented by the $kG$-module $\overline{U}$. 

\begin{subrem}
\label{rem:uniliftpsl1}
It follows from the proof of Proposition \ref{prop:psl1} that if $Y=\begin{array}{c}T_i\\T_0\\T_j\\T_0\end{array}$ where $i\neq j$ in $\{1,2\}$, then the universal mod $2$ deformation of $Y$ is represented by the uniserial $kG$-module $\overline{U}$ of length $4\cdot 2^{d-2}$ with descending radical series
$$(T_i,T_0,T_j,T_0,T_i,T_0,T_j,T_0,\ldots,T_i,T_0,T_j,T_0).$$
Moreover the uniserial $kG$-module $\overline{U'}\cong \overline{U}/Y$ of length $4\cdot (2^{d-2}-1)$ defines a lift of $Y$ over $k[t]/(t^{2^{d-2}-1})$.
\end{subrem}

\subsection{Universal deformation rings modulo $2$ for blocks as in $\S\ref{ss:psl2}$ and $\S\ref{ss:a7}$}
\label{ss:defmod2rest}
\setcounter{equation}{0}
\setcounter{figure}{0}

In this subsection we prove analogous results to Proposition \ref{prop:psl1} for blocks  that are Morita equivalent to blocks as in \S\ref{ss:psl2} and \S\ref{ss:a7}. We start with  blocks as in \S\ref{ss:psl2}.
\begin{subprop}
\label{prop:psl2}
Let $B$ be a block of $kG$ which is Morita equivalent to $\Lambda=kQ/I$ with $Q$ and $I$ as in $\S\ref{ss:psl2}$.  Suppose $2^d$ is the order of the defect groups of $B$ with $d\ge 3$. Denote the three simple $B$-modules by $T_0$, $T_1$ and $T_2$, where $T_i$ corresponds to $S_i$, for $i\in\{0,1,2\}$, under the Morita equivalence.
Let $\mathfrak{C}$ be a component of the stable Auslander-Reiten quiver of $B$ containing a module with endomorphism ring $k$. Then $\mathfrak{C}$ contains a simple module, or a uniserial module of length $2$.
\begin{enumerate}
\item[i.] Suppose $\mathfrak{C}$ contains $T_0$. Then $\Omega(\mathfrak{C})$ contains $T_1$ and $T_2$, and $\mathfrak{C}$ and $\Omega(\mathfrak{C})$ are both of type $\mathbb{Z}A_\infty^\infty$. All modules $M$ in $\mathfrak{C}\cup\Omega(\mathfrak{C})$ have stable endomorphism ring equal to $k$ and $R(G,M)/2R(G,M)\cong k$.
\item[ii.] Let $i\neq j$ be in $\{1,2\}$, and suppose $\mathfrak{C}$ contains $T_{0,i}=\begin{array}{c}T_0\\T_i\end{array}$. Then $\mathfrak{C}$ is a $3$-tube with $T_{0,i}$ and $T_{j,0}=\begin{array}{c}T_j\\T_0\end{array}$ belonging to its boundary, and $\mathfrak{C}=\Omega(\mathfrak{C})$. If $M$ lies in $\mathfrak{C}$ having stable endomorphism ring $k$, then $M\in\{T_{0,i},\Omega^2(T_{0,i}),\Omega^4(T_{0,i})\}$ up to isomorphism, and $R(G,M)/2R(G,M)\cong k$.
\item[iii.] Suppose $Y\in\left\{\begin{array}{c}T_1\\T_2\end{array},\begin{array}{c}T_2\\T_1\end{array}\right\}$ and $\mathfrak{C}$ contains $Y$. Then $\mathfrak{C}$ and $\Omega(\mathfrak{C})$ are both of type $\mathbb{Z}A_\infty^\infty$, and $\mathfrak{C}=\Omega(\mathfrak{C})$ exactly when $d=3$. If $M$ lies in $\mathfrak{C}\cup\Omega(\mathfrak{C})$ having stable endomorphism ring $k$, then $M$ is isomorphic to $\Omega^j(Y)$ for some integer $j$, and $R(G,M)/2R(G,M)\cong k[t]/(t^{2^{d-2}})$. 
\end{enumerate}
The only components of the stable Auslander-Reiten quiver of $B$ containing modules with stable endomorphism ring $k$ are the ones in $(i)-(iii)$.
\end{subprop}

To prove Proposition \ref{prop:psl2} we need the following result which is proved similarly to Lemma \ref{lem:psl1}.
\begin{sublem}
\label{lem:psl2}
Let $\Lambda=kQ/I$ with $Q$ and $I$ as in $\S\ref{ss:psl2}$. 
For $i\in\{0,1,2\}$, let $\mathfrak{C}_i$ be the component of the stable Auslander-Reiten quiver of $\Lambda$ containing $S_i$.
Let $M$ be an indecomposable $\Lambda$-module with $\mathrm{End}_\Lambda(M)=k$. Then $M$ either lies in $\mathfrak{C}_i$ for some $i$, or $M$ is uniserial of length $2$.
\end{sublem}

\medskip

\noindent
\textit{Proof of Proposition $\ref{prop:psl2}$.}
Let $\mathfrak{C}$ be a component of the stable Auslander-Reiten quiver of $B$ containing a module with endomorphism ring $k$. By Lemma \ref{lem:psl2}, $\mathfrak{C}$ contains a simple module or a uniserial module of length $2$.

Let $B_0$ be the principal block of $k\mathrm{PSL}_2(\mathbb{F}_q)$ where $q\equiv 1\mod 4$ and the Sylow $2$-subgroups of $\mathrm{PSL}_2(\mathbb{F}_q)$ have order $2^d$.
By \cite{linckel}, $B$ and $B_0$ are derived equivalent. By \cite[Cor. 5.5]{rickard1}, this means that there is a stable equivalence of Morita type between $B$ and $B_0$ given by a $B$-$B_0$-bimodule $\Xi$. By Lemma \ref{lem:weakstabmordef} and Remark \ref{rem:stabmor}, if $V$ is a finitely generated $B_0$-module with stable endomorphism ring $k$ and $V'=\Xi\otimes_{B_0}V$, then $V'$ has stable endomorphism ring $k$ and $R(\mathrm{PSL}_2(\mathbb{F}_q),V)/2R(\mathrm{PSL}_2(\mathbb{F}_q),V)\cong R(G,V')/2R(G,V')$. Moreover, $V'\cong V''\oplus P$ as $kG$-modules where $P$ is projective, $V''$ is indecomposable and $R(G,V'')\cong R(G,V')$.

Because of Remark \ref{rem:stablegood}, to complete the proof of Proposition \ref{prop:psl2}, we need to find the components of the stable Auslander-Reiten quiver of $B_0$ and of $B$, respectively, that correspond to each other under the functor $\Xi\otimes_{B_0}-$, and we need to match up certain modules in these components. Note that by Remark \ref{rem:stablegood}, $\Xi\otimes_{B_0}-$ commutes with the Heller operator $\Omega$. Let $\mathfrak{C}$ be a component of the stable Auslander-Reiten quiver of $B$ containing a module with endomorphism ring $k$.

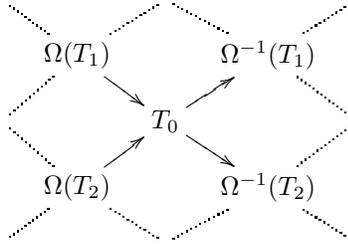
\begin{figure}[ht] \hrule \caption{\label{fig:ct0} The stable Auslander-Reiten component near $T_0$.}
$$\xymatrix @-1.2pc{
&&&&\\
&\Omega(T_1)\ar[rd]\ar@{.}[ld]\ar@{.}[lu]\ar@{.}[ru]&&\Omega^{-1}(T_1)\ar@{.}[rd]\ar@{.}[lu]\ar@{.}[ru]&\\
&&T_0\ar[ru]\ar[rd]&&\\
&\Omega(T_2)\ar[ru]\ar@{.}[rd]\ar@{.}[ld]\ar@{.}[lu]&&\Omega^{-1}(T_2)\ar@{.}[rd]\ar@{.}[ld]\ar@{.}[ru]&\\
&&&&
}$$
\hrule
\end{figure}

Suppose first $\mathfrak{C}$ contains $T_0$. Then $\mathfrak{C}$ is of type $\mathbb{Z}A_\infty^\infty$. Near $T_0$, $\mathfrak{C}$ looks 
as in Figure \ref{fig:ct0}.
Hence $\Omega(\mathfrak{C})$ contains $T_1$ and $T_2$.
By Proposition \ref{prop:psl1}, it follows that $\mathfrak{C}$ and $\Omega(\mathfrak{C})$ correspond to the components of the stable Auslander-Reiten quiver of $B_0$ as in Proposition \ref{prop:psl1}(i). This proves part (i) of Proposition \ref{prop:psl2}.

Let now $i\ne j$ be in $\{1,2\}$, and suppose $\mathfrak{C}$ contains $T_{0,i}$. Then $\mathfrak{C}$ is a $3$-tube and $T_{0,i}$ belongs to its boundary. Since  $\Omega^2(T_{j,0})=T_{0,i}$ and $\Omega(T_{0,i})=\Omega^{-2}(T_{0,i})$,  $\mathfrak{C}$ also contains $T_{j,0}$ and $\mathfrak{C}=\Omega(\mathfrak{C})$.
By Proposition \ref{prop:psl1}, it follows that for $i=1,2$, the components $\mathfrak{C}$ correspond to the components of the stable Auslander-Reiten quiver of $B_0$ as in Proposition \ref{prop:psl1}(ii). This proves part (ii) of Proposition \ref{prop:psl2}.

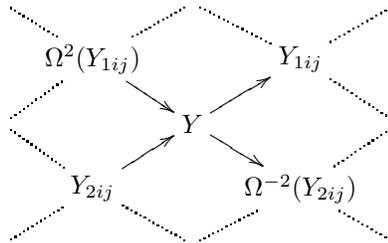
\begin{figure}[ht] \hrule \caption{\label{fig:ctij} The stable Auslander-Reiten component near $Y$.}
$$\xymatrix @-1.2pc{
&&&&\\
&\Omega^2(Y_{1ij})\ar[rd]\ar@{.}[ld]\ar@{.}[lu]\ar@{.}[ru]&&Y_{1ij}\ar@{.}[rd]\ar@{.}[lu]\ar@{.}[ru]&\\
&&Y\ar[ru]\ar[rd]&&\\
&Y_{2ij}\ar[ru]\ar@{.}[rd]\ar@{.}[ld]\ar@{.}[lu]&&\Omega^{-2}(Y_{2ij})\ar@{.}[rd]\ar@{.}[ld]\ar@{.}[ru]&\\
&&&&
}$$
\hrule
\end{figure}

Finally, let $i\neq j$ in $\{1,2\}$, and suppose $\mathfrak{C}$ contains $Y=\begin{array}{c}T_i\\T_j\end{array}$. Then $\mathfrak{C}$ is a component of type $\mathbb{Z}A_\infty^\infty$ and $Y$ is a module of minimal length in $\mathfrak{C}$. Near $Y$, $\mathfrak{C}$ looks 
as in Figure \ref{fig:ctij}.
where $Y_{1ij}$ and $Y_{2ij}$ correspond to string modules of $\Lambda=kQ/I$ with $Q$ and $I$ as in $\S\ref{ss:psl2}$ as follows:
$$Y_{1ij}=\begin{array}{c@{}c@{}c@{}c}T_i&&T_0&\\&T_j&&T_i\end{array},\qquad
Y_{2ij}=\begin{array}{c@{}c@{}c@{}c}T_i&&&\\&T_j&&T_0\\&&T_i\end{array}.$$
For all $d\ge 3$, $Y_{1ij}$ has a non-trivial endomorphism factoring through $T_i$ which does not factor through a projective $B$-module, which means that the $k$-dimension of the stable endomorphism ring of $Y_{1ij}$ is at least $2$. If $d=3$ then $Y_{2ij}=\Omega(Y)$, and so $\Omega(\mathfrak{C})=\mathfrak{C}$. For $d>3$, $Y_{2ij}$ has a non-trivial endomorphism factoring through $T_i$ which does not factor through a projective $B$-module, and thus the $k$-dimension of the stable endomorphism ring of $Y_{2ij}$ is at least $2$. By Proposition \ref{prop:psl1} and by Remark \ref{rem:psl1uni}, it follows that for $i\neq j$ in $\{1,2\}$, the components $\mathfrak{C}$ correspond to the components of the stable Auslander-Reiten quiver of $B_0$ as in Proposition \ref{prop:psl1}(iii). 
This implies part (iii) of Proposition \ref{prop:psl2}.

Since we have matched up all components of the stable Auslander-Reiten quiver of $B_0$ containing modules with stable endomorphism ring $k$ to certain components of the stable Auslander-Reiten quiver of $B$, it follows that all other components of the stable Auslander-Reiten quiver of $B$ do not contain any modules with stable endomorphism ring $k$. This completes the proof of Proposition \ref{prop:psl2}.

\begin{subrem}
\label{rem:uniliftpsl2}
Similarly to Remark \ref{rem:uniliftpsl1}, it follows that if $Y=\begin{array}{c}T_i\\T_j\end{array}$ where $i\neq j$ in $\{1,2\}$, then the universal mod $2$ deformation of $Y$ is represented by the uniserial $kG$-module $\overline{U}$ of length $2\cdot 2^{d-2}$ with descending radical series
$$(T_i,T_j,T_i,T_j,\ldots,T_i,T_j).$$
Moreover the uniserial $kG$-module $\overline{U'}\cong \overline{U}/Y$ of length $2\cdot (2^{d-2}-1)$ defines a lift of $Y$ over $k[t]/(t^{2^{d-2}-1})$.
\end{subrem}

\medskip
We next turn to blocks as in \S\ref{ss:a7}.
\begin{subprop}
\label{prop:a7}
Let $B$ be a block of $kG$ which is Morita equivalent to $\Lambda=kQ/I$ with $Q$ and $I$ as in $\S\ref{ss:a7}$. 
Then the order of the defect groups of $B$ is $2^d=8$.
Denote the three simple $B$-modules by $T_0$, $T_1$ and $T_2$, where $T_i$ corresponds to $S_i$, for $i\in\{0,1,2\}$, under the Morita equivalence, and define $Y_{01102}$ to be the $B$-module corresponding to the $\Lambda$-string-module 
with underlying string $\beta\alpha^{-1}\gamma\eta$.
Let $\mathfrak{C}$ be a component of the stable Auslander-Reiten quiver of $B$. If $\mathfrak{C}$ contains a module with endomorphism ring $k$, then $\mathfrak{C}$ contains a simple module or a uniserial module of length $4$, or $\Omega(\mathfrak{C})$ contains $T_0$.
\begin{enumerate}
\item[i.] Suppose $\mathfrak{C}$ contains $T_0$. Then $\mathfrak{C}$ and $\Omega(\mathfrak{C})$ are both of type $\mathbb{Z}A_\infty^\infty$. All modules $M$ in $\mathfrak{C}\cup\Omega(\mathfrak{C})$ have stable endomorphism ring equal to $k$ and $R(G,M)/2R(G,M)\cong k$.
\item[ii.] Let $i\neq j$ be in $\{1,2\}$ and define $T_{i,0,j,0}=\begin{array}{c}T_i\\T_0\\T_j\\T_0\end{array}$ and $T_{0,j,0,i}=\begin{array}{c}T_0\\T_j\\T_0\\T_i\end{array}$. Suppose $\mathfrak{C}$ contains $T_{i,0,j,0}$. Then $\mathfrak{C}$ is a $3$-tube with $T_{i,0,j,0}$ and $T_{0,j,0,i}$ belonging to its boundary, and $\mathfrak{C}=\Omega(\mathfrak{C})$. If $i=2$, $\mathfrak{C}$ also contains $T_2$.
In both cases, if $M$ lies in $\mathfrak{C}$ having stable endomorphism ring $k$, then $M\in\{T_{i,0,j,0},\Omega^2(T_{i,0,j,0}),\Omega^4(T_{i,0,j,0})\}$ up to isomorphism, and $R(G,M)/2R(G,M)\cong k$.
\item[iii.] Suppose $Y\in \{T_1,Y_{01102}\}$ and $\mathfrak{C}$ contains $Y$. 
Then $\mathfrak{C}=\Omega(\mathfrak{C})$ is of type $\mathbb{Z}A_\infty^\infty$. If $Y=Y_{01102}$ then $\mathfrak{C}$ contains no $B$-modules with endomorphism ring $k$. In both cases, if $M$ lies in $\mathfrak{C}$ having stable endomorphism ring $k$, then $M$ is isomorphic to $\Omega^j(Y)$ for some integer $j$, and $R(G,M)/2R(G,M)\cong k[t]/(t^2)$. 
\end{enumerate}
The only components of the stable Auslander-Reiten quiver of $B$ containing modules with stable endomorphism ring $k$ are the ones in $(i)-(iii)$.
\end{subprop}

The proof of Proposition \ref{prop:a7} is similar to the proof of Proposition \ref{prop:psl2} using the following result in place of Lemma \ref{lem:psl2}.
\begin{sublem}
\label{lem:a7}
Let $\Lambda=kQ/I$ with $Q$ and $I$ as in $\S\ref{ss:a7}$.
For $i\in\{0,1,2\}$, let $\mathfrak{C}_i$ be the component of the stable Auslander-Reiten quiver of $\Lambda$ containing $S_i$.
Let $M$ be an indecomposable $\Lambda$-module with $\mathrm{End}_\Lambda(M)=k$. Then $M$ either lies in $\mathfrak{C}_0\cup\Omega(\mathfrak{C}_0)$, or $M$ lies in $\mathfrak{C}_i$ for an $i\in\{1,2\}$, or $M$ is uniserial of length $4$.
\end{sublem}

\begin{subrem}
\label{rem:unilifta7}
Similarly to Remark \ref{rem:uniliftpsl1}, it follows that the universal mod $2$ deformation of $T_1$ is represented by the uniserial $kG$-module $\overline{U}_1$ 
with radical series $(T_1,T_1)$.

The universal mod $2$ deformation of $Y_{01102}$ can be obtained as follows. There is a non-split short exact sequence of $kG$-modules
$$0\to \Omega(Y_{01102}) \to P_{T_1}\oplus P_{T_2} \to Y_{01102} \to 0$$
where $P_{T_1}$ (resp. $P_{T_2}$) is the projective indecomposable $kG$-module with top $T_1$ (resp. $T_2$).
Moreover, there is a natural surjection $\Omega(Y_{01102})\to Y_{01102}$ with kernel isomorphic to $T_2$. Thus there is a non-split short exact sequence of $kG$-modules
\begin{equation}
\label{eq:obrother}
0\to Y_{01102}\to P_{T_1}\oplus P_{T_2}/T_2 \to Y_{01102}\to 0
\end{equation}
representing a non-zero element in $\mathrm{Ext}^1_{kG}(Y_{01102},Y_{01102})$. Hence the universal mod $2$ deformation of $Y_{01102}$ is represented by the $kG$-module $P_{T_1}\oplus P_{T_2}/T_2$. 
\end{subrem}


\section{Universal deformation rings}
\label{s:udr}
\setcounter{equation}{0}
\setcounter{figure}{0}

We assume the notations from \S\ref{s:dihedralsylow} and \S\ref{s:defmod2}.
In this section, we determine the universal deformation rings of all finitely generated $kG$-modules which belong to $B$ and have stable endomorphism ring equal to $k$. Since the case $d=2$ has been done in \cite{bl}, we assume throughout this section that $d\ge 3$. 

In \S\ref{ss:mainthm}, we state our main Theorem \ref{thm:main} and some consequences concerning local complete intersections (see Corollary \ref{cor:lci}). In \S\ref{ss:3tubes}, we analyze the modules belonging to the boundaries of $3$-tubes. In \S\ref{ss:proof}, we then use these results together with the results from \S\ref{ss:ordinarydihedral} to prove Theorem \ref{thm:main}.

\subsection{The Main Theorem}
\label{ss:mainthm}
\setcounter{equation}{0}
\setcounter{figure}{0}

\begin{thm}
\label{thm:main}
Let $G$ be a finite group, and let $B$ be a block of $kG$ with dihedral defect group $D$ of order $2^d$ where $d\ge 3$. Assume that $B$ is Morita equivalent to the principal $2$-modular block of a finite simple group. Then $B$ is Morita equivalent to a block as in $\S\ref{ss:psl1}$, $\S\ref{ss:psl2}$ or $\S\ref{ss:a7}$. Let $V$ be  a finitely generated indecomposable $B$-module with stable endomorphism ring $k$, and let $\mathfrak{C}$ be the component of the stable Auslander-Reiten quiver of $B$ containing $V$. 
\begin{enumerate}
\item[i.] If $\mathfrak{C}$ or $\Omega(\mathfrak{C})$ is as in part $(i)$ of Propositions $\ref{prop:psl1}$, $\ref{prop:psl2}$ or $\ref{prop:a7}$, then $R(G,V)$ is isomorphic to a quotient ring of $W$.
\item[ii.] If $\mathfrak{C}$ or $\Omega(\mathfrak{C})$ is as in part $(ii)$ of Propositions $\ref{prop:psl1}$, $\ref{prop:psl2}$ or $\ref{prop:a7}$, then $R(G,V)$ is isomorphic to $k$.
\item[iii.] If $\mathfrak{C}$ or $\Omega(\mathfrak{C})$ is as in part $(iii)$ of Propositions $\ref{prop:psl1}$, $\ref{prop:psl2}$ or $\ref{prop:a7}$, then $R(G,V)$ is isomorphic to $$W[[t]]/(p_d(t)(t-2),2\,p_d(t))$$ as a $W$-algebra where $p_d(t)\in W[t]$ is as in Definition $\ref{def:seemtoneed}$.
\end{enumerate}
In all cases $(i) - (iii)$, $R(G,V)$ is isomorphic to a subquotient ring of $WD$.
The only components of the stable Auslander-Reiten quiver of $B$ containing modules with stable endomorphism ring $k$ are the ones in $(i) - (iii)$.
\end{thm}

\begin{subrem}
\label{rem:liftsoverw}Using Lemma \ref{lem:Wlift}, we can prove that if $V$ belongs to $\mathfrak{C}$ as in part (i) of Theorem \ref{thm:main} then $R(G,V)\cong W$. 
\end{subrem}

\begin{subcor}
\label{cor:lci}
Assuming the notations of Theorem $\ref{thm:main}$, suppose $\mathfrak{C}$ is as in part $(iii)$ of Propositions $\ref{prop:psl1}$, $\ref{prop:psl2}$ or $\ref{prop:a7}$. 
Let $V$ be a finitely generated indecomposable $kG$-module in $\mathfrak{C}$ with stable endomorphism ring $k$. Then $R(G,V)$ is not a complete intersection ring. 

In particular, there is an infinite series of finite groups $G$ and indecomposable $kG$-modules $V$ such that $R(G,V)$ is not a complete intersection. This series is given by $G=\mathrm{PSL}_2(\mathbb{F}_q)$ for $q$ an odd prime power and $8$ dividing $\#G$, and $V=\Omega^i(V')$ for $i\in\mathbb{Z}$, where $V'$ is a uniserial $kG$-module belonging to the principal block of $kG$ with  radical series length $4$ $($resp. radical series length $2$ and no composition factor isomorphic to the trivial simple module$)$ in case $q\equiv 1\mod 4$ $($resp. $q\equiv 3\mod 4$$)$.
Moreover, in case $G=A_7$, the unique $($up to isomorphism$)$ irreducible $kG$-module $V$ of dimension $14$ provides an example of an irreducible $V$ such that $R(G,V)$ is not a complete intersection.
\end{subcor}

\subsection{Modules at the boundaries of three-tubes}
\label{ss:3tubes}
\setcounter{equation}{0}
\setcounter{figure}{0}

We first summarize the known facts about the modules belonging to the boundaries of $3$-tubes in the stable Auslander-Reiten quiver of one of the blocks $B$ under consideration. These facts can be found e.g. in \cite[Chapter V]{erd}.

\begin{subfacts}
\label{facts:brauererdmann}
Let $G$ be a finite group, and let $B$ be a block of $kG$ with dihedral defect group $D$ of order $2^d$ where $d\ge 3$. Assume that $B$ contains exactly three isomorphism classes of simple $kG$-modules. Then the stable Auslander-Reiten quiver of $B$ has exactly two $3$-tubes. Suppose $V$ is a finitely generated indecomposable $B$-module belonging to the boundary of a $3$-tube. Then the vertices of $V$ are Klein four groups. Let $K$ be a vertex of $V$.
\begin{enumerate}
\item[i.] The quotient group $N_G(K)/C_G(K)$ is isomorphic to a symmetric group $S_3$.
\item[ii.] There is a block $b$ of $kN_G(K)$ with $b^G=B$ such that the Green correspondent $fV$ of $V$ belongs to the boundary of a $3$-tube in the stable Auslander-Reiten quiver of $b$. Moreover, $b$ is Morita equivalent to $kS_4$ modulo the socle.
\end{enumerate}
\end{subfacts}

\begin{sublem}
\label{lem:1}
Suppose $K$ and $b$ are as in Facts $\ref{facts:brauererdmann}$. Let $H=N_G(K)$, let $C=C_G(K)$, and let $N$ be a normal subgroup of $H$ of index $2$ containing $C$. 
\begin{enumerate}
\item[i.] There is a unique block $b_1$ of $kN$ with defect group $K$ which is covered by $b$. Moreover, $b_1^H=b$ and $b_1$ is Morita equivalent to $kA_4$.
\item[ii.] Let $g\in H$ with $\langle gN\rangle = H/N$. Then there is a simple $b_1$-module $S_0$ with $g(S_0)\cong S_0$, where $g(S_0)$ denotes the $kN$-module such that  $gxg^{-1}$ acts on $g(S_0$) the same way as $x\in N$ acts on $S_0$. The other two representatives of non-isomorphic simple $b_1$-modules $S_1$ and $S_2$ satisfy $g(S_1)\cong S_2$ and $g(S_2)\cong S_1$.
\item[iii.] The stable Auslander-Reiten quiver of $b_1$ has two $3$-tubes whose modules all have vertex $K$. Each $b_1$-module at the end of one $3$-tube has as source the band module $S_\omega=M(XY^{-1},\omega,1)$, and each $b_1$-module at the end of the other $3$-tube has as source the band module  $S_{\omega^2}=M(XY^{-1},\omega^2,1)$, where $XY^{-1}$ denotes the single band of $k K$ and $\omega$ is a primitive cube root of unity in $k$.
\end{enumerate}
\end{sublem}

\begin{proof}
It follows from \cite[Thm. 15.1]{alp} that there is a unique block $b_1$ of $kN$ with defect group $K$ which is covered by $b$, and that $b_1^H=b$. 
Since $N/C_N(K)$ has order $3$, it follows e.g. from \cite[Proof of Prop. V.2.14]{erd} that $b_1$ is Morita equivalent to $kA_4$. This implies part (i).

Now we turn to part (ii). Since  there are three isomorphism classes of simple $b_1$-modules and $[H:N]=2$, there is a simple $b_1$-module $S_0$ with $g(S_0)\cong S_0$.
So we only need to show that $g(S_1)$ is not isomorphic to $S_1$. To get a contradiciton, suppose that $g(S_1)\cong S_1$, and hence also $g(S_2)\cong S_2$. For $i=0,1,2$, consider the $kH$-module $X_i=\mathrm{Ind}_N^H(S_i)$. Then by Mackey's Theorem $(X_i)_N\cong S_i\oplus g(S_i)\cong S_i\oplus S_i$. In particular, $X_i$ is a $b$-module. If $Y_i$ is a simple $b$-module in the socle of $X_i$ then $(Y_i)_N$ is a submodule of $(X_i)_N\cong S_i\oplus S_i$. Since $S_0$, $S_1$ and $S_2$ are pairwise non-isomorphic, it follows that $Y_0$, $Y_1$ and $Y_2$ are pairwise non-isomorphic. But this is a contradiction to $b$ having only two non-isomorphic simple modules. Hence $g(S_1)\cong S_2$ and vice versa.

Finally, we turn to part (iii). It follows e.g. from \cite[Proof of Thm. V.4.1]{erd} that all the $b_1$-modules in a $3$-tube have the same vertex and that this vertex must be a Klein four group. Since $K$ is a defect group of $b_1$, $K$ is a vertex. Let $\tilde{\nu}\in N$ such that $\langle \tilde{\nu}C\rangle =N/C$. Then $\tilde{\nu}$ acts on $K$ by conjugation which induces an automorphism $\nu$ of $K$ of order $3$. Let $A=kK * \langle \nu\rangle$ be the skew group ring of $kK$ with $\langle \nu\rangle$. Then $A\cong kA_4$ (see e.g. \cite[Cor. V.2.4.1]{erd}). Let $\beta$ be a block of $k C$ which is covered by $b_1$. Then by \cite[Proof of Prop. V.2.14]{erd}, $\beta$ is Morita equivalent to $k K$. Let $P$ be the unique projective indecomposable $\beta$-module (up to isomorphism) and identify $P$ with an inner direct summand of $kC$. Then $\mathrm{End}_{kC}(P)=kK$. Since $K$ is central in $C$ and we assume $P\subseteq kC$, right multiplication by any element of $kK$ defines a $kC$-module endomorphism of $P$. If $P^N=\mathrm{Ind}_C^N(P)$ then  $\mathrm{End}_{kN}(P^N)\cong A=kK * \langle \nu\rangle$ by \cite[Proof of Prop. V.2.14]{erd}. Moreover, it follows from \cite[\S V.2.9]{erd} and \cite[Proof of Prop. V.2.14]{erd} that the right action of $kK$ on $P$ by right multiplication extends to a right action of $A$ on $P^N$. In particular, the elements of $kK$ act by right multiplication on $P^N$ when identifying $P^N$ with an inner direct summand of $kN$. In \cite{conlon}, the indecomposable $A$-modules are described as direct summands of induced indecomposable $kK$-modules. This shows that the direct sum of the $A$-modules belonging to the boundary of a $3$-tube of the stable Auslander-Reiten quiver of $A$ is isomorphic to $A\otimes_{kK}S_\rho$ where for one $3$-tube $\rho=\omega$ and for the other $3$-tube $\rho=\omega^2$ and $S_\rho$ is as in part (iii) of the statement of the lemma. Since $P^N$ gives a Morita equivalence between $b_1$ and $A$, the functor $P^N\otimes_A-$ sends the $A$-modules belonging to the boundary of a $3$-tube of the stable Auslander-Reiten quiver of $A$ to $b_1$-modules which belong to the boundary of a $3$-tube of the stable Auslander-Reiten quiver of $b_1$. Thus the direct sum of the $b_1$-modules belonging to the boundary of a $3$-tube is isomorphic to
$$P^N\otimes_{A} (A \otimes_{kK}\,S_\rho) \cong P^N\otimes_{kK}\,S_\rho$$
where $kK$ acts on $P^N$ by right multiplication and $\rho=\omega$ or $\omega^2$. Since we view $P^N$ as an inner direct summand of $kN$, it follows that the direct sum of the $b_1$-modules belonging to the boundary of a $3$-tube is isomorphic to a direct summand of $kN\otimes_{kK}\,S_\rho =\mathrm{Ind}_K^N\,S_\rho$. But this means that $S_\rho$ is a source of  each of the $b_1$-modules belonging to the boundary of the $3$-tube under consideration.  This completes the proof of Lemma \ref{lem:1}.
\end{proof}

\begin{sublem}
\label{lem:2}
Assume the notations from Lemma $\ref{lem:1}$. Let $E_0$ and $E$ be the two non-isomorphic simple $b$-modules such that $U_E=\begin{array}{c}E\\E\end{array}$ belongs to the boundary of the $3$-tube of the stable Auslander-Reiten quiver of $b$. Then $U_E$ is isomorphic to the induced module $\mathrm{Ind}_N^H \,M$ for some $b_1$-module $M$ which belongs to the boundary of a $3$-tube of the stable Auslander-Reiten quiver of $b_1$. Moreover, $K$ is a vertex of $U_E$, and the source of $U_E$ is $S_E=S_\omega$ which is conjugate to $S_{\omega^2}$ in $N_H(K)=H$. In particular, the restriction of $S_E$ to any proper subgroup of $K$ is projective.
\end{sublem}

\begin{proof}
For $i\in\{0,1,2\}$, let $X_i=\mathrm{Ind}_N^H\,S_i$. By Mackey's Theorem, $(X_i)_N\cong S_i\oplus g(S_i)$, where $g$ is as in part (ii) of Lemma \ref{lem:1}. Hence $(X_0)_N\cong S_0\oplus S_0$ and $(X_i)_N\cong S_1\oplus S_2$ for $i=1,2$. In particular, for $i\in\{0,1,2\}$, $X_i$ is a $b$-module. Let $Y_i$ be a simple $b$-module in the socle of $X_i$. Since $g(Y_i)\cong Y_i$, it follows that $X_1\cong X_2$ is simple. Call this simple $b$-module $Y_{12}$. By \cite[Lemmas 8.5 and 8.6]{alp},  induced modules of projective modules are projective, and the induced modules of the composition factors of a $b_1$-module $T$ provide the composition factors of the induced module $\mathrm{Ind}_N^H\, T$. Since the projective cover of $S_0$ has composition factors $S_0,S_1,S_2,S_0$ and since the projective covers of $E_0$ and $E$ have 6, resp. 5, composition factors, it follows that $\mathrm{Ind}_N^H\,S_0$ has composition factors $Y_0$ and $Y_0$. This implies that $Y_0=E_0$ and $Y_{12}=E$.

Now consider the indecomposable $b_1$-module $M=\begin{array}{c}S_1\\S_2\end{array}$, which belongs to the boundary of a $3$-tube of the stable Auslander-Reiten quiver of $b_1$. By \cite[Lemma 8.6(5)]{alp}, it follows that $\mathrm{Ind}_N^H\,M$ satisfies a non-split short exact sequence 
$$0\to E \to \mathrm{Ind}_N^H\,M \to E \to 0,$$
which implies $U_E=\begin{array}{c}E\\E\end{array}\cong \mathrm{Ind}_N^H\,M$.

It follows e.g. from \cite[Proof of Thm. V.4.1]{erd} that all the modules in $3$-tubes of the stable Auslander-Reiten quiver of $b$ have Klein four groups as vertices. By Lemma \ref{lem:1}(iii), $S_\rho$ is a source of $M$ where $\rho=\omega$ or $\omega^2$. Hence $U_E$ is a direct summand of $\mathrm{Ind}_N^H (\mathrm{Ind}_K^N\, S_\rho) \cong \mathrm{Ind}_K^H\, S_\rho$, which means that $K$ is a vertex of $U_E$ and $S_\rho$ is a source of $U_E$. But similarly, we also get that $S_{\rho^2}$ is a source of $U_E$. Because sources are unique up to conjugation in $N_H(K)=H$, $S_\omega$ is conjugate to $S_{\omega^2}$ in $H$. Since $S_E=S_\omega$ is $2$-dimensional over $k$ and each proper subgroup of $K$ of order $2$ acts non-trivially on $S_E$, the restriction of $S_E$ to any proper subgroup of $K$ is projective. This completes the proof of Lemma \ref{lem:2}.
\end{proof}

\begin{subprop}
\label{prop:3tubedef}
Assume the notations of Lemmas $\ref{lem:1}$ and $\ref{lem:2}$. If $U$ is an indecomposable $b$-module belonging to the boundary of the $3$-tube of the stable Auslander-Reiten quiver of $b$, then $U$ has a universal deformation ring with $R(H,U)\cong k$.
\end{subprop}

\begin{proof}
It follows that $U$ lies in the $\Omega$-orbit of the indecomposable $b$-module $U_E$. By Lemma \ref{lem:2},  $U_E\cong \mathrm{Ind}_N^H\,M$ for some $b_1$-module $M$ which belongs to the boundary of a $3$-tube of the stable Auslander-Reiten quiver of $b_1$. Since $b_1$ is Morita equivalent to $kA_4$, it follows from \cite[Prop. 3.4]{bl} that $\underline{\mathrm{End}}_{kN}(M)= k$ and $R(N,M)\cong k$. Since $\underline{\mathrm{End}}_{kH}(U_E)= k$ and
$$\mathrm{dim}_k \mathrm{Ext}_{kH}^1(U_E,U_E) = 0 = \mathrm{dim}_k \mathrm{Ext}_{kN}^1(M,M),$$
Proposition \ref{prop:induceddef} implies that $R(H,U_E)\cong R(N,M)\cong k$. By Lemma \ref{lem:defhelp}, this implies Proposition \ref{prop:3tubedef}.
\end{proof}

\begin{subcor}
\label{cor:3tubedef}
Assume the notations of Theorem $\ref{thm:main}$ and suppose $\mathfrak{C}$ is as in part $(ii)$ of Theorem $\ref{thm:main}$. Let $V$ be a finitely generated $kG$-module in $\mathfrak{C}$ with stable endomorphism ring $k$. Then $R(G,V)\cong k$.
\end{subcor}

\begin{proof}
It follows from part (ii) of Propositions $\ref{prop:psl1}$, $\ref{prop:psl2}$ or $\ref{prop:a7}$ that $V$ is a finitely generated indecomposable $B$-module belonging to the boundary of a $3$-tube. Hence we may use the notations of Facts \ref{facts:brauererdmann}, Lemmas \ref{lem:1}, \ref{lem:2} and Proposition \ref{prop:3tubedef}. We claim that
\begin{equation}
\label{eq:yuckyuck}
\mathrm{Ind}_H^G\,fV \cong V \oplus \mathrm{\;projectives}.
\end{equation}
Since $fV$ belongs to the boundary of the $3$-tube of the stable Auslander-Reiten quiver of $b$, it lies in the $\Omega$-orbit of $U_E$ from Lemma \ref{lem:2}. Since the Green correspondence commutes with $\Omega$ (see e.g. \cite[Prop. 20.7]{alp}), it is enough to show $(\ref{eq:yuckyuck})$ in case $fV=U_E$. Using Green correspondence, we know $\mathrm{Ind}_H^G U_E \cong V\oplus X$ where $X$ is relatively $\mathfrak{x}$-projective and the groups in $\mathfrak{x}$ have the form $sKs^{-1} \cap K$ for some $s\in G$, $s\not\in H$. Suppose there is an indecomposable summand $Y$ of $X$ which has a non-trivial vertex $Q$. Then $Q$ has order $2$. Because $\mathrm{Ind}_H^G\, U_E$ is a direct summand of $\mathrm{Ind}_K^G\,S_E$, we get that $Y_Q$ is a direct summand of 
$$(\mathrm{Ind}_K^G\,S_E)_Q\cong \bigoplus_{t\in Q\backslash G/K} \mathrm{Ind}_{Q\cap t K t^{-1}}^Q (t(S_E)_{Q\cap t K t^{-1}}).$$
Since each term $t(S_E)_{Q\cap t K t^{-1}}$  is projective by Lemma \ref{lem:2}, it follows that $Y_Q$ is projective. But $Y$ is a direct summand of $\mathrm{Ind}_Q^G (Y_Q)$, hence $Y$ is projective. This is a contradiction to $Q$ being a non-trivial group. This  implies $(\ref{eq:yuckyuck})$. By Proposition \ref{prop:induceddef} it follows that  $R(G,V)\cong R(H,fV)$. Hence Corollary \ref{cor:3tubedef} follows from Proposition \ref{prop:3tubedef}.
\end{proof}

\subsection{Proof of Theorem $\ref{thm:main}$}
\label{ss:proof}
\setcounter{equation}{0}
\setcounter{figure}{0}

Part (i) and part (ii) of Theorem \ref{thm:main} and the very last statement about the components of the stable Auslander-Reiten quiver containing modules with stable endomorphism ring $k$ follow from Propositions \ref{prop:psl1}, \ref{prop:psl2} or \ref{prop:a7} together with Corollary \ref{cor:3tubedef}.
Since it follows from Lemma \ref{lem:dihedraltrick} that $W[[t]]/(p_d(t)(t-2),2\,p_d(t))$ is isomorphic to a subquotient ring of $WD$ 
when $p_d(t)\in W[t]$ is as in Definition $\ref{def:seemtoneed}$, it only remains to prove part (iii) of Theorem \ref{thm:main}.

Let $\mathfrak{C}$ be as in part (iii) of Theorem \ref{thm:main}, and let $V$ be  a finitely generated indecomposable $kG$-module with stable endomorphism ring $k$ belonging to $\mathfrak{C}$. By part (iii) of Propositions \ref{prop:psl1}, \ref{prop:psl2} or \ref{prop:a7}, $R(G,V)/2R(G,V)\cong k[t]/(t^{2^{d-2}})$.

We first look at the case when $B$ is as in \S\ref{ss:psl1}. By Proposition \ref{prop:psl1}(iii) and Remark \ref{rem:psl1uni}, it is enough to consider $V$ where $V$ is either $Y_1=\begin{array}{c}T_1\\T_0\\T_2\\T_0\end{array}$ or $Y_2=\begin{array}{c}T_2\\T_0\\T_1\\T_0\end{array}$. Let $V=Y_1$. (The case $V=Y_2$ is proved similarly.) By Remark \ref{rem:uniliftpsl1}, the universal mod $2$ deformation of $Y_1$ is represented by the uniserial $kG$-module $\overline{U}$ with descending radical series
$$(T_1,T_0,T_2,T_0,T_1,T_0,T_2,T_0,\ldots,T_1,T_0,T_2,T_0)$$
of length $4\cdot 2^{d-2}$. Moreover, the uniserial $kG$-module $\overline{U'}=\overline{U}/Y_1$ of length $4\cdot (2^{d-2}-1)$ defines a lift of $Y_1$ over $k[t]/(t^{2^{d-2}-1})$. We show that there is a surjective continuous $W$-algebra homomorphism $\tau: R\to R'$ where $R'=W[[t]]/(p_d(t))$ is as in Definition \ref{def:seemtoneed}. We use the notations from \S\ref{ss:ordinarydihedral}. It follows from the decomposition matrix in Figure \ref{fig:decompsl1} and \cite[Prop. (23.7)]{CR} that there is a $W$-pure $WG$-sublattice $X'$ of the projective indecomposable $WG$-module $P_1^W$ with top $T_1$ such that $U'=P_1^W/X'$ has $F$-character
$$\sum_{\ell=2}^{d-1}\rho_\ell=\sum_{i=1}^{2^{d-2}-1}\chi_{5,i}.$$
Then $U'/2U'$ is an indecomposable $kG$-module with top $T_1$ which has the same decomposition factors as $\overline{U'}$ and is thus isomorphic to $\overline{U'}$. Hence $U'$ is a lift of $\overline{U'}$ over $W$. We want to show that $U'$ is in fact a lift of $Y_1$ over $R'$. We first prove that $R'$ is isomorphic to a $W$-subalgebra of $\mathrm{End}_{WG}(U')$. Let $\sigma$ be an element of order $2^{d-1}$ in $D$, and let $t(C)$ be the class sum of the conjugacy class $C$ of $\sigma$ in $G$. Since $t(C)$ lies in the center of $WG$, multiplication by $t(C)$ defines a $WG$-module endomorphism of $U'$. Since $\mathrm{End}_{WG}(U')$ can naturally be identified with a subring of $\mathrm{End}_{FG}(F\otimes_W U')$, $t(C)$ acts on $U'$ as multiplication by a scalar in $F\otimes_W R'$. This scalar can be read off from the action of $t(C)$ on $F\otimes_W U'=\bigoplus_{\ell=2}^{d-1}V_\ell$. By (\ref{eq:thatsit}), there exists a unit $\omega$ in $W$ such that for $2\le \ell \le d-1$, the action of $t(C)$ on $V_\ell$ is given as multiplication by
$$\omega\cdot (\zeta_{2^{d-1}}^{2^{d-1-\ell}}+\zeta_{2^{d-1}}^{-2^{d-1-\ell}})$$
when we identify $\mathrm{End}_{FG}(V_\ell)$ with $\mathrm{End}_{F(\zeta_{2^\ell}+\zeta_{2^\ell}^{-1})G}(X_\ell)$ for an absolutely irreducible $F(\zeta_{2^\ell}+\zeta_{2^\ell}^{-1})G$-constituent $X_\ell$ of $V_\ell$ with character $\chi_{5,2^{d-1-\ell}}$. Recall from Remark \ref{rem:ohyeah} that $R'$ can be identified with the $W$-subalgebra of $\bigoplus_{\ell=2}^{d-1} W[\zeta_{2^\ell}+\zeta_{2^\ell}^{-1}]$ generated by $(\zeta_{2^\ell}^{r_\ell}+\zeta_{2^\ell}^{-r_{\ell}})_{\ell=2}^{d-1}$ for any sequence $(r_\ell)_{\ell=2}^{d-1}$ of odd numbers. This implies that $R'$ is isomorphic to a $W$-subalgebra of $\mathrm{End}_{WG}(U')$, and hence $U'$ is an $R'G$-module. We next prove that $U'$ is free as an $R'$-module. Since $U'$ is finitely generated as a $W$-module, it is also finitely generated as an $R'$-module. Since $R'$ is a local ring, it follows by Nakayama's Lemma that any $k$-basis $\{\overline{b}_1,\ldots, \overline{b}_s\}$ of $U'/\mathrm{max}(R')U'\cong Y_1$ can be lifted to a set $\{b_1,\ldots, b_s\}$ of generators of $U'$ over $R'$. Since $F\otimes_W U'$ is a free $(F\otimes_W R')$-module of rank $s$, it follows that $b_1,\ldots, b_s$ are linearly independent over $R'$. Since $\mathrm{End}_{FG}(F\otimes_W U')\cong F\otimes_WR'$, this then implies that $\mathrm{End}_{WG}(U')\cong R'$. Hence $U'$ is a lift of $Y_1$ over $R'$. We therefore have a continuous $W$-algebra homomorphism $\tau:R\to R'$ relative to $U'$. Since $U'/2U'$ is indecomposable as a $kG$-module, $\tau$ must be surjective. 
By Lemma \ref{lem:herewegoagain},  it follows that $R(G,V)\cong W[[t]]/(p_d(t)(t-2c),a2^mp_d(t))$ for some $c\in W$, $a\in\{0,1\}$ and $0 < m\in\mathbb{Z}$. If $a=0$ (resp. $a=1$), the natural projection $R(G,V)\to W[[t]]/(p_d(t)(t-2c))$ (resp. $R(G,V)\to (W/2^m W)[[t]]/(p_d(t) (t-2c))$) gives a lift of $\overline{U}$, when regarded as a $kG$-module, over $W$ (resp. over $W/2^mW$). But $\overline{U}\cong \Omega^{-1}(T_1)$, which implies $R(G,\overline{U})\cong k$ by Corollary \ref{cor:3tubedef}. Hence $a=1$ and $m=1$, and part (iii) of Theorem \ref{thm:main} follows in case $B$ is as in \S\ref{ss:psl1}.

The case when $B$ is as in \S\ref{ss:psl2} is proved similarly to the case when $B$ is as in \S\ref{ss:psl1}, using Proposition \ref{prop:psl2}(iii), Remark \ref{rem:uniliftpsl2} and the decomposition matrix in Figure \ref{fig:decompsl2}.

Finally we look at the case when $B$ is as in \S\ref{ss:a7}. By Proposition \ref{prop:a7}(iii), it is enough to consider $V\in\{T_1,Y_{01102}\}$. 
If we show that $V$ has a lift over $W$, then it follows by Lemma \ref{lem:herewegoagain}, using Remark \ref{rem:unilifta7} and Corollary \ref{cor:3tubedef}, that $R(G,V)\cong W[[t]]/(t(t-2),2\,t)$. 
Using the decomposition matrix in Figure \ref{fig:decoma7} and \cite[Prop. (23.7)]{CR} we see that there is a lift of $T_1$ over $W$.
To see that $Y_{01102}$ has a lift over $W$, we consider the two uniserial $kG$-modules
$Z_1=\begin{array}{c}T_1\\T_0\end{array}$ and $Z_2= \begin{array}{c}T_2\\T_0\\T_1\end{array}$.
Then there is a non-split short exact sequence of $kG$-modules
$$0\to Z_2\to Y_{01102}\to Z_1\to 0$$
with $\mathrm{Ext}^1_{kG}(Z_1,Z_2)=k$. It follows from the decomposition matrix in Figure  \ref{fig:decoma7} and \cite[Prop. (23.7)]{CR} that there exists a lift $X_i$ of $Z_i$ over $W$ for $i=1,2$. Moreover, the $F$-character of $F\otimes_W X_1$ (resp. of $F\otimes_W X_2$) is $\chi_2$ (resp. $\chi_4$). In particular, $\mathrm{Hom}_{FG}(F\otimes_W X_1,F\otimes_WX_2)=0$. Since $\mathrm{Hom}_{kG}(Z_1,Z_2)=k$, it follows from Lemma \ref{lem:Wlift} that there is a lift of $Y_{01102}$ over $W$. 
This completes the proof of Theorem \ref{thm:main}.


\section{Stable endomorphism rings and Ext groups}
\label{s:stableend}
\setcounter{equation}{0}
\setcounter{figure}{0}

Assume the notations of \S\ref{s:defmod2}. In this section we complete the proof of  Proposition \ref{prop:psl1}, by determining which components of the stable Auslander-Reiten quiver of a block $B$ as in \S\ref{ss:psl1} contain modules with stable endomorphism ring $k$ and by determining the Ext groups for these modules.
Since there are stable equivalences of Morita type between $B$ and blocks as in \S\ref{ss:psl2} or \S\ref{ss:a7}, we determine at the same time which components of the stable Auslander-Reiten quiver of the blocks in \S\ref{ss:psl2} and in \S\ref{ss:a7} contain modules with stable endomorphism ring $k$ and also determine the Ext groups for these modules.

Recall that $B$ is Morita equivalent to the basic algebra $\Lambda$ of a special biserial algebra, so we can use the description of indecomposable $\Lambda$-modules as string and band modules (see \S\ref{s:prelimstring}). 
In \S \ref{ss:stringhom}, we give a description of the homomorphisms between string modules as determined in \cite{krau} and provide a short-hand notation for such homomorphisms. We also give a criterion which helps determine stable homomorphisms between string modules. 
In \S\ref{ss:kendo}, we prove Lemmas \ref{lem:psl1endos0}, \ref{lem:psl1ci} and \ref{lem:psl1uni}. In \S\ref{ss:infinity}, we consider all components of the stable Auslander-Reiten quiver of type $\mathbb{Z}A_\infty^\infty$ and prove that the only such components containing a module with stable endomorphism ring $k$ are precisely those components containing a module $M$ with $\mathrm{End}_\Lambda(M)=k$ or $\mathrm{End}_\Lambda(\Omega(M))=k$. In \S\ref{ss:bands}, we consider the components of the stable Auslander-Reiten quiver which are $1$-tubes and prove that no $1$-tube contains any modules with stable endomorphism ring $k$.

We freely use \S \ref{s:prelimstring} without always explicitly referring to particular results. We especially use the phrase ``canonical $k$-basis'' for string modules as introduced in \S\ref{sss:strings} to be able to readily write down homomorphisms between string modules.


\subsection{Homomorphisms between string modules}
\label{ss:stringhom}
\setcounter{equation}{0}
\setcounter{figure}{0}

\begin{subrem}
\label{rem:stringhoms}
Let $\Lambda=kQ/I$ be a basic special biserial algebra, and let $M(S)$ (resp. $M(T)$) be a string module with canonical $k$-basis $\{x_u\}_{u=0}^m$ (resp. $\{y_v\}_{v=0}^n$) relative to the representative $S$ (resp. $T$).
Suppose $C$ is a string such that
\begin{enumerate}
\item[i.] $S\sim_s S_1CS_2$ with ($S_1$ of length $0$ or $S_1=S_1'\zeta_1$) and ($S_2$ of length $0$ or $S_2=\zeta_2^{-1} S_2'$), where $S_1,S_1',S_2,S'_2$ are strings and $\zeta_1,\zeta_2$ are arrows in $Q$; and
\item[ii.] $T\sim_sT_1CT_2$ with ($T_1$ of length $0$ or $T_1=T'_1\xi_1^{-1} $) and ($T_2$ of length $0$ or $T_2=\xi_2 T'_2$), where $T_1,T'_1,T_2,T'_2$ are strings and $\xi_1,\xi_2$ are arrows in $Q$.
\end{enumerate} 
Then there exists a $\Lambda$-module homomorphism $\sigma_C:M(S)\to M(T)$ which factors through $M(C)$ and which sends 
each element of $\{x_u\}_{u=0}^m$ either to zero or to an element of $\{y_v\}_{v=0}^n$,
according to the relative position of $C$ in $S$ and $T$, respectively.
If e.g. $S=s_1s_2\cdots s_m$, $T=t_1t_2\cdots t_n$, and $C=s_{i+1}s_{i+2}\cdots s_{i+\ell} = t_{j+\ell}^{-1}t_{j+\ell-1}^{-1}\cdots t_{j+1}^{-1}$,
then 
$$\sigma_C(x_{i+t})=y_{j+\ell-t} \mbox{ for } 0\le t\le \ell, \mbox{ and } \sigma_C(x_u)=0\mbox{ for all other $u$.}$$
Note that there may be several choices of $S_1,S_2$ (resp. $T_1,T_2$) in (i) (resp. (ii)). In other words, there may be more than one homomorphism factoring through $M(C)$. By \cite{krau}, every $\Lambda$-module homomorphism $\sigma:M(S)\to M(T)$ is a $k$-linear combination of homomorphisms which factor through string modules corresponding to strings $C$ satisfying (i) and (ii). 
\end{subrem}

The following definition provides a short-hand notation for homomorphisms between string modules, relative to fixed choices of canonical $k$-bases.

\begin{subdfn}
\label{def:helpmor}
Let $\Lambda=kQ/I$ be a basic special biserial algebra. Suppose $X=M(S)$ (resp. $Y=M(T)$) is a string module with canonical $k$-basis $\{x_u\}_{u=0}^m$ (resp. $\{y_v\}_{v=0}^n$) relative to the representative $S$ (resp. $T$). 
\begin{enumerate}
\item[i.] Suppose there exist $0\le i\le m$, $0\le j\le n$ and $0\le \ell\le \mathrm{min}\{m-i,j\}$ such that $\alpha:X\to Y$ defined by
$$\alpha(x_{i+t})=y_{j-\ell+t} \mbox{ for } 0\le t\le \ell, \mbox{ and } \alpha(x_u)=0 
\mbox{ for all other $u$}$$
is a $\Lambda$-module homomorphism. Then we denote $\alpha$ by
$\mathrm{hom}_{(X,Y)}^{++}(x_i,y_j,\ell)$.
\item[ii.] Suppose there exist $0\le i\le m$, $0\le j\le n$ and $0\le \ell\le \mathrm{min}\{m-i,n-j\}$ such that $\beta:X\to Y$ defined by
$$\beta(x_{i+t})=y_{j+\ell-t} \mbox{ for } 0\le t\le \ell, \mbox{ and } \beta(x_u)=0 
\mbox{ for all other $u$}$$
is a $\Lambda$-module homomorphism. Then we denote $\beta$ by
$\mathrm{hom}_{(X,Y)}^{+-}(x_i,y_j,\ell)$.
\item[iii.] Suppose there exist $0\le i\le m$, $0\le j\le n$ and $0\le \ell\le \mathrm{min}\{i,j\}$ such that $\gamma:X\to Y$ defined by
$$\gamma(x_{i-t})=y_{j-\ell+t} \mbox{ for } 0\le t\le \ell, \mbox{ and } \gamma(x_u)=0 
\mbox{ for all other $u$}$$
is a $\Lambda$-module homomorphism. Then we denote $\gamma$ by
$\mathrm{hom}_{(X,Y)}^{-+}(x_i,y_j,\ell)$.
\item[iv.] Suppose there exist $0\le i\le m$, $0\le j\le n$ and $0\le \ell\le \mathrm{min}\{i,n-j\}$ such that $\delta:X\to Y$ defined by
$$\delta(x_{i-t})=y_{j+\ell-t} \mbox{ for } 0\le t\le \ell, \mbox{ and } \delta(x_u)=0 
\mbox{ for all other $u$}$$
is a $\Lambda$-module homomorphism. Then we denote $\delta$ by
$\mathrm{hom}_{(X,Y)}^{--}(x_i,y_j,\ell)$.
\end{enumerate}
If $X=Y$ then we write $\mathrm{end}_X$ instead of $\mathrm{hom}_{(X,Y)}$.
\end{subdfn}

Note that in many of our applications of Definition \ref{def:helpmor}, $x_i$ corresponds to a source in the linear quiver $Q_S$ defining $M(S)$ and $y_j$ corresponds to a sink in the linear quiver $Q_T$ defining $M(T)$.

We now illustrate cases (i) and (ii) of Definition \ref{def:helpmor}. 
Suppose the string $S$ (resp. $T$) is the word $S=s_1s_2\cdots s_m$ (resp. $T=t_1t_2\cdots t_n$). Then $\alpha=\mathrm{hom}_{(X,Y)}^{++}(x_i,y_j,\ell)$ factors as 
$$\xymatrix @-1.2pc {
\alpha:&X=M(S)\ar[rr]^(.375){\pi}&& M(s_{i+1}s_{i+2}\cdots s_{i+\ell}) \ar[rr]^(.475){\cong}
&&M(t_{j-\ell+1}t_{j-\ell+2} \cdots t_j)\ar[rr]^(.6){\iota}&&M(T)=Y
}$$
where $s_{i+1}s_{i+2}\cdots s_{i+\ell}=t_{j-\ell+1}t_{j-\ell+2} \cdots t_j$ as words, and $\pi$ is the canonical projection and $\iota$ the canonical injection.
On the other hand, $\beta=\mathrm{hom}_{(X,Y)}^{+-}(x_i,y_j,\ell)$ factors as
$$\xymatrix @-1.2pc {
\beta:&X=M(S)\ar[rr]^(.375){\pi}&& M(s_{i+1}s_{i+2}\cdots s_{i+\ell}) \ar[rr]^(.475){\cong}
&&M(t_{j+\ell}^{-1} t_{j+\ell-1}^{-1}\cdots t_{j+1}^{-1})\ar[rr]^(.6){\iota}&&M(T)=Y
}$$
where $s_{i+1}s_{i+2}\cdots s_{i+\ell}=t_{j+\ell}^{-1} t_{j+\ell-1}^{-1}\cdots t_{j+1}^{-1}$ as words.

\medskip
 
The following is an easy combinatorial Lemma which helps determine stable homomorphisms. 
 
\begin{sublem}
\label{lem:needthis} 
Let $\Lambda=kQ/I$ be a basic special biserial algebra.
Suppose $0<\mu\in\mathbb{Z}$ and $0\le a<\mu$.
Let $M(S)$ and $M(T)$ be  string modules with canonical $k$-bases $\{x_u\}$ and $\{y_v\}$ relative to the representatives $S$ and $T$, respectively. Suppose that $\{h_i\}_{i=1}^s$ $($resp. $\{f_j\}_{j=1}^{t}$$)$ are subsets of $\{x_u\}$ $($resp. $\{y_v\}$$)$. Let $\epsilon\in\{\pm 1, 0\}$ be the sign of $(t-s)$. For $1\le i\le s$ and $1\le j\le t$ with $j-i\equiv a\mod\mu$, assume that the map $$\lambda_{i,j}:M(S)\to M(T)\mbox{ defined by $\lambda_{i,j}(h_i)=f_j$
and $\lambda_{i,j}(x_u)=0$ for $x_u\neq h_i$}$$
is a $\Lambda$-module homomorphism.
Suppose that for all $1\le i\le s-1$ and $1\le j\le t-1$ with $j-i\equiv a\mod\mu$,
$\alpha_{i,j}=\lambda_{i,j}+\lambda_{i+1,j+1}$ factors through a projective $\Lambda$-module. Suppose further that either
\begin{enumerate}
\item[i.] $\epsilon\ge 0$, and each $\lambda_{1,j}$ $($$1+\delta_1\le j\le t$, $j-1\equiv a\mod\mu$$)$ and each $\lambda_{s,j}$ $($$1\le j\le t-\delta_2$, $j-s\equiv a\mod \mu$$)$ factors through a projective $\Lambda$-module for $\{\delta_1,\delta_2\}=\{1,\epsilon\}$; or
\item[ii.] $\epsilon\le 0$, and each $\lambda_{i,1}$ $($$1-\delta_1\le i\le s$, $1-i\equiv a\mod \mu$$)$ and each $\lambda_{i,t}$ $($$1\le i\le s+\delta_2$, $t-i\equiv a\mod \mu$$)$ factors through a projective $\Lambda$-module for $\{\delta_1,\delta_2\}=\{-1,\epsilon\}$.
\end{enumerate}
Then for all $1\le i\le s$ and $1\le j\le t$ with $j-i\equiv a\mod\mu$, $\lambda_{i,j}$ factors through a projective $\Lambda$-module. 
\end{sublem}


\subsection{Components containing modules with endomorphism ring $k$}
\label{ss:kendo}
\setcounter{equation}{0}
\setcounter{figure}{0}

In this subsection, we prove Lemmas \ref{lem:psl1endos0}, \ref{lem:psl1ci} and \ref{lem:psl1uni}.
Let  $\Lambda=kQ/I$ with $Q$ and $I$ as in $\S\ref{ss:psl1}$. 

\medskip

\noindent
\textit{Proof of Lemma $\ref{lem:psl1endos0}$.}
Let $\mathcal{C}_0$ be the component of the stable Auslander-Reiten quiver of $\Lambda$ containing $S_0$, and let $M$ be a $\Lambda$-module belonging to $\mathcal{C}_0\cup\Omega(\mathcal{C}_0)$. We need to show that $\underline{\mathrm{End}}_\Lambda(M)=k$ and $\mathrm{Ext}^1_\Lambda(M,M)=0$. 

\begin{figure}[ht] \hrule \caption{\label{fig:a1n} The module $A_{1,n}$.}
$$\def\objectstyle{\scriptstyle} \def\labelstyle{\scriptstyle}
\vcenter{\xymatrix @-1.2pc @C=-.001pc{
&&&0
\ar[ld]_{\gamma}\ar[rd]^{\delta}&&&&\\
&&1\ar[ld]_{\beta}&&2\ar[rd]^{\eta}\\ 
&0
&&&&0\ar[rd]^{\gamma}\\
&&&&&&1\ar[rd]^{\beta}\\
&&&&&&&0\ar@{.}[rd]\\
&&&&&&&&0\ar[rd]_{\delta}&&&&0
\ar[ld]_{\gamma}\ar[rd]^\delta\\
&&&&&&&&&2\ar[rd]_{\eta}&&1\ar[ld]^{\beta}&&2\ar[rd]^{\eta}\\ 
&&&&&&&&&&0
&&&&0\ar[rd]^{\gamma}\\
&&&&&&&&&&&&&&&1\ar@{.}[rd]\\
&&&&&&&&&&&&&&&&0\ar[rd]_{\delta}&&&&&&0
\ar[ld]_{\gamma}\ar[rd]^\delta\\
&&&&&&&&&&&&&&&&&2\ar[rd]_{\eta}&&&&1\ar[ld]^{\beta}&&2\ar[rd]^{\eta}\\
&&&&&&&&&&&&&&&&&&0
&\cdots\cdots&0
&&&&0\ar[rd]^{\gamma}\\
&&&&&&&&&&&&&&&&&&&&&&&&&1\ar@{.}[rd]\\
&&&&&&&&&&&&&&&&&&&&&&&&&&0\ar[rd]^{\delta}\\
&&&&&&&&&&&&&&&&&&&&&&&&&&&2\ar[rd]^{\eta}\\
&&&&&&&&&&&&&&&&&&&&&&&&&&&&0\ar[rd]^{\gamma}\\
&&&&&&&&&&&&&&&&&&&&&&&&&&&&&1
 }}$$
\hrule
\end{figure}
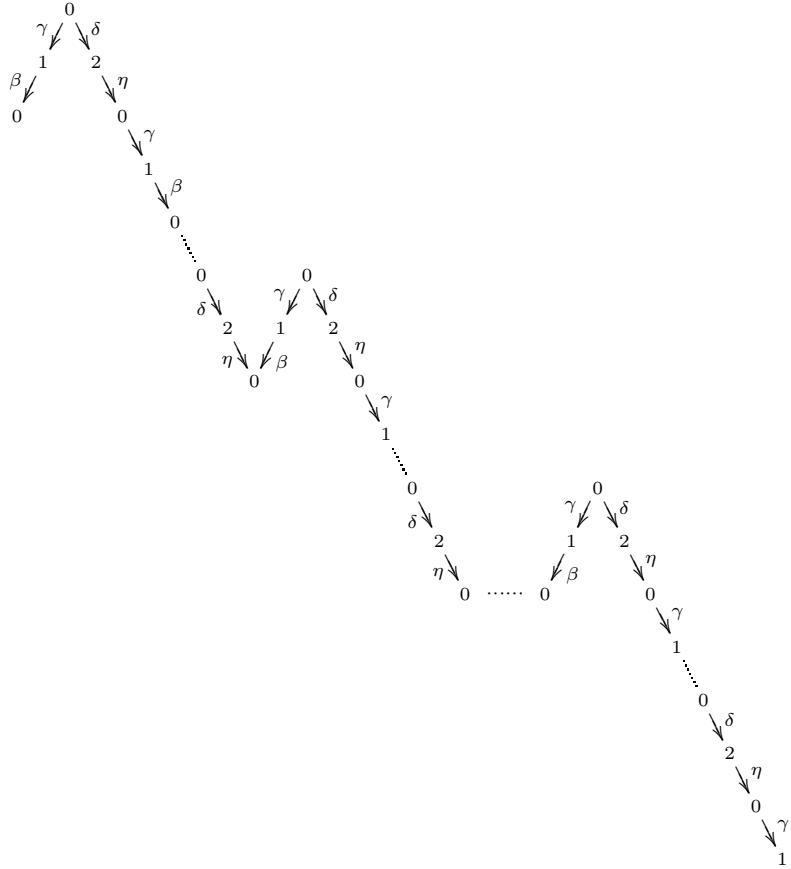

Using the description of the components of the stable Auslander-Reiten quiver of $\Lambda$ as in \S\ref{ss:arcomps}, we see that $\mathfrak{C}_0$ is of type $\mathbb{Z}A_\infty^\infty$. Moreover, using hooks and cohooks (see \S\ref{ss:arcomps}) we obtain the following.
There is an $i\in\mathbb{Z}$ such that $\Omega^i(M)$ is isomorphic to one of the following string modules: 
$$\begin{array}{c@{\;=\;}l}
\multicolumn{2}{c}{S_0, M(\beta), M(\eta), \mbox{ or for $n\ge 1$}}\\[2ex]
A_{1,n}&M\left(\left(\beta\gamma(\delta^{-1}\eta^{-1}\gamma^{-1}\beta^{-1})^{2^{d-2}-1}\delta^{-1}\eta^{-1}\right)^n\gamma^{-1}\right),\\[2ex]
A_{2,n}&M\left(\left(\beta\gamma(\delta^{-1}\eta^{-1}\gamma^{-1}\beta^{-1})^{2^{d-2}-1}\delta^{-1}\eta^{-1}\right)^n\right),\\[2ex]
A_{3,n}&M\left(\left(\beta\gamma(\delta^{-1}\eta^{-1}\gamma^{-1}\beta^{-1})^{2^{d-2}-1}\delta^{-1}\eta^{-1}\right)^n\beta\right),\\[2ex]
A'_{1,n}&M\left(\left(\eta\delta(\gamma^{-1}\beta^{-1}\delta^{-1}\eta^{-1})^{2^{d-2}-1}\gamma^{-1}\beta^{-1}\right)^n\delta^{-1}\right),\\[2ex]
A'_{2,n}&M\left(\left(\eta\delta(\gamma^{-1}\beta^{-1}\delta^{-1}\eta^{-1})^{2^{d-2}-1}\gamma^{-1}\beta^{-1}\right)^n\right),\\[2ex]
A'_{3,n}&M\left(\left(\eta\delta(\gamma^{-1}\beta^{-1}\delta^{-1}\eta^{-1})^{2^{d-2}-1}\gamma^{-1}\beta^{-1}\right)^n\eta\right).
\end{array}$$
It is straightforward to check that $\underline{\mathrm{End}}_\Lambda(M)=k$ and $\mathrm{Ext}^1_\Lambda(M,M)=0$ for $M\in\{S_0, M(\beta), M(\eta)\}$. We now demonstrate how to show this for $M=A_{1,n}$ for $n\ge 1$. The other cases are proved similarly.
The module $A_{1,n}$ (or more precisely the quiver defining it) is given in Figure \ref{fig:a1n}.
Let $\{x_r\}_{r=0}^{n2^d+1}$ be the corresponding canonical $k$-basis for $A_{1,n}$  (see \S\ref{sss:strings}). 
Then there are $n$ sources $h_1,\ldots,h_n$ in the quiver of $A_{1,n}$, corresponding to direct summands of $\mathrm{top}(A_{1,n})$, namely $h_i=x_{(i-1)2^d+2}$ for $1\le i\le n$. There are $n+1$ sinks $f_1,\ldots,f_{n+1}$ in the quiver of $A_{1,n}$, corresponding to direct summands of $\mathrm{soc}(A_{1,n})$, namely $f_j=x_{(j-1)2^d}$ for $1\le j\le n$ and $f_{n+1}=x_{n2^d+1}$.
By Remark \ref{rem:stringhoms}, each endomorphism of $A_{1,n}$ is a $k$-linear combination of
the elements of
\begin{equation}
\label{eq:endosa1n}
\{\mathrm{id}_{A_{1,n}}\}\;\bigcup\;\{\lambda_{i,j},\mu_i,\rho_{i,\ell}^s,
\sigma_i^s\;|\;1\le i,j\le n,2\le \ell\le n,1\le s\le 2^{d-2}-1\}
\end{equation}
where each of these endomorphisms is defined as follows, using Definition \ref{def:helpmor}:
\begin{eqnarray}
\label{eq:maybe1}
\lambda_{i,j}&=&\mathrm{end}_{A_{1,n}}^{++}(h_i,f_j,0);\\
\mu_i&=&\mathrm{end}_{A_{1,n}}^{-+}(h_i,f_{n+1},1);  \nonumber\\
\rho_{i,\ell}^s&=&\mathrm{end}_{A_{1,n}}^{++}(h_i,f_\ell,4s-2);  \nonumber\\
\sigma_i^s&=&\mathrm{end}_{A_{1,n}}^{++}(h_i,f_{n+1},4s-1).  \nonumber
\end{eqnarray}
Let $\alpha_{i,j}=\lambda_{i,j}+\lambda_{i+1,j+1}$ for $1\le i,j\le n-1$. Then it follows that $\alpha_{i,j}$ ($1\le i,j\le n-1$), $\lambda_{1,j}$ ($2\le j\le n$), and $\lambda_{n,j}$ ($1\le j\le n$) each factors through the projective $\Lambda$-module $P_0$. Hence by Lemma \ref{lem:needthis}, $\lambda_{i,j}$ ($1\le i,j\le n$) factors through a projective $\Lambda$-module.
Let $\beta_i=\lambda_{i-1,n}+\mu_i$ for $2\le i\le n$. Then it follows that $\beta_i$ ($2\le i\le n$) and $\mu_1$ each factors through $P_0$. Hence $\mu_i$ ($1\le i\le n$) factors through a projective $\Lambda$-module.
It also follows that $\rho_{i,\ell}^s$ ($1\le i\le n, 2\le \ell\le n, 1\le s\le 2^{d-2}-1$), and $\sigma_i^s$ ($1\le i\le n, 1\le s\le 2^{d-2}-1$) each factors through $P_0$. Hence $\underline{\mathrm{End}}_\Lambda(A_{1,n})=k$ for all $n\ge 1$.

Using that $\mathrm{Ext}^1_{\Lambda}(A_{1,n},A_{1,n})=\underline{\mathrm{Hom}}_\Lambda(\Omega(A_{1,n}),A_{1,n})$ and analyzing the quiver defining $\Omega(A_{1,n})$, one shows in a similar fashion that $\mathrm{Ext}^1_{\Lambda}(A_{1,n},A_{1,n})=0$ for all $n\ge 1$.
This completes the proof of Lemma \ref{lem:psl1endos0}.

\medskip

\noindent
\textit{Proof of Lemma $\ref{lem:psl1ci}$.}
Let $i\in\{1,2\}$ and let $\mathcal{C}_i$ be the component of the stable Auslander-Reiten quiver of $\Lambda$ containing $S_i$. We need to show that $\mathcal{C}_i$ is a $3$-tube with $S_i$ belonging to its boundary, and $\mathcal{C}_i=\Omega(\mathcal{C}_i)$. Moreover, if $M$ belongs to $\mathcal{C}_i$ and has stable endomorphism ring $k$, we need to show that $M$ is in the $\Omega$-orbit of $S_i$ and $\mathrm{Ext}^1_\Lambda(M,M)=0$.

We prove this for $i=1$. (The case $i=2$ is treated similarly.)
The component $\mathfrak{C}_1$ is a $3$-tube with boundary consisting of
$S_1$, $\Omega^2(S_1)=P_1/S_1$, and $\Omega^4(S_1)=\Omega(S_1)$. Since $\Omega$ maps the boundary to itself, it follows that $\mathfrak{C}_1=\Omega(\mathfrak{C}_1)$. It is straightforward to check that $\underline{\mathrm{End}}_{\Lambda}(M)=k$ and $\mathrm{Ext}^1_{\Lambda}(M,M)=0$ for $M$ belonging to the boundary of $\mathfrak{C}_1$. It remains to show that all other modules $X$ belonging to $\mathfrak{C}_1$ but not to its boundary have stable endomorphism ring of $k$-dimension larger than $1$.

Using hooks and cohooks (see \S\ref{ss:arcomps}), we obtain the following.
If $X$ belongs to $\mathfrak{C}_1$ but not to its boundary, then there exist integers $j$, $n$ with $n\ge 1$ such that $\Omega^j(X)$ is isomorphic to one of the following string modules: 
$$\begin{array}{c@{\;=\;}l}
X_{1,n}&M\left(\gamma(\delta^{-1}\eta^{-1}\gamma^{-1}\beta^{-1})^{2^{d-2}-1}\delta^{-1}\eta^{-1}\left(\beta\gamma(\delta^{-1}\eta^{-1}\gamma^{-1}\beta^{-1})^{2^{d-2}-1}\delta^{-1}\eta^{-1}\right)^{n-1}\gamma^{-1}\right),\\[2ex]
X_{2,n}&M\left(\gamma(\delta^{-1}\eta^{-1}\gamma^{-1}\beta^{-1})^{2^{d-2}-1}\delta^{-1}\eta^{-1}\left(\beta\gamma(\delta^{-1}\eta^{-1}\gamma^{-1}\beta^{-1})^{2^{d-2}-1}\delta^{-1}\eta^{-1}\right)^{n-1}\right),\\[2ex]
X_{3,n}&M\left(\gamma(\delta^{-1}\eta^{-1}\gamma^{-1}\beta^{-1})^{2^{d-2}-1}\delta^{-1}\eta^{-1}\left(\beta\gamma(\delta^{-1}\eta^{-1}\gamma^{-1}\beta^{-1})^{2^{d-2}-1}\delta^{-1}\eta^{-1}\right)^{n-1}\beta\right).
\end{array}$$
Let $\{x_r^{1,n}\}_{r=0}^{n2^d}$ be the canonical $k$-basis for $X_{1,n}$. Then $\mathrm{end}_{X_{1,n}}^{-+}(x_1^{1,n},x_{n2^d}^{1,n},1)$ does not factor through any projective $\Lambda$-module.
Let $\{x_r^{2,n}\}_{r=0}^{n2^d-1}$ (resp. $\{x_r^{3,n}\}_{r=0}^{n2^d}$) be the canonical $k$-basis for $X_{2,n}$ (resp. $X_{3,n}$). Then, for $\ell\in\{2,3\}$, $\mathrm{end}_{X_{\ell,n}}^{++}(x_1^{\ell,n},x_{n2^d-1}^{\ell,n},0)$
does not factor through any projective $\Lambda$-module.
This completes the proof of Lemma \ref{lem:psl1ci}.

\medskip

\noindent
\textit{Proof of Lemma $\ref{lem:psl1uni}$.}
Let $\mathfrak{C}$ be a component of the stable Auslander-Reiten quiver of $\Lambda$ containing a uniserial module $X$ of length $4$. We need to show that $\mathfrak{C}$ and $\Omega(\mathfrak{C})$ are both of type $\mathbb{Z}A_\infty^\infty$, and $\mathfrak{C}=\Omega(\mathfrak{C})$ exactly when $d=3$. Moreover, if $M$ belongs to $\mathfrak{C}\cup\Omega(\mathfrak{C})$ and has stable endomorphism ring $k$, we need to show that $M$ is in the $\Omega$-orbit of $X$ and $\mathrm{Ext}^1_\Lambda(M,M)=k$.

There are four uniserial $\Lambda$-modules of length $4$:
$$X_1=\begin{array}{c}1\\0\\2\\0\end{array}, \;
\Omega^2(X_1)= \begin{array}{c}0\\2\\0\\1\end{array}, \;
X_2=\begin{array}{c}2\\0\\1\\0\end{array},\;
\Omega^2(X_2)= \begin{array}{c}0\\1\\0\\2\end{array}.$$
Let $\mathfrak{C}$ be the component containing $X_1$. (The case when $\mathfrak{C}$ contains $X_2$ is treated similarly.)
We have seen in Lemma \ref{lem:psl1} that $X_1$ has endomorphism ring $k$. Moreover, 
$\Omega(X_1)$ is uniserial of length $2^d-3$ with descending radical series $(S_1,S_0,S_2,S_0,\ldots,S_1,S_0,S_2,S_0,S_1)$,
and thus
$\mathrm{Ext}^1_\Lambda(X_1,X_1)=
\underline{\mathrm{Hom}}_\Lambda(\Omega(X_1),X_1)=k$.

Using the description of the components of the stable Auslander-Reiten quiver of $\Lambda$ as in \S\ref{ss:arcomps}, we see that $\mathfrak{C}$ is of type $\mathbb{Z}A_\infty^\infty$. 
Moreover, using hooks and cohooks we can describe the modules $M$ in $\mathfrak{C}$ not belonging to the $\Omega$-orbit of $X_1$. From this description it follows that $\Omega(X_1)$ lies in $\mathfrak{C}$, and hence $\mathfrak{C}=\Omega(\mathfrak{C})$, if and only if $d=3$. Similarly to the proof of Lemma \ref{lem:psl1ci} we then construct for each remaining $M$ an endomorphism which does not factor through any projective $\Lambda$-module.
This completes the proof of Lemma \ref{lem:psl1uni}.


\subsection{Components of type $\mathbb{Z}A_\infty^\infty$}
\label{ss:infinity}
\setcounter{equation}{0}
\setcounter{figure}{0}

We next consider all components of the stable Auslander-Reiten quiver of type $\mathbb{Z}A_\infty^\infty$. We prove that the only such components containing a module with stable endomorphism ring $k$ are precisely those components containing a module $M$ with $\mathrm{End}_\Lambda(M)=k$ or $\mathrm{End}_\Lambda(\Omega(M))=k$. 

\begin{figure}[ht] \hrule \caption{\label{fig:arbcomp} The stable Auslander-Reiten component near $M(S)$.}
$$\xymatrix @-1.2pc{
&&&&\\
&M({}_{c'}S)\ar[rd]\ar@{.}[ld]\ar@{.}[lu]\ar@{.}[ru]&&M(S_{h'})\ar@{.}[rd]\ar@{.}[lu]\ar@{.}[ru]&\\
&&M(S)\ar[ru]\ar[rd]&&\\
&M(S_{c'})\ar[ru]\ar@{.}[rd]\ar@{.}[ld]\ar@{.}[lu]&&M({}_{h'}S)\ar@{.}[rd]\ar@{.}[ld]\ar@{.}[ru]&\\
&&&&
}$$
\hrule
\end{figure}
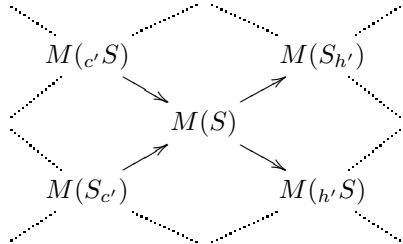

We use the notation of hooks and cohooks from \S\ref{ss:arcomps}. Moreover, if $S$ is a string of positive length in a $\mathbb{Z}A_\infty^\infty$-component, we use  $S_{h'}$ (resp. ${}_{h'}S$) to denote the string obtained from $S$ by either adding a hook or subtracting a cohook on the right side (resp. left side) of $S$. We also use $S_{c'}$ (resp. ${}_{c'}S$) to denote the string obtained from $S$ by either adding a cohook or subtracting a hook on the right side (resp. left side) of $S$. This means that near $M(S)$ the stable Auslander-Reiten component looks 
as in Figure \ref{fig:arbcomp}.
Note that if $S$ has minimal length in its stable Auslander-Reiten component then $S_{h'}=S_h$, ${}_{h'}S={}_hS$, $S_{c'}=S_c$ and ${}_{c'}S={}_cS$. If none of the projective $\Lambda$-modules is uniserial and $S$ has minimal length, then also $S_{h'\cdots h'}=S_{h\cdots h}$, ${}_{h'\cdots h'}S={}_{h\cdots h}S$, $S_{c'\cdots c'}=S_{c\cdots c}$ and ${}_{c'\cdots c'}S={}_{c\cdots c}S$.
The following Lemma is straightforward.

\begin{sublem}
\label{lem:upsidedown}
Let $\Lambda=kQ/I$ be a symmetric special biserial algebra with the following properties:
\begin{enumerate}
\item[a.] the quiver $Q$ contains no double arrows;
\item[b.] for all $v_1,v_2\in Q_0$ and $\alpha\in Q_1$ with $s(\alpha)=v_1$ and $e(\alpha)=v_2$, there exists $\tau(\alpha)\in Q_1$ with $s(\tau(\alpha))=v_2$ and $e(\tau(\alpha))=v_1$;
\item[c.] if $P$ is a uniserial projective $\Lambda$-module, then $P$ is a string module corresponding to a directed string $\alpha_1\alpha_2\cdots \alpha_\ell$, and
$$\alpha_1\alpha_2\cdots\alpha_{\ell} = \tau(\alpha_{\ell})\tau(\alpha_{\ell-1})\cdots \tau(\alpha_1);$$
if $P$ is a non-uniserial projective $\Lambda$-module, then $P$ corresponds to a relation $p_1=p_2$ in $I$ for two paths $p_1=\alpha_1\alpha_2\cdots\alpha_{\ell_1}$ and $p_2=\beta_1\beta_2\cdots\beta_{\ell_2}$ in $kQ$, and
$$\{\alpha_1\alpha_2\cdots\alpha_{\ell_1},\;\beta_1\beta_2\cdots\beta_{\ell_2}\} = \{\tau(\alpha_{\ell_1})\tau(\alpha_{\ell_1-1})\cdots\tau(\alpha_1),\;\tau(\beta_{\ell_2})\tau(\beta_{\ell_2-1})\cdots \tau(\beta_1)\}.$$
\end{enumerate}
Let $S=w_1w_2\cdots w_n$ be a string of length $n\ge 1$, and define $\tau(S)=\tau(w_1)^{-1}\tau(w_2)^{-1}\cdots \tau(w_n)^{-1}$. Then $\tau(S)$ is a string.
Moreover, suppose $\{x_r\}_{r=0}^n$ $($resp. $\{y_r\}_{r=0}^n$$)$ is the canonical $k$-basis of $M(S)$ $($resp. $M(\tau(S))$$)$ relative to the representative $S$ $($resp. $\tau(S)$$)$. Let $\epsilon_1,\epsilon_2\in\{+,-\}$, and let $0\le u,v,\ell\le n$. Define $\rho(+)=-$ and $\rho(-)=+$. If $\mathrm{end}_{M(S)}^{\epsilon_1\, \epsilon_2}(x_u,x_v,\ell)$ is a $\Lambda$-endomorphism of $M(S)$, then $\mathrm{end}_{M(\tau(S))}^{\rho(\epsilon_2)\, \rho(\epsilon_1)}(y_v,y_u,\ell)$ is a $\Lambda$-endomorphism of $M(\tau(S))$.
Moreover, $\mathrm{end}_{M(S)}^{\epsilon_1\, \epsilon_2}(x_u,x_v,\ell)$ factors through a projective $\Lambda$-module $P$ if and only if $\mathrm{end}_{M(\tau(S))}^{\rho(\epsilon_2)\, \rho(\epsilon_1)}(y_v,y_u,\ell)$ factors through $P$.
\end{sublem}

We are now able to prove the following result.
\begin{subprop}
\label{prop:endospsl1}
Let $\Lambda=kQ/I$ where $Q$ and $I$ are as in $\S\ref{ss:psl1}$. Then the components of the stable Auslander-Reiten quiver of type $\mathbb{Z}A_\infty^\infty$ containing a module with stable endomorphism ring $k$ are precisely the components containing a module $M$ with $\mathrm{End}_\Lambda(M)=k$ or $\mathrm{End}_\Lambda(\Omega(M))=k$.
\end{subprop}

\begin{proof}
Let $\mathfrak{C}$ be a component of type $\mathbb{Z}A_\infty^\infty$ of the stable Auslander-Reiten quiver of $\Lambda$ such that $\mathfrak{C}\cup\Omega(\mathfrak{C})$ contains no simple $\Lambda$-module and no uniserial $\Lambda$-module of length $4$, i.e. by Lemma \ref{lem:psl1},  $\mathfrak{C}\cup\Omega(\mathfrak{C})$ contains no $\Lambda$-module with endomorphism ring $k$.

Let $X$ be a $\Lambda$-module of minimal length in $\mathfrak{C}$. In particular, $X$ (or more precisely the string defining $X$) cannot start or end in a peak (resp. in a  deep).  Since $P_1$ and $P_2$ are uniserial, this means that $X$ cannot have any of the following forms
\begin{eqnarray}
\label{eq:oyvey1}
X=\begin{array}{c@{}c}1\\&0\; \cdots\end{array},\quad
X=\begin{array}{c@{}c}2\\&0\; \cdots\end{array},\quad
X=\begin{array}{c@{}c}&1\\ \cdots\; 0\end{array},\quad
X=\begin{array}{c@{}c}&2\\ \cdots\; 0\end{array},\\ 
X=\begin{array}{c@{}c}&0\; \cdots\\1\end{array},\quad
X=\begin{array}{c@{}c}&0\; \cdots\\2\end{array},\quad
X=\begin{array}{c@{}c} \cdots\; 0\\&1\end{array},\quad
X=\begin{array}{c@{}c} \cdots\; 0\\&2\end{array}.\nonumber
\end{eqnarray}
Let $\{z_r\}_{r=0}^{\ell(X)}$ be the canonical $k$-basis of the string module $X$ relative to the chosen representative.  

Suppose first $X$ is uniserial. By (\ref{eq:oyvey1}), $X=M(S)$ where $S$ is one of the following strings for $1\le n\le 2^{d-2}-1$:
\begin{eqnarray}
\label{eq:oyvey2}
S_{1,1,n}&=&\left(\gamma^{-1}\beta^{-1}\delta^{-1}\eta^{-1}\right)^n\gamma^{-1}\beta^{-1};\\
S_{1,2,n}&=&\left(\gamma^{-1}\beta^{-1}\delta^{-1}\eta^{-1}\right)^n;\nonumber\\
S_{2,2,n}&=&\left(\delta^{-1}\eta^{-1}\gamma^{-1}\beta^{-1}\right)^n\delta^{-1}\eta^{-1};\nonumber\\
S_{2,1,n}&=&\left(\delta^{-1}\eta^{-1}\gamma^{-1}\beta^{-1}\right)^n.\nonumber
\end{eqnarray}
Note that $M(S_{1,1,0})$ (resp. $M(S_{2,2,0})$) lies in the same component of the Auslander-Reiten quiver as a uniserial $\Lambda$-module of length $4$. 
If $i\neq j$ in $\{1,2\}$, then $\Omega(M(S_{i,i,2^{d-2}-2}))$ lies in the same Auslander-Reiten component as a  uniserial $\Lambda$-module of length $4$, $M(S_{i,i,2^{d-2}-1})$ lies in the same Auslander-Reiten component as $S_j$, and $\Omega(M(S_{i,j,2^{d-2}-1}))$ lies in the same Auslander-Reiten component as $S_0$.
On the other hand, if $n<2^{d-2}-2$, then $\mathrm{end}_{M(S_{i,i,n})}^{++}(z_0,z_{4n+2},4n-2)$ does not factor through a projective $\Lambda$-module for $i\in\{1,2\}$; and, if $n<2^{d-2}-1$, then $\mathrm{end}_{M(S_{i,j,n})}^{++}(z_0,z_{4n},4n-4)$ does not factor through a projective $\Lambda$-module for $i\neq j$ in $\{1,2\}$. By considering hooks and cohooks, we see that in both cases, all $\Lambda$-modules in the Auslander-Reiten component containing $X$ have stable endomorphism ring of dimension at least $2$.

Suppose now $X$ is not uniserial. By Lemma \ref{lem:upsidedown}, it is enough to consider $X$ where 
$$X=\begin{array}{c@{}c}0\\&1\; \cdots\end{array},\quad\mbox{ or }\quad
X=\begin{array}{c@{}c}0\\&2\; \cdots\end{array}.$$
Let $(\xi_1,\xi_2)=(\eta,\beta)$.
If $X=M(S)$ where $S=S_{i,i,n}\xi_i\cdots$ for $i\in\{1,2\}$ and $1\le n <2^{d-2}-1$, then $\mathrm{end}_{X}^{++}(z_0,z_{4n+2},4n-2)$ does not factor through any projective $\Lambda$-module, and the same is true for $M(S_{h'\cdots h'})$ and $M(S_{c'\cdots c'})$. Similarly, if $X=M(S)$ where $S=S_{i,j,n}\xi_j\cdots$ for $i\neq j$ in $\{1,2\}$ and $1\le n<2^{d-2}-1$, then $\mathrm{end}_{X}^{++}(z_0,z_{4n},4n-4)$ does not factor through any projective $\Lambda$-module, and the same is true for $M(S_{h'\cdots h'})$ and $M(S_{c'\cdots c'})$. If $X=M(S_{i,i,2^{d-2}-1}\xi_iC_i)$ where $i\in\{1,2\}$ and $C_i$ is some string, then $X$ lies in the same Auslander-Reiten component as $M(C_i)$. Hence $X$ is not of minimal length in $\mathfrak{C}$. The two remaining possibilities for $X$ are
$$X=M(S_{i,i,0}\xi_i\cdots) \mbox{ for $i\in\{1,2\}$, } \quad \mbox{ or } \quad
X=M(S_{i,j,2^{d-2}-1}\xi_j\cdots) \mbox{ for $i\neq j$ in $\{1,2\}$.}$$
Let $i\neq j$ in $\{1,2\}$. Then $\Omega(M(S_{i,j,2^{d-2}-1}\xi_j\cdots))$ lies in the same Auslander-Reiten component as a string module of the form $M(\tau(S_{i,i,0}\xi_i\cdots))$, where $\tau$ is as in Lemma \ref{lem:upsidedown}. Hence, by Lemma \ref{lem:upsidedown}, it is enough to consider $X$ of the form $X=M(S_{i,i,0}\xi_i\cdots)$ for $i\in\{1,2\}$. 

Let $i=1$. (The case $i=2$ is done similarly.) Since $X$ has minimal length in the Auslander-Reiten component $\mathfrak{C}$, this means $X=M(SC)$ where $S$ is one of the following strings:
\begin{eqnarray}
\label{eq:oyvey3}
S^{1,2}_n&=&\gamma^{-1}\beta^{-1}(S_{1,2,n})^{-1}, \quad1\le n\le 2^{d-2}-1;\\
S^{2,2}_n&=&\gamma^{-1}\beta^{-1}(S_{2,2,n})^{-1}, \quad0\le n\le 2^{d-2}-1; \nonumber
\end{eqnarray}
and $C$ is a string such that the following holds. If $C$ has positive length then $C=\zeta^{-1}\cdots$ for the appropriate arrow $\zeta$ in $Q$, and if $S=S^{2,2}_{2^{d-2}-1}$ then $C$ cannot have length $0$. 

If $S=S^{1,2}_n$, $1\le n\le 2^{d-2}-1$, then $\mathrm{end}_X^{++}(z_0,z_2,0)$ does not factor through a projective $\Lambda$-module. The same is true for $M((SC)_{c'\cdots c'})$ and for $M((SC)_{h'})$, and, if $n<2^{d-2}-1$ or $C$ has positive length, also for $M((SC)_{h'\cdots h'})$. If $SC=S=S^{1,2}_{2^{d-2}-1}$, then $(SC)_{h'\cdots h'}=S^{1,2}_{2^{d-2}-1}\eta\delta \gamma^{-1}\beta^{-1}\delta^{-1} \cdots$. Thus $\mathrm{end}_{M((SC)_{h'\cdots h'})}^{++}(z_{2^d},z_2,2)$ does not factor through a projective $\Lambda$-module. 

If $S=S^{2,2}_n$, $0\le n\le 2^{d-2}-2$, then $\mathrm{end}_X^{++}(z_0,z_2,0)$ does not factor through a projective $\Lambda$-module. The same is also true for $M((SC)_{c'\cdots c'})$ and for $M((SC)_{h'\cdots h'})$. 

Now suppose $S=S^{2,2}_{2^{d-2}-1}$. Since $X$ has minimal length in its Auslander-Reiten component, $C=\gamma^{-1}\beta^{-1}C'$ for some string $C'$. If $C'$ has length $0$ or $C'=\delta^{-1}\cdots$, then $\mathrm{end}_X^{++}(z_{2^d},z_2,2)$ does not factor through a projective $\Lambda$-module. The same is true for $M((SC)_{c'\cdots c'})$, and, if $C'$ has positive length, also for $M((SC)_{h'\cdots h'})$. If $C'$ has length $0$, i.e. $SC=S^{2,2}_{2^{d-2}-1}\gamma^{-1}\beta^{-1}$, then $(SC)_{h}=(SC)\eta$ and $(SC)_{h'\cdots h'}=
(SC)\eta\delta\gamma^{-1}\cdots$. Thus $\mathrm{end}_{M((SC)_{h})}^{++}(z_{2^d},z_3,3)$ and $\mathrm{end}_{M((SC)_{h'\cdots h'})}^{++}(z_{2^d},z_4,4)$
do not factor through a projective $\Lambda$-module.
Now suppose $C'=\eta\cdots$. Then $X$ has the form $X=M(S^{2,2}_{2^{d-2}-1}SC_1)$ where $S$ is one of the strings in (\ref{eq:oyvey3}) and $C_1$ has the same properties as the properties of $C$ described below (\ref{eq:oyvey3}). 
Hence one continues using similar arguments as above, and thus concludes that 
the Auslander-Reiten components containing any of these $X$ contain no $\Lambda$-module with stable endomorphism ring $k$.
\end{proof}


\subsection{One-tubes}
\label{ss:bands}
\setcounter{equation}{0}
\setcounter{figure}{0}

Finally, we consider the components of the stable Auslander-Reiten quiver which are $1$-tubes. We prove that no $1$-tube contains any modules with stable endomorphism ring $k$.
Since for the blocks in question all string modules lie either in components of type $\mathbb{Z}A_\infty^\infty$ or in $3$-tubes, all the modules in $1$-tubes are band modules. We use the description of band modules from \S\ref{ss:stringind}.

\begin{subrem}
\label{rem:krause}
Let $\Lambda$ be a special biserial algebra.
It follows from \cite{krau} that if $B$ is a band, $\lambda\in k-\{0\}$ and $n\ge 2$ is an  integer, then $\underline{\mathrm{End}}_\Lambda(M(B,\lambda,n))$ has dimension
at least $2$. To be more precise, the endomorphisms coming from the circular quiver associated to the band $B$ are parametrized by upper triangular $n\times n$ matrices with equal entries along each diagonal; and these endomorphisms cannot factor through a projective $\Lambda$-module. So if $n\ge 2$, then the $k$-dimension of the space spanned by these endomorphisms is $n\ge 2$.
\end{subrem}

\begin{subdfn}
\label{def:bandendohelp}
Let $\Lambda=kQ/I$ be a special biserial algebra. Let $B$ be a band for $\Lambda$, $\lambda\in k-\{0\}$ and let $M_{B,\lambda}=M(B,\lambda,1)$. 
Suppose $S$ is a string such that 
\begin{enumerate}
\item[i.] $B\sim_rST_1$ with $T_1=\xi_1^{-1} T'_1\xi_2$, where $T_1,T'_1$ are strings and $\xi_1,\xi_2$ are arrows in $Q$; and 
\item[ii.] $B\sim_r ST_2$ with $T_2=\zeta_1 T'_2\zeta_2^{-1}$, where $T_2,T'_2$ are strings and $\zeta_1,\zeta_2$ are arrows in $Q$.
\end{enumerate} 
Then by \cite{krau} there exists an endomorphism of $M_{B,\lambda}$ which factors through $M(S)$. We will call such an endomorphism to be of string type $S$. Note that there may be several choices of $T_1$ (resp. $T_2$) in (i) (resp. (ii)). In other words, there may be more than one endomorphism of string type $S$. By \cite{krau}, every endomorphism of $M_{B,\lambda}$ is a $k$-linear combination of the identity morphism and of endomorphisms of string type $S$ for suitable choices of strings $S$ satisfying (i) and (ii). 
\end{subdfn}

\begin{subdfn}
\label{def:piece}
Let $\Lambda=kQ/I$ be a special biserial algebra, and let $B$ be a band for $\Lambda$. 
\begin{enumerate}
\item[i.] 
We call a string $C$ a substring of $B$ if $B\sim_r  C C'$ for some string $C'$. 
\item[ii.]
A substring $S$ of $B$ is called a  top-socle piece of $B$ if 
\begin{enumerate}
\item[a.] $S=\alpha_\ell^{-1}\cdots \alpha_2^{-1}\alpha_1^{-1}$ for $\ell\ge 1$ and arrows $\alpha_1,\alpha_2,\ldots,\alpha_\ell$ in $Q$, and 
\item[b.] $B\sim_r ST$ for some string $T$ where $T=\xi T' \zeta$ and $\xi,\zeta$ are arrows in $Q$. 
\end{enumerate}
Note that $B\sim_r C_0C_1^{-1}C_2C_3^{-1}\cdots C_s^{-1}$ where $s\ge 1$ is odd and $C_0,C_1,\ldots,C_s$ are top-socle pieces of $B$. 
\item[iii.]
If $S$ is a top-socle piece of $B$ and $E$ is a simple $\Lambda$-module which is isomorphic to the top (resp. socle) of $M(S)$, we also say that $S$ has top (resp. socle) isomorphic to $E$.
\end{enumerate}
\end{subdfn}

In the proof of the following Proposition, we will often use the equality sign instead of the more precise $\sim_r$.

\begin{subprop}
\label{prop:bandspsl1}
Let $\Lambda=kQ/I$ where $Q$ and $I$ are as in $\S\ref{ss:psl1}$. Then each module $M$ which lies in a component of the stable Auslander-Reiten quiver which is a $1$-tube satisfies $\mathrm{dim}_k\,\underline{\mathrm{End}}_\Lambda(M)\ge 2$.
\end{subprop}

\begin{proof}
Let $B$ be a band for $\Lambda$ such that possibly $\underline{\mathrm{End}}_\Lambda(M(B,\lambda,n))=k$. By Remark \ref{rem:krause}, we only need to consider $M_{B,\lambda}=M(B,\lambda,1)$ for $\lambda\in k-\{0\}$.

It follows from the shape of the quiver $Q$ that the only top-socle pieces that can occur in $B$ must have top and socle isomorphic to $S_0$. This means they have the form, using the notation from (\ref{eq:oyvey2}):
$$S_{i,i,n} \mbox{ for } 0\le n\le 2^{d-2}-1,\quad\mbox{ and }
S_{i,j,n} \mbox{ for } 1\le n\le 2^{d-2}-1,\quad\mbox{ where $i\neq j$ in $\{1,2\}$.}$$
Let now $i\neq j$ in $\{1,2\}$.
If $S_{i,i,n}$ for $0< n < 2^{d-2}-1$ (resp.  $S_{i,j,n}$ for $1< n < 2^{d-2}-1$) is a top-socle piece of $B$, then $M_{B,\lambda}$ has an endomorphism of string type $S_{i,i,n-1}$ (resp. of string type $S_{i,j,n-1}$) which does not factor through a projective module. This implies that the only top-socle pieces that can occur in $B$ are (again using the notation from (\ref{eq:oyvey2}))
\begin{equation}
\label{eq:topsoclepsl1}
S_{i,i,0}, S_{i,i,2^{d-2}-1},\quad\mbox{ and } S_{i,j,1}, S_{i,j,2^{d-2}-1}\quad\mbox{ where $i\neq j$ in $\{1,2\}$.}
\end{equation}

Let $i\neq j$ in $\{1,2\}$, and suppose $S_{i,j,1}$ (resp. $S_{i,j,2^{d-2}-1}$) is a top-socle piece of $B$. If $S_{j,i,1}$ is also a top-socle piece of $B$, then $M_{B,\lambda}$ has an endomorphism of string type $S_{j,j,0}$ which does not factor through a projective module. On the other hand, if $S_{j,i,2^{d-2}-1}$ is also a top-socle piece of $B$, then $M_{B,\lambda}$ has an endomorphism of string type $S_{j,j,0}$ (resp. $S_{j,j,2^{d-2}-2}$) which does not factor through a projective module. This means that if $S_{i,j,1}$ (resp. $S_{i,j,2^{d-2}-1}$) is a top-socle piece of $B$, then neither $S_{j,i,1}$ nor $S_{j,i,2^{d-2}-1}$ can be a top-socle piece of $B$. 

We claim that this makes it impossible for either $S_{i,j,1}$ or $S_{i,j,2^{d-2}-1}$ to occur as top-socle piece of $B$. We demonstrate this for the case $(i,j)=(1,2)$. Let $\ell\in\{1,2^{d-2}-1\}$, and suppose $S_{i,j,\ell}$ is a top-socle piece of $B$. This means $B=S_{1,2,\ell}C$ for some string $C=C_1^{-1}C_2C_3^{-1}\cdots C_s^{-1}$ where $s\ge 1$ is odd and $C_1,C_2,\ldots,C_s$ are top-socle pieces of $B$. But then $C_u$ must be of the form $S_{1,1,1}$ or $S_{1,1,2^{d-2}-1}$ for odd $u$, and of the form $S_{2,2,1}$ or $S_{2,2,2^{d-2}-1}$ for even $u$. In particular, $C_s$ is of the form $S_{1,1,1}$ or $S_{1,1,2^{d-2}-1}$. But this contradicts that $B$ is a band, since $C_s^{-1}S_{1,2,\ell}$ is not a valid word. 

Hence, by (\ref{eq:topsoclepsl1}), the only top-socle pieces that can occur in $B$ are
\begin{equation}
\label{eq:topsoclepsl2}
S_{i,i,0}, S_{i,i,2^{d-2}-1}\quad\mbox{for $i\in\{1,2\}$.}
\end{equation}

Let $i\neq j$ in $\{1,2\}$, and suppose $S_{i,i,0}$ is a top-socle piece of $B$. Then $B=S_{i,i,0}C$ for some string $C=C_1^{-1}C_2C_3^{-1}\cdots C_s^{-1}$ where $s\ge 1$ is odd and $C_1,C_2,\ldots,C_s$ are top-socle pieces of $B$. 

We claim that $C_1$ and $C_s$ must both be $S_{j,j,2^{d-2}-1}$. This can be shown as follows. Suppose first $C_s=S_{j,j,0}$. Then it follows that $C_{s-1}=S_{i,i,2^{d-2}-1}$ and $C_1=S_{j,j,2^{d-2}-1}$, since otherwise $M_{B,\lambda}$ has an endomorphism of string type $1_0$ (i.e. the string of length $0$ corresponding to the vertex $0$) which does not factor through a projective module. But then, to ensure that $B$ is a band, one of two things must be true: Either there is an odd $u_0$ with $C_{u_0}^{-1}C_{u_0+1}=S_{j,j,2^{d-2}-1}^{-1}S_{i,i,2^{d-2}-1}$, in which case $M_{B,\lambda}$ has an endomorphism of string type $S_{j,j,0}^{-1}S_{i,i,0}$ which does not factor through a projective module. Or there is an odd $u_0$ with $C_{u_0-1}C_{u_0}^{-1}=S_{i,i,0}S_{j,j,0}^{-1}$, in which case $C_{u_0+1}=S_{i,i,2^{d-2}-1}$ which means $M_{B,\lambda}$ has an endomorphism of string type $S_{i,i,0}S_{j,j,0}^{-1}S_{i,i,0}$ which does not factor through a projective module. Hence we get a contradiction, which means $C_s=S_{j,j,2^{d-2}-1}$. The same argument shows that if $u<s$ is odd then at least one of $C_u$ and $C_{u+1}$ must be in $\{S_{1,1,2^{d-2}-1},S_{2,2,2^{d-2}-1}\}$. Suppose now that  $C_1=S_{j,j,0}$. Then we must have $C_2=S_{i,i,2^{d-2}-1}$. Moreover, there must exist an odd $u_0$ such that $C_{u_0-1}C_{u_0}^{-1}=S_{i,i,2^{d-2}-1}S_{j,j,2^{d-2}-1}^{-1}$. But this means that $M_{B,\lambda}$ has an endomorphism of string type $S_{i,i,0}S_{j,j,0}^{-1}$ which does not factor through a projective module. Thus $C_1$ and $C_s$ are both $S_{j,j,2^{d-2}-1}$. 

Since $B$ cannot be the power of a smaller word, we see that if $s>1$ then $M_{B,\lambda}$ has an endomorphism of string type $S_{i,i,0}S_{j,j,2^{d-2}-1}^{-1}S_{i,i,0}$ which does not factor through a projective module. On the other hand, if $s=1$ then $B=S_{i,i,0}S_{j,j,2^{d-2}-1}^{-1}$, and $M_{B,\lambda}$ has an endomorphism of string type $1_0$ (i.e. the string of length $0$ corresponding to the vertex $0$) which does not factor through a projective module.

Hence neither $S_{1,1,0}$ nor $S_{2,2,0}$ is a top-socle piece of $B$. Thus the only possible band is $B=S_{1,1,2^{d-2}-1}S_{2,2,2^{d-2}-1}^{-1}$. But then $M_{B,\lambda}$ has an endomorphism of string type $S_{1,2,2^{d-2}-1}$ which does not factor through a projective module. So in all cases, $\underline{\mathrm{End}}_\Lambda(M_{B,\lambda})$ has $k$-dimension at least $2$, which completes the proof of Proposition \ref{prop:bandspsl1}.
\end{proof}


\section{Background: Special biserial algebras}
\label{s:prelimstring}
\setcounter{equation}{0}
\setcounter{figure}{0}

In this section, we give a short introduction to special biserial algebras. For more background material, we refer to \cite{buri}.
Let $k$ be an algebraically closed field of characteristic $p>0$, let $Q$ be a finite quiver and let $I$ be an admissible ideal in the path algebra $kQ$.

\begin{dfn} 
\label{def:stringalg}
A finite dimensional basic $k$-algebra $\Lambda=kQ/I$ is called special biserial
if  the following conditions are satisfied:
\begin{enumerate}
\item[i.] Any vertex of $Q$ is starting point (resp. end point) of at most two arrows.
\item[ii.] For a given arrow $\beta$ in $Q$, there is at most one arrow $\gamma$  with $\beta\gamma \not\in I$, and there is at most one arrow $\alpha$ with $\alpha\beta\not\in I$.
\end{enumerate}
If additionally $I$ is generated by paths, $\Lambda$ is called a string algebra.
\end{dfn}

If $\Lambda$ is a special biserial algebra and $\mathcal{P}$ is a full set of representatives of the projective indecomposable $\Lambda$-modules which are also injective and not uniserial, then $\bar{\Lambda}=\Lambda/(\oplus_{P\in{\mathcal{P}}}\mathrm{soc}(P))$ is a string algebra. Furthermore, the indecomposable $\bar{\Lambda}$-modules are exactly the indecomposable $\Lambda$-modules which are not isomorphic to any $P\in{\mathcal{P}}$.


\subsection{Indecomposable modules for string algebras}
\label{ss:stringind}
\setcounter{equation}{0}
\setcounter{figure}{0}

Let $\Lambda=kQ/I$ be a basic string algebra. Then all indecomposable $\Lambda$-modules are either string or band modules (see e.g. \cite[\S3]{buri}). The definitions are as follows.

Given an arrow $\beta$ in $Q$ with starting point $s(\beta)$ and end point $e(\beta)$, denote by $\beta^{-1}$ a formal inverse of $\beta$. In particular, $s(\beta^{-1})=e(\beta)$, $e(\beta^{-1})=s(\beta)$, and $(\beta^{-1})^{-1}=\beta$.
A word $w$ is a sequence $w_1\cdots w_n$, where $w_i$ is either an arrow or a formal inverse such that $s(w_i)=e(w_{i+1})$ for $1\leq i \leq n-1$. Define $s(w)=s(w_n)$, $e(w)=e(w_1)$ and $w^{-1}=w_n^{-1}\cdots w_1^{-1}$. For each vertex $u$ in $Q$ there exists an empty word $1_u$ of length $0$ with $e(1_u)=s(1_u)=u$ and $(1_u)^{-1}=1_u$. Denote the set of all words by $\mathcal{W}$, and the set of all 
non-empty words $w$ with $e(w)=s(w)$ by ${\mathcal{W}}_r$. In the following, Greek letters inside words always denote arrows.

\subsubsection{Strings and string modules}
\label{sss:strings}
Let $\sim_s$ be the equivalence relation on $\mathcal{W}$ with $w\sim_s w'$ if and only if $w=w'$ or $w^{-1}=w'$. Then strings are representatives $w\in{\mathcal{W}}$ of the equivalence classes under $\sim_s$ with the following property: Either $w=1_u$ or $w=w_1\cdots w_n$, where $w_i \neq w_{i+1}^{-1}$ for $1\leq i\leq n-1$ and no subword of $w$ or its formal inverse belongs to $I$.

Let $C=w_1\cdots w_n$ be a string of length $n$ and let $Q_C$ be the linear quiver
$$Q_C = \xymatrix @-1.2pc {
\cdot \ar@{-}[r]^(.3){w_1} & \cdot\; \cdots  \;\cdot\ar@{-}[r]^(.7){w_n}&\cdot
}$$
where $\xymatrix @-1.2pc {\cdot \ar@{-}[r]^(.3){w_i}&\cdot \; = \; \cdot&\cdot\ar[l]_(.3){\beta}}$ if $w_i=\beta$ is an arrow, and
$\xymatrix @-1.2pc {\cdot \ar@{-}[r]^(.3){w_i}&\cdot \; = \; \cdot\ar[r]^(.7){\beta}&\cdot}$ if $w_i=\beta^{-1}$ is a formal inverse. 
Then the representation of $Q_C$ which assigns to each vertex the vector space $k$ and to each arrow the identity map defines an indecomposable $\Lambda$-module, called the string module $M(C)$ corresponding to the string $C$. 
More precisely, there is a $k$-basis $\{z_0,z_1,\ldots, z_n\}$ of $M(C)$ such that the action of $\Lambda$ on $M(C)$ is given by the following representation $\varphi_C:\Lambda\to\mathrm{Mat}(n+1,k)$. Let $v(i)=e(w_{i+1})$ for $0\leq i\leq n-1$ and $v(n)=s(w_n)$. Then for each vertex $u$ and for each arrow $\alpha$ in $Q$
$$\varphi_C(u)(z_i) = \left\{ \begin{array}{c@{\quad,\quad}l}
z_i & \mbox{if $v(i)=u$}\\ 0 & \mbox{else} \end{array} \right\}\; \mbox{ and } \;
\varphi_C(\alpha)(z_i) = \left\{ \begin{array}{c@{\quad,\quad}l}
z_{i-1} & \mbox{if $w_i=\alpha$}\\ z_{i+1} & \mbox{if $w_{i+1}=\alpha^{-1}$}\\
0 & \mbox{else}
\end{array} \right\} .$$
We will call $\varphi_C$ the canonical representation and $\{z_0,z_1,\ldots,z_n\}$ the canonical $k$-basis for $M(C)$ relative to the representative $C$. Note that $M(C)\cong M(C^{-1})$. 

The string modules $M(1_u)$, $u$ a vertex of $Q$, correspond bijectively to the isomorphism classes of the simple $\Lambda$-modules. We say a string $C=w_1\cdots w_n$ is directed if all $w_i$ are arrows. For each vertex $u$ of $Q$, there exist at most two directed strings of maximal length starting in $u$. Let these be $C_1$ and $C_2$. Then the projective indecomposable $\Lambda$-module $P(u)$ is the string module $M(C_1C_2^{-1})$. Dually, the injective indecomposable module $E(u)$ is the string module $M(D_1^{-1}D_2)$ where $D_1$ and $D_2$ are the directed strings of maximal length ending in $u$.

\subsubsection{Bands and band modules}
\label{sss:bands}
Let $w=w_1\cdots w_n\in {\mathcal{W}}_r$. Then, for $0\leq
i\leq n-1$, the $i$-th rotation of $w$ is defined to be the word $\rho_i(w)=w_{i+1}\cdots w_n w_1 \cdots w_i$. Let $\sim_r$ be the equivalence relation on ${\mathcal{W}}_r$ such that
$w\sim_r w'$ if and only if $w=\rho_i(w')$ for some $i$ or $w^{-1}=\rho_j(w')$ for some $j$. Then bands are representatives $w\in {\mathcal{W}}_r$ of the equivalence classes under $\sim_r$ with the following property: 
$w=w_1\cdots w_n$, $n\ge 1$, with $w_i\neq w_{i+1}^{-1}$ and $w_n\neq w_1^{-1}$, such that $w$ is not a power of a smaller word, and, for all positive integers $m$, no subword of $w^m$ or its formal inverse belongs to $I$.

Let $B=w_1\cdots w_n$ be a band of length $n$. We may assume that $w_1$ is an arrow, by rotating and possibly inverting. Let $Q_B^c$ be the circular quiver
$$Q_B^c = \vcenter
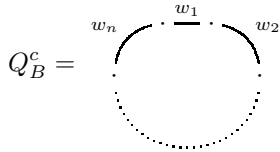
{\xymatrix @-1.2pc {
&\cdot\ar@{-}[r]^{w_1}&\cdot&\\ \cdot \ar@{-}@/^/[ru]^{w_n}\ar@{.}@(d,d)[rrr]&&&\cdot\ar@{-}@/_/[lu]_{w_2}\\ \\}}$$
where $\xymatrix @-1.2pc @C=0pt{\cdot \ar@{-}[rrr]^{w_i}&&&\cdot & = & \cdot&&&&\cdot\ar@/_/[llll]_{\beta}}$ points counter-clockwise if $w_i=\beta$ is an arrow,
and $\xymatrix @-1.2pc @C=0pt{\cdot \ar@{-}[rrr]^{w_i}&&&\cdot & = &\cdot\ar@/^/[rrrr]^{\beta}&&&&\cdot}$ points clockwise if $w_i=\beta^{-1}$ is a formal inverse. 
Let $m>0$ be an integer and $\lambda\in k^*$. Then the representation of $Q_B^c$ which assigns to each vertex the vector space $k^m$, to $w_1$ the indecomposable Jordan matrix $J_m(\lambda)$, and to $w_i$, $2\le i\le n$, the identity map defines an indecomposable $\Lambda$-module, called the band module $M(B,\lambda,m)$ corresponding to $B$, $\lambda$ and $m$. Note that for all $i,j$
$$M(B,\lambda,m)\cong M(\rho_i(B),\lambda,m)\cong M(\rho_j(B)^{-1},\lambda,m).$$


\subsection{Auslander-Reiten components}
\label{ss:arcomps}
\setcounter{equation}{0}
\setcounter{figure}{0}

Let $\Lambda=kQ/I$ be a basic string algebra. Then in each component of the Auslander-Reiten quiver of $\Lambda$ there are either only string modules or only band modules. The band modules all lie in $1$-tubes. The string modules can lie in periodic components or in non-periodic components. 
We now describe the irreducible morphisms between string modules using hooks and cohooks.  Let $S$ be a string.

We say that $S$ starts on a peak provided there is no arrow $\beta$ with $S\beta$ a string, and that $S$ starts in a deep provided there is no arrow $\gamma$ with $S\gamma^{-1}$ a string. Dually, we say that $S$ ends on a peak provided there is no arrow $\beta$ with $\beta^{-1}S$ a string, and that $S$ ends in a deep provided there is no arrow $\gamma$ with $\gamma S$ a string.

If $S$ does not start on a peak, there is a unique arrow $\beta$ and a unique maximal directed string $M$ such that $S_h=S\beta M^{-1}$ starts in a deep. We say $S_h$ is obtained from $S$ by adding a hook on the right side. Dually, if $S$ does not end on a peak, there is a unique arrow $\beta$ and a unique directed string $M$ such that ${}_hS=M\beta^{-1}S$ ends in a deep. We say ${}_hS$ is obtained from $S$ by adding a hook on the left side.

If $S$ does not start in a deep, there is a unique arrow $\gamma$ and a unique maximal directed string $N$ such that $S_c=S\gamma^{-1}N$ starts on a peak. We say $S_c$ is obtained from $S$ by adding a cohook on the right side. Dually, if $S$ does not end in a deep, there is a unique arrow $\gamma$ and a unique directed string $N$ such that ${}_cS=N^{-1}\gamma S$ ends on a peak. We say ${}_cS$ is obtained from $S$ by adding a cohook on the left side.

All irreducible morphisms between string modules are either canonical injections
$$M(S)\to M(S_h),\quad \mbox{or} \quad M(S)\to M({}_hS),$$
or canonical projections
$$M(S_c)\to M(S),\quad \mbox{or} \quad M({}_cS)\to M(S).$$

If $\Lambda$ is a basic special biserial algebra which is self-injective, then $\Lambda/\mathrm{soc}(\Lambda)$ is a string algebra. Moreover, the stable Auslander-Reiten quiver of $\Lambda$ is equal to the Auslander-Reiten quiver of $\Lambda/\mathrm{soc}(\Lambda)$. In case $\Lambda$ is Morita equivalent to a block with dihedral defect groups, the periodic components containing string modules are either $1$-tubes or $3$-tubes, and the non-periodic components all have type $\mathbb{Z}A_\infty^\infty$.


\end{document}